\documentclass{article}

\usepackage{microtype}
\usepackage{graphicx}
\usepackage{subfigure}
\usepackage{booktabs} 

\usepackage{hyperref}

\usepackage{algorithm}
\usepackage{algorithmic}

\usepackage{hyperref}
\usepackage{url}

\usepackage{pifont}

\usepackage{epsfig}
\usepackage{amssymb}
\usepackage{amsmath}
\usepackage{amsthm}
\usepackage{amsfonts}
\usepackage{bbding}
\usepackage{array}
\usepackage{paralist}
\usepackage{xargs}                      
\usepackage{caption}

\newtheorem{thm}{Theorem}
\newtheorem{lem}{Lemma}

\newtheorem{cor}{Corollary}

\newtheorem{ass}{Assumption}
\newtheorem{remark}{Remark}

\newcolumntype{C}[1]{>{\centering\let\newline\\\arraybackslash\hspace{0pt}}m{#1}}

\newcommand\tagthis{\addtocounter{equation}{1}\tag{\theequation}}

\DeclareMathOperator{\Ocal}{\mathcal{O}}

\newcommand{\Exp}[1]{\mathbb{E}\left[#1\right]}           
\newcommand{\E}{\mathbb{E}} 
\newcommand{\set}[1]{\left\{#1\right\}}
\newcommand{\sets}[1]{\{#1\}}

\newcommand{\norms}[1]{\Vert#1\Vert}
\newcommand{\dom}[1]{\mathrm{dom}\left(#1\right)}
\newcommand{\iprods}[1]{\langle#1\rangle}

\usepackage{xcolor}

\newcommand{\R}{\mathbb{R}}

\newcommand{\Eproof}{\hfill $\square$}

\newcommand{\zero}[1]{{\boldsymbol{0}}}

\usepackage[accepted]{icml2021}

\icmltitlerunning{SMG: A Shuffling Gradient-Based Method with Momentum}

\begin{document}

\twocolumn[
\icmltitle{SMG: A Shuffling Gradient-Based Method with Momentum}



\icmlsetsymbol{equal}{*}

\begin{icmlauthorlist}
\icmlauthor{Trang H. Tran}{to1}
\icmlauthor{Lam M. Nguyen}{to2}
\icmlauthor{Quoc Tran-Dinh}{to3}
\end{icmlauthorlist}

\icmlaffiliation{to1}{School of Operations Research and Information Engineering, Cornell University, Ithaca, NY, USA.}
\icmlaffiliation{to2}{IBM Research, Thomas J. Watson Research Center, Yorktown Heights, NY, USA.}
\icmlaffiliation{to3}{Department of Statistics and Operations Research, The University of North Carolina at Chapel Hill, NC, USA}

\icmlcorrespondingauthor{Lam M. Nguyen}{LamNguyen.MLTD@ibm.com}

\icmlkeywords{Machine Learning, ICML}

\vskip 0.3in
]



\printAffiliationsAndNotice{}  

\begin{abstract}
We combine two advanced ideas widely used in optimization for machine learning: \textit{shuffling} strategy and \textit{momentum} technique to develop a novel shuffling gradient-based method with momentum, coined \textbf{S}huffling \textbf{M}omentum \textbf{G}radient (SMG),  for non-convex finite-sum optimization problems. 
While our method is inspired by momentum techniques, its update is fundamentally different from existing momentum-based methods.
We establish state-of-the-art convergence rates of SMG for any shuffling strategy using either constant or diminishing learning rate under standard assumptions (i.e. \textit{$L$-smoothness} and \textit{bounded variance}).
When the shuffling strategy is fixed, we develop another new algorithm that is similar to existing momentum methods,
and prove the same convergence rates for this algorithm under the $L$-smoothness and bounded gradient assumptions. 
We demonstrate our algorithms via numerical simulations on standard datasets and compare them with existing shuffling methods.
Our tests have shown encouraging performance of the new algorithms.
\end{abstract}

\section{Introduction}\label{sec_intro}
Most training tasks in supervised learning are boiled down to solving the following finite-sum minimization:
\begin{equation}\label{ERM_problem_01}
    \min_{w \in \mathbb{R}^d} \Big\{ F(w) := \frac{1}{n} \sum_{i=1}^n f(w ; i)  \Big\},
\end{equation}
where $f(\cdot; i) : \R^d\to\R$ is a given smooth and possibly nonconvex function for $i \in [n] := \set{1,\cdots,n}$.

Problem \eqref{ERM_problem_01} looks simple, but covers various convex and nonconvex applications in machine learning and statistical learning, including, but not limited to, logistic regression, multi-kernel learning, conditional random fields, and neural networks. 
Especially, \eqref{ERM_problem_01} covers the \textit{empirical risk minimization} as a special case. 
Solution methods for approximately solving \eqref{ERM_problem_01} have been widely studied in the literature under different sets of assumptions.
The most common approach is perhaps stochastic gradient-type (SGD) methods \citep{RM1951,ghadimi2013stochastic,Bottou2018,Nguyen2018_sgdhogwild,Nguyen2019_sgd_new_aspects} and their variants.

\textbf{Motivation.}
While SGD and its variants rely on randomized sampling strategies with replacement, gradient-based methods using without-replacement strategies are often easier and faster to implement.
Moreover, practical evidence \citep{bottou2009curiously} has shown that they usually produce a faster decrease of the training loss.
Randomized shuffling strategies (also viewed as sampling without replacement) allow the algorithm to use exactly one function component $f(\cdot; i)$ at each epoch compared to SGD, which has only statistical convergence guarantees (e.g., in expectation or with high probability).
However, very often, the analysis of shuffling methods is more challenging than SGD due to the lack of statistical independence.

In the deterministic case, single permutation (also called shuffle once, or single shuffling) and incremental gradient methods can be considered as special cases of the shuffling gradient-based methods we study in this paper. 
One special shuffling strategy is randomized reshuffling, which is broadly used in practice, where we use a different random permutation at each epoch.
%
Alternatively, in recent years, it has been shown that many gradient-based methods with momentum update can notably boost the convergence speed both in theory and practice \citep{Nesterov2004,dozat2016incorporating,wang2020scheduled}.
These methods have been widely used in both convex and nonconvex settings, especially, in deep learning community.
Remarkably, Nesterov's accelerated method \citep{Nesterov1983} has made a revolution in large-scale convex optimization in the last two decades, and has been largely exploited in nonconvex problems.
The developments we have discussed here motivate us to raise the following research question:
\begin{center}
\textit{Can we combine both shuffling strategy and momentum scheme to develop new provable gradient-based algorithms for handling \eqref{ERM_problem_01}?}
\end{center}


{In this paper, we answer this question affirmatively by proposing a novel algorithm called Shuffling Momentum Gradient (SMG).} We establish its convergence guarantees 
for different shuffling strategies, and in particular, randomized reshuffling strategy. 
We also investigate different variants of our method.

\textbf{Our contribution.}
To this end, our contributions in this paper  can be summarized as follows.
\begin{compactitem}
    \item[(a)] 
    
    We develop {a novel shuffling gradient-based method with momentum (Algorithm~\ref{sgd_momentum_shuffling_01} in Section~\ref{sec:shuffling_momentum})} for approximating a stationary point of the nonconvex minimization problem \eqref{ERM_problem_01}. 
    Our algorithm covers any shuffling strategy ranging from deterministic to randomized, including  incremental, single shuffling, and randomized reshuffling variants.
    
    \item[(b)] We establish the convergence of our method in the nonconvex setting and achieve the state-of-the-art $\Ocal\left(1/{T^{2/3}}  \right)$ 
    convergence rate under standard assumptions (i.e. the $L$-smoothness and bounded variance conditions), 
    where $T$ is the number of epochs. 
    For randomized reshuffling strategy, we can improve our convergence rate up to $\Ocal\left(1/({n^{1/3}T^{2/3})} \right)$. 
    \item[(c)]
    We study different strategies for selecting learning rates (LR), including constant, diminishing, exponential, and cosine scheduled learning rates.
    In all cases, we prove the same convergence rate of the corresponding variants without any additional assumption.
    
    \item[(d)] 
    When a single shuffling strategy is used, we show that a momentum strategy can be incorporated directly at each iteration of the shuffling gradient method to obtain a different variant as presented in      Algorithm~\ref{sgd_momentum_shuffling2}.
    We analyze the convergence of this algorithm and achieve the same $\Ocal(1/T^{2/3})$ epoch-wise convergence rate, but under a bounded gradient assumption instead of the bounded variance as {for the SMG algorithm.}
    
\end{compactitem}


Our $\Ocal(1/T^{2/3})$ convergence rate is the best known so far for shuffling gradient-type methods in nonconvex optimization \citep{nguyen2020unified,mishchenko2020random}.
However, like \citep{mishchenko2020random}, our {SMG method} only requires a generalized bounded variance assumption (Assumption~\ref{as:A1}(c)), which is weaker and more standard than the bounded component gradient assumption used in existing works.
Algorithm~\ref{sgd_momentum_shuffling2} uses the same set of assumptions as in \citep{nguyen2020unified} to achieve the same rate, but has a momentum update.
For the randomized reshuffling strategy, our $\Ocal\left(1/({n^{1/3}T^{2/3})} \right)$ convergence rate also matches the rate of the without-momentum algorithm in  \citep{mishchenko2020random}.
It leads to the total of iterations $nT = \Ocal(\sqrt{n}\varepsilon^{-3})$.

We emphasize that, in many existing momentum variants, the momentum $m_i^{(t)}$ is updated recursively at each iteration as $m_{i+1}^{(t)} := \beta m_i^{(t)} + (1-\beta) g_i^{(t)}$ for a given weight $\beta \in (0, 1)$.
This update shows that the momentum $m_{i+1}^{(t)}$ incorporates all the past gradient terms $g_i^{(t)}, g_{i-1}^{(t)}, g_{i-1}^{(t)}, \dots, g_0^{(t)}$ with exponential decay weights $1, \beta, \beta^2, \dots, \beta^i$, respectively.
However, in shuffling methods, 
the convergence guarantee is often obtained in epoch.
Based on this observation, we modify the classical momentum update in the shuffling method as shown in Algorithm~\ref{sgd_momentum_shuffling_01}.
More specifically, the momentum term $m_0^{(t)}$ is fixed at the beginning of each epoch, and an auxiliary sequence $\{v_i^{(t)}\}$ is introduced to keep track of the gradient average to update the momentum term in the next epoch. 
This modification makes  Algorithm~\ref{sgd_momentum_shuffling_01} fundamentally different from existing momentum-based methods.
This new algorithm still achieves $\Ocal(1/T^{2/3})$ epoch-wise convergence rate under standard assumptions.
To the best of our knowledge, our work is the first analyzing convergence rate guarantees of shuffling-type gradient methods with momentum under standard assumptions.

%

Besides Algorithm~\ref{sgd_momentum_shuffling_01}, we also exploit recent advanced strategies for selecting learning rates, including exponential and cosine scheduled learning rates.
These two strategies have shown state-of-the-art performance in practice \citep{smith2017cyclical,loshchilov10sgdr,li2020exponential}. 
Therefore, it is worth incorporating them in shuffling methods.

\textbf{Related work.}
Let us briefly review the most related works to our methods studied in this paper.

\textbf{\textit{Shuffling gradient-based methods.}}
Shuffling gradient-type methods for solving \eqref{ERM_problem_01} have been widely studied in the literature in recent years \citep{bottou2009curiously,Gurbuzbalaban2019,shamir2016without,haochen2019random,nguyen2020unified} for both convex and nonconvex settings.
It was empirically investigated in a short note \citep{bottou2009curiously} and also discussed in \citep{bottou2012stochastic}. 
These methods have also been implemented in several software packages such as TensorFlow and PyTorch, broadly used in machine learning \citep{tensorflow2015-whitepaper,pytorch}. 

In the strongly convex case, shuffling methods have been extensively studied in  \citep{ahn2020sgd,Gurbuzbalaban2019,haochen2019random,safran2020good,nagaraj2019sgd,rajput2020closing,nguyen2020unified,mishchenko2020random} under different assumptions.
The best known convergence rate in this case is  $\Ocal(1/(nT)^2 + 1/(nT^3))$, which matches the lower bound rate studied in  \citep{safran2020good} up to some constant factor.
Most results in the convex case are for the incremental gradient variant, which are studied in \citep{nedic2001incremental,nedic2001convergence}.
Convergence results of shuffling methods on the general convex case are investigated in \citep{shamir2016without,mishchenko2020random}, where \citep{mishchenko2020random} provides a unified approach to cover different settings.
The authors in \citep{ying2017convergence} combine a randomized shuffling and a variance reduction technique (e.g., SAGA \citep{SAGA} and SVRG \citep{SVRG}) to develop a new variant. They show a linear convergence rate for strongly convex problems but using an energy function, which is unclear how to convert it into known convergence criteria.


In the nonconvex case, \citep{nguyen2020unified} first shows $\Ocal(1/T^{2/3})$ convergence rate for a general class of shuffling gradient methods under the $L$-smoothness and bounded gradient assumptions on \eqref{ERM_problem_01}.
This analysis is then extended in \citep{mishchenko2020random} to a more relaxed assumption.
The authors in \citep{meng2019convergence} study different distributed SGD variants with shuffling for  strongly convex, general convex, and nonconvex problems.
An incremental gradient method for weakly convex problems is investigated in \citep{li2019incremental}, where the authors show $\Ocal(1/T^{1/2})$ convergence rate as in standard SGD.
To the best of our knowledge, the best known rate of shuffling gradient methods for the nonconvex case under standard assumptions is  $\Ocal(1/T^{2/3})$ as shown in \citep{nguyen2020unified,mishchenko2020random}.

Our Algorithm~\ref{sgd_momentum_shuffling_01} developed in this paper is a nontrivial momentum variant of the general shuffling method in \citep{nguyen2020unified} but our analysis uses a standard bounded variance assumption instead of bounded gradient one.

\textbf{\textit{Momentum-based methods.}}
Gradient methods with momentum were studied in early works for convex problems such as heavy-ball, inertial, and Nesterov's accelerated gradient methods \citep{polyak1964some,Nesterov2004}.
Nesterov's accelerated method is the most influent scheme and achieves optimal convergence rate for convex problems.
While momentum-based methods are not yet known to improve theoretical convergence rates in the nonconvex setting, they show significantly encouraging performance in practice \citep{dozat2016incorporating,wang2020scheduled}, especially in the deep learning community.
However, the momentum strategy has not yet been exploited in shuffling methods.

\textbf{\textit{Adaptive learning rate schemes.}}
Gradient-type methods with adaptive learning rates such as AdaGrad \citep{AdaGrad} and Adam \citep{Kingma2014} have shown state-of-the-art performance in several optimization applications.
Recently, many adaptive schemes for learning rates have been proposed such as diminishing \citep{nguyen2020unified}, exponential scheduled \citep{li2020exponential}, and cosine scheduled 
\citep{smith2017cyclical,loshchilov10sgdr}. 
In \citep{li2020exponential}, the authors analyze convergence guarantees for the exponential and cosine learning rates in SGD.
These adaptive learning rates have also been empirically studied in the literature, especially in machine learning tasks.
Their convergence guarantees have also been investigated accordingly under certain assumptions. 

Although the analyses for adaptive learning rates are generally non-trivial, our theoretical results in this paper are flexible enough to cover various different learning rates. 
However, we only exploit the diminishing, exponential scheduled, and cosine scheduled schemes for our shuffling methods with momentum.
We establish that in the last two cases, our algorithm still achieves state-of-the-art $\Ocal(T^{-2/3})$ (and possibly up to $\Ocal(n^{-1/3}T^{-2/3})$) epoch-wise rates.

\textbf{Content.}
The rest of this paper is organized as follows.
Section~\ref{sec:shuffling_momentum} describes our novel {method, Shuffling Momentum Gradient} (Algorithm~\ref{sgd_momentum_shuffling_01}).
Section~\ref{sec:convergence_analysis} investigates its convergence rate under different shuffling-type strategies and different learning rates.
Section~\ref{sec:single_shuffling} proposes an algorithm {with traditional momentum update} for single shuffling strategy.
Section~\ref{sec:experiments} presents our numerical simulations.
Due to space limit, the convergence analysis of our methods, all technical proofs, and additional experiments are deferred to the Supplementary Document (Supp. Doc.).
\section{Shuffling Gradient-Based Method with Momentum}\label{sec:shuffling_momentum}
In this section, we describe our new shuffling gradient algorithm with momentum in Algorithm~\ref{sgd_momentum_shuffling_01}. 
We also compare our algorithm with existing methods and discuss its per-iteration complexity and possible modifications, e.g., mini-batch.

\begin{algorithm}[hpt!]
  \caption{Shuffling Momentum Gradient (SMG)}
  \label{sgd_momentum_shuffling_01}
\begin{algorithmic}[1]
  \STATE {\bfseries Initialization:} 
  Choose $\tilde{w}_0 \in\R^d$ and set $\tilde{m}_0 := \textbf{0}$.{\!\!}
  \FOR{$t := 1,2,\cdots,T $}
  \STATE Set $w_0^{(t)} := \tilde{w}_{t-1}$; $m_0^{(t)} := \tilde{m}_{t-1}$; and $v_0^{(t)} := \textbf{0}$; \label{alg:A1_step1}
  \STATE Generate an arbitrarily deterministic or random permutation  $\pi^{(t)}$ of $[n]$; \label{alg:A1_step2}
  \FOR{$i := 0,\cdots,n-1$}
   \STATE Query $ g_{i}^{(t)} := \nabla f ( w_{i}^{(t)} ; \pi^{(t)} ( i + 1 ) )$; \label{alg:A1_step3}
    \STATE Choose $\eta^{(t)}_i := \frac{\eta_t}{n}$ and update \label{alg:A1_step4}
    \begin{equation}\label{eq:main_update}
    \left\{\begin{array}{lcl}
        m_{i+1}^{(t)} & := & \beta m_{0}^{(t)} + (1-\beta) g_{i}^{(t)} \vspace{1ex}\\
        v_{i+1}^{(t)} & := &  v_{i}^{(t)} + \frac{1}{n} g_{i}^{(t)} \vspace{1ex}\\
        w_{i+1}^{(t)} & := & w_{i}^{(t)} - \eta_i^{(t)} m_{i+1}^{(t)};
    \end{array}\right.
    \end{equation}
  \ENDFOR
  \STATE Set $\tilde{w}_t := w_{n}^{(t)}$ and $\tilde{m}_t := v_{n}^{(t)}$; \label{alg:A1_step5}
  \ENDFOR
  \STATE\textbf{Output:} 
Choose $\hat{w}_T \in \{ \tilde{w}_{0},\cdots,\tilde{w}_{T-1}\}$ at random with probability $\mathbb{P}[\hat{w}_T = \tilde{w}_{t-1}] = \frac{\eta_t}{\sum_{t=1}^{T} \eta_t}$.
\label{alg:A1_step6}
\end{algorithmic}
\end{algorithm}

Clearly, if $\beta = 0$, then Algorithm~\ref{sgd_momentum_shuffling_01} reduces to the standard shuffling gradient method as in \citep{nguyen2020unified,mishchenko2020random}.
{Since each inner iteration of our method uses one component $f(\cdot,i)$, we use $\eta_i^{(t)} = \frac{\eta_t}{n}$ in Algorithm~\ref{sgd_momentum_shuffling_01} to make it consistent with one full gradient, which consists of $n$ gradient components.
This form looks different from the learning rates used in previous work for SGD \cite{shamir2016without,mishchenko2020random}, however, it does not necessary make our learning rate smaller. 
In the same order of training samples $n$, our learning rate matches the one in \cite{mishchenko2020random} as well as the state-of-the-art complexity results. 
In fact, the detailed learning rate $\eta_t$ used in Algorithm~\ref{sgd_momentum_shuffling_01} will be specified in Section~\ref{sec:convergence_analysis}. }

\textbf{\textit{Comparison.}}
Unlike existing momentum methods where $m_i^{(t)}$ in \eqref{eq:main_update} is updated recursively as $m_{i+1}^{(t)} := \beta m_i^{(t)} + (1-\beta) g_i^{(t)}$ for $\beta \in (0, 1)$, we instead fix the first term $m_0^{(i)}$ in \eqref{eq:main_update} at each epoch.
It is only updated at the end of each epoch by averaging all the gradient components $\{ g_{i}^{(t)}\}_{i=0}^{n-1}$ evaluated in such an epoch.
To avoid storing these gradients, we introduce an auxiliary variable $v_i^{(t)}$ to keep track of the gradient average. 
Here, we fix the weight $\beta$ for the sake of our analysis, but it is possible to make $\beta$ adaptive.

Our new method is based on the following observations.  First, while SGD generally uses an unbiased estimator for the full gradient, shuffling gradient methods often do not have such a nice property, making them difficult to analyze. Due to this fact, updating the momentum at each inner iteration does not seem preferable since it could make the estimator deviate from the true gradient. Therefore, we consider updating the momentum after each epoch. Second, unlike the traditional momentum with exponential decay weights,
our momentum term $m_0^{(t)}$ is an equal-weighted average of all the past gradients in an epoch, but evaluated at different points $w_i^{(t-1)}$ in the inner loop.
Based on these observations, the momentum term should aid the full gradient $\nabla{F}$ instead of its component $\nabla{f}(\cdot; i)$, leading to the update at the end of each epoch.  

It is also worth noting that our algorithm is fundamentally different from variance reduction methods. 
For instance, SVRG and SARAH variants require evaluating full gradient at each snapshot point while our method always uses single component gradient. 
Hence, our algorithm does not require a full gradient evaluation at each outer iteration, and our momentum term does not require full gradient of $f$.

When the learning rate $\eta_i^{(t)}$ is fixed at each epoch as $\eta_i^{(t)} := \frac{\eta_t}{n}$, where $\eta_t > 0$, we can derive from \eqref{eq:main_update} that
\begin{equation*}
w_{i+1}^{(t)} = w_i^{(t)} - \frac{(1-\beta)\eta_t}{n}g_i^{(t)} + \frac{\beta\eta_t}{n(1-\beta)\eta_{t-1}}e_t,
\end{equation*}
where $e_t := \tilde{w}_{t-1} - \tilde{w}_{t-2} + \beta\eta_{t-1}\tilde{m}_{t-2}$.
Here, $e_t$ plays a role as a momentum or an inertial term, but it is different from the usual momentum term. 
However, we still name Algorithm~\ref{sgd_momentum_shuffling_01} the \textit{Shuffling Momentum Gradient} since it is inspired by momentum methods.

\textbf{\textit{Per-iteration complexity.}}
The per-iteration complexity of Algorithm~\ref{sgd_momentum_shuffling_01} is almost the same as in standard shuffling gradient schemes, see, e.g., \citep{shamir2016without}.
It needs to store two additional vectors $m_{i+1}^{(t)}$ and $v_i^{(t)}$, and performs two vector additions and three scalar-vector multiplications. 
Such additional costs are very mild.

\textbf{\textit{Shuffling strategies.}}
Our convergence guarantee in Section~\ref{sec:convergence_analysis} holds for any permutation $\pi^{(t)}$ of $\{1,2,\cdots, n\}$, including deterministic and randomized ones. 
Therefore, our method covers any shuffling strategy, including incremental, single shuffling, and randomized reshuffling variants as special cases. 
We highlight that our convergence result for the randomized reshuffling variant is significantly better than the general ones as we will explain in detail in Section~\ref{sec:convergence_analysis}.

\textbf{\textit{Mini-batch.}}
Algorithm~\ref{sgd_momentum_shuffling_01} also works with mini-batches, and our theory can be slightly adapted to establish the same convergence rate for mini-batch variants.
However, it remains unclear to us how mini-batches can improve the convergence rate guarantee of Algorithm~\ref{sgd_momentum_shuffling_01} in the general case where the shuffling strategy is not specified.



\section{Convergence Analysis}\label{sec:convergence_analysis}
We analyze the convergence of Algorithm~\ref{sgd_momentum_shuffling_01} under standard assumptions, which is organized as follows.

\subsection{Technical Assumptions}\label{subsec:assumptions}
Our analysis relies on the following assumptions:
\begin{ass}\label{as:A1}
Problem \eqref{ERM_problem_01} satisfies:
\begin{compactitem}
\item[$\mathrm{(a)}$]\textbf{(Boundedness from below)} 
$\mathrm{dom}(F) \neq \emptyset$ and $F$ is bounded from below, i.e. $F_{*} := \displaystyle\inf_{w\in\R^d}F(w) > -\infty$.

\item[$\mathrm{(b)}$]\textbf{($L$-Smoothness)} 
$f(\cdot; i)$ is $L$-smooth for all $i \in [n]$, i.e. there exists a universal constant $L > 0$ such that, for all $w,w' \in \dom{F}$, it holds that
\begin{equation}\label{eq:Lsmooth_basic}
\norms{ \nabla f(w;i) - \nabla f(w';i)} \leq L \norms{ w - w'}. 
\end{equation}

\item[$\mathrm{(c)}$]\textbf{(Generalized bounded variance)} 
There exist two non-negative and finite constants $\Theta$ and $\sigma$ such that for any $w \in \dom{F}$ we have
\begin{equation}\label{eq:bounded_variance}
\hspace{-3ex}   
 \frac{1}{n} \sum_{i=1}^n \norms{\nabla f(w; i) - \nabla F(w)}^2 \leq \Theta \| \nabla F(w) \|^2 + \sigma^2. 
 \hspace{-1ex}
\end{equation}
\end{compactitem}
\end{ass}
Assumption~\ref{as:A1}(a) is required in any algorithm to guarantee the well-definedness of \eqref{ERM_problem_01}.
The $L$-smoothness \eqref{eq:Lsmooth_basic} is standard in gradient-type methods for both stochastic and deterministic algorithms.
{From this assumption, we have for any $w, w' \in \dom{F}$ (see  \citep{Nesterov2004}):}
\begin{equation}\label{eq:Lsmooth}
F(w) \leq F(w') + \iprods{\nabla F(w'), w-w'} + \frac{L}{2}\norms{w-w'}^2. 
\end{equation}
The condition \eqref{eq:bounded_variance}  in Assumption~\ref{as:A1}(c) reduces to the standard bounded variance condition if $\Theta = 0$. 
Therefore,  \eqref{eq:bounded_variance} is more general than the bounded variance assumption, which is often used in stochastic optimization. 
Unlike recent existing works on  momentum SGD and shuffling  \citep{chen2018convergence,nguyen2020unified}, we do not assume the bounded gradient assumption on each $f(\cdot; i)$ in Algorithm~\ref{sgd_momentum_shuffling_01}  (see Assumption~\ref{as:A2}).
This condition is stronger than \eqref{eq:bounded_variance}.

\subsection{Main Result 1 and Its Consequences}\label{subsec:main_result}
 Our first main result   is the following convergence theorem for Algorithm~\ref{sgd_momentum_shuffling_01} under any shuffling strategy.
 
\begin{thm}\label{thm_momentum_variance}
Suppose that Assumption~\ref{as:A1} holds for \eqref{ERM_problem_01}.
Let $\{w_i^{(t)}\}_{t=1}^{T}$ be generated by Algorithm \ref{sgd_momentum_shuffling_01} with a fixed momentum weight $0\leq \beta < 1$ and an epoch learning rate $\eta_i^{(t)} := \frac{\eta_t}{n}$ for every $t \geq 1$. 
Assume that $\eta_0 = \eta_1$, $\eta_t \geq \eta_{t+1}$, and $0 < \eta_t \leq \frac{1}{L\sqrt{K}}$ for $t \geq 1$, where $K := \max \left\{ \frac{5}{2}, \frac{9(5 -3\beta) (\Theta + 1)}{1-\beta}\right\}$. 
Then, it holds that
\begin{align} \label{eq_thm_momentum_variance}
&
\mathbb{E}\big[\| \nabla F( \hat{w}_T )  \|^2 \big]  =  \displaystyle\frac{1}{\sum_{t=1}^T \eta_t}  \displaystyle\sum_{t=1}^T \eta_t \norms{ \nabla F( \tilde{w}_{t-1} )}^2  \\
& \leq \frac{ 4[F( \tilde{w}_0 ) - F_* ]}{(1-\beta)\sum_{t=1}^T \eta_t} +  \displaystyle\frac{9 \sigma^2 L^2 (5 -3\beta)}{(1-\beta)} \left( \frac{\sum_{t=1}^T \eta_{t-1}^3}{\sum_{t=1}^T \eta_t} \right). \nonumber
\end{align}
\end{thm}


\begin{remark}[Convergence guarantee]
 When a deterministic permutation $\pi^{(t)}$ is used, our convergence rate can be achieved in a deterministic sense.
However, to unify our analysis, we will express our convergence guarantees in expectation, where the expectation is taken over all the randomness generated by $\pi^{(t)}$ and $\hat{w}_T$ up to $T$ iterations. Since we can choose permutations $\pi^{(t)}$ either deterministically or randomly, our  bounds in the sequel will hold either deterministically or with probability $1$ (w.p.1), respectively.
Without loss of generality, we write these results in expectation.
\end{remark}

Next, we derive two direct consequences of Theorem~\ref{thm_momentum_variance} by choosing constant and diminishing learning rates. 

\begin{cor}[\textbf{Constant learning rate}]\label{co:constant_LR}
Let us fix the number of epochs $T \geq 1$, and  choose a constant learning rate $\eta_t := \frac{\gamma}{T^{1/3}}$ for some $\gamma > 0$ such that $\frac{\gamma}{T^{1/3}} \leq \frac{1}{L\sqrt{K}}$ for $t\geq 1$ in Algorithm~\ref{sgd_momentum_shuffling_01}.
Then, under the conditions of Theorem~\ref{thm_momentum_variance}, $\mathbb{E}\big[\| \nabla F( \hat{w}_T )  \|^2 \big] $ is upper bounded by
\begin{align*} 
    \frac{1}{T^{2/3}}  \left(\frac{ 4[F( \tilde{w}_0 ) - F_* ]}{(1-\beta)\gamma} +  \frac{9 \sigma^2  (5 -3\beta)L^2\gamma^2}{(1-\beta)}  \right). 
\end{align*}
Consequently, the convergence rate of Algorithm~\ref{sgd_momentum_shuffling_01} is $\mathcal{O}(T^{-2/3})$ in epoch.
\end{cor}

With a constant LR as in Corollary \ref{co:constant_LR}, the convergence rate of Algorithm~\ref{sgd_momentum_shuffling_01} is exactly expressed as
\begin{equation*}
\Ocal\left(\frac{[F( \tilde{w}_{0} ) - F_*] + \sigma^2 }{T^{2/3}} \right),
\end{equation*}
which matches the best known rate in the literature \citep{mishchenko2020random,nguyen2020unified} in term of $T$ for general shuffling-type strategies.

\begin{cor}[\textbf{Diminishing learning rate}]\label{co:diminishing_LR}
Let us choose a diminishing learning rate $\eta_t := \frac{\gamma}{(t + \lambda)^{1/3}}$ for some $\gamma > 0$ and $\lambda \geq 0$ for $t\geq 1$ such that $\eta_1 := \frac{\gamma}{(1 + \lambda)^{1/3}} \leq \frac{1}{L\sqrt{K}}$ in Algorithm~\ref{sgd_momentum_shuffling_01}.
Then, under the conditions of Theorem~\ref{thm_momentum_variance}, after $T$ epochs with $T \geq 2$, we have
\begin{equation*} 
    \mathbb{E}\big[\| \nabla F( \hat{w}_T )  \|^2 \big]    \leq \frac{C_1 + C_2 \log (T - 1 +\lambda)}{ (T + \lambda)^{2/3} - (1 + \lambda)^{2/3}},
\end{equation*}
where $C_1$ and $C_2$ are respectively given by
\begin{align*}
    C_1 &:= \frac{4\left[F(  \tilde{w}_{0} )- F_*\right]}{(1-\beta)\gamma }  + \frac{18 \sigma^2 L^2 (5 -3\beta)\gamma^2}{(1-\beta)(1+\lambda)} \quad \text{and} \\
    C_2 &:= \frac{9 \sigma^2 L^2 (5 -3\beta)\gamma^2}{(1-\beta)}.
\end{align*}
Consequently, the convergence rate of Algorithm~\ref{sgd_momentum_shuffling_01} is $\mathcal{O}(T^{-2/3}\log(T))$ in epoch.
\end{cor}

The diminishing LR $\eta_t := \frac{\gamma}{(t + \lambda)^{1/3}}$ allows us to use large learning rate values at early epochs compared to the constant case.
However, we loose a $\log(T)$ factor in the second term of our worst-case convergence rate bound.


We also derive the following two consequences of Theorem~\ref{thm_momentum_variance} for exponential and cosine scheduled learning rates. 



\textbf{\textit{Exponential scheduled learning rate.}}
Given an epoch budget $T \geq 1$, and two positive constants $\gamma > 0$ and $\rho > 0$, we consider the following exponential LR, see \citep{li2020exponential}:
\begin{equation}\label{exponential_learning_rate_01} 
    \eta_t :=  \frac{\gamma \alpha^t}{T^{1/3}}, \quad \text{where}\ \alpha := \rho^{1/T} \in (0, 1).
\end{equation}
The following corollary shows the convergence of Algorithm~\ref{sgd_momentum_shuffling_01} using this LR without any additional assumption.

\begin{cor}\label{co:exponential_LR}
Let $\{w_i^{(t)}\}_{t=1}^{T}$ be generated by Algorithm~\ref{sgd_momentum_shuffling_01} with $\eta_i^{(t)} := \frac{\eta_t}{n}$, where $\eta_t$ is in \eqref{exponential_learning_rate_01} such that $0 < \eta_t \leq \frac{1}{L\sqrt{K}}$.
Then, under Assumption~\ref{as:A1}, we have 
\begin{align*} 
  \mathbb{E}\big[\| \nabla F( \hat{w}_T )  \|^2 \big] \leq  \frac{ 4[F( \tilde{w}_0 ) - F_* ]}{(1-\beta)\gamma \rho T^{2/3}} +  \frac{9 \sigma^2 L^2 (5 -3\beta)\gamma^2}{(1-\beta)\rho T^{2/3}}. 
\end{align*}
\end{cor}



\textbf{\textit{Cosine annealing learning rate:}}
%
%
Alternatively, given an epoch budget $T \geq 1$, and a positive constant $\gamma > 0$, we consider the following cosine LR for Algorithm~\ref{sgd_momentum_shuffling_01}:
\begin{equation}\label{cos_learning_rate_01} 
    \eta_t := \frac{\gamma}{T^{1/3}}\left(1+ \cos \frac{t \boldsymbol{\pi}}{T} \right), \quad t = 1,2,\cdots, T. 
\end{equation}
This LR is adopted from \citep{loshchilov10sgdr,smith2017cyclical}.
However, different from these works, we fix our learning rate at each epoch instead of updating it at every iteration as in \citep{loshchilov10sgdr,smith2017cyclical}.

\begin{cor}\label{co:cos_LR}
Let $\{w_i^{(t)}\}_{t=1}^{T}$ be generated by Algorithm~\ref{sgd_momentum_shuffling_01} with $\eta_i^{(t)} := \frac{\eta_t}{n}$, where $\eta_t$ is given by \eqref{cos_learning_rate_01}  such that $0 < \eta_t \leq \frac{1}{L\sqrt{K}}$.
Then, under Assumption~\ref{as:A1}, and for $T\geq 2$, $\mathbb{E}\big[\| \nabla F( \hat{w}_T )  \|^2 \big]$ is upper bounded by
\begin{equation*} 
    \frac{1}{ T^{2/3}} \left( \frac{8[F( \tilde{w}_0 ) - F_* ]}{(1-\beta)\gamma } +  \frac{144 \sigma^2 (5 -3\beta)L^2  \gamma^2}{(1-\beta)} \right). 
\end{equation*}
\end{cor}
The scheduled LRs \eqref{exponential_learning_rate_01} and \eqref{cos_learning_rate_01} still preserve our best known convergence rate $\mathcal{O}(T^{-2/3})$.
Note that the exponential learning rates are available in both TensorFlow and PyTorch, while cosine learning rates are also used in PyTorch.

\subsection{Main Result 2: Randomized Reshuffling Variant}



A variant of Algorithm~\ref{sgd_momentum_shuffling_01} is called a \textbf{randomized reshuffling variant} if at each iteration $t$, the generated permutation $\pi^{(t)} = \big(\pi^{(t)}(1), \cdots, \pi^{(t)}(n)\big)$ is uniformly sampled at random without replacement from $\set{1,\cdots, n}$.
Since the randomized reshuffling strategy is extremely popular in practice, we analyze Algorithm~\ref{sgd_momentum_shuffling_01} under this strategy.

\begin{thm}\label{thm_momentum_variance_RR}
Suppose that Assumption~\ref{as:A1} holds for (1).
Let $\{w_i^{(t)}\}_{t=1}^{T}$ be generated by Algorithm \ref{sgd_momentum_shuffling_01} under a \textbf{randomized reshuffling strategy}, a fixed momentum weight $0\leq \beta < 1$, and an epoch learning rate $\eta_i^{(t)} := \frac{\eta_t}{n}$ for every $t \geq 1$. 
Assume that $\eta_t \geq \eta_{t+1}$ and $0 < \eta_t \leq \frac{1}{L\sqrt{D}}$ for $t \geq 1$, where $D = \max \left(\frac{5}{3}, \frac{6(5 -3\beta) (\Theta +n)}{n(1-\beta)}\right)$ and $\eta_0 = \eta_1$. 
Then
\begin{align} \label{eq_thm_momentum_variance_RR}
    &\hspace{-2ex}\mathbb{E}\big[\| \nabla F( \hat{w}_T )  \|^2 \big]  
    =  \frac{1}{\sum_{t=1}^T \eta_t} \sum_{t=1}^{T} \eta_t \Exp{  \norms{ \nabla F( \tilde{w}_{t-1} )}^2}  \\
    &\leq \frac{4 \left[F(  \tilde{w}_{0} )- F_*\right]}{(1-\beta)\sum_{t=1}^T \eta_t}  +  \frac{6 \sigma^2(5 -3\beta)L^2}{n (1-\beta)} \left( \frac{\sum_{t=1}^{T}\eta_{t-1}^3}{\sum_{t=1}^T \eta_t}\right). \nonumber 
\end{align}

\end{thm}

We can derive the following two consequences.

\begin{cor}[\textbf{Constant learning rate}]\label{co:constant_LR_RR}
Let us fix the number of epochs $T \geq 1$, and  choose a constant learning rate $\eta_t := \frac{\gamma n^{1/3}}{T^{1/3}}$ for some $\gamma > 0$ such that $\frac{\gamma n^{1/3}}{T^{1/3}} \leq \frac{1}{L\sqrt{D}}$ for $t\geq 1$ in Algorithm~\ref{sgd_momentum_shuffling_01}.
Then, under the conditions of Theorem~\ref{thm_momentum_variance_RR}, $\mathbb{E}\big[\| \nabla F( \hat{w}_T )  \|^2 \big] $ is upper bounded by
\begin{equation*} 
    \frac{1}{n^{1/3} T^{2/3}} \left(\frac{4\left[F(  \tilde{w}_{0} )- F_*\right]}{(1-\beta)\gamma}  +  \frac{6 \sigma^2(5 -3\beta)L^2\gamma^2}{(1-\beta)}  \right).
\end{equation*}

Consequently, the convergence rate of Algorithm~\ref{sgd_momentum_shuffling_01} is $\mathcal{O}(n^{-1/3}T^{-2/3})$ in epoch.
\end{cor}
%
With a constant LR as in Corollary \ref{co:constant_LR_RR}, the convergence rate of Algorithm~\ref{sgd_momentum_shuffling_01} is improved to
\begin{equation*}
\Ocal\left(\frac{[F( \tilde{w}_{0} ) - F_*] +\sigma^2 }{n^{1/3}T^{2/3}} \right),
\end{equation*}
which matches the best known rate as in the randomized reshuffling scheme, see, e.g., \citep{mishchenko2020random}.

In this case, the total number of iterations $\mathcal{T}_{\textrm{tol}} := nT$ is $\mathcal{T}_{\textrm{tol}} = \Ocal{(\sqrt{n}\varepsilon^{-3})}$ to obtain $\mathbb{E}\big[\| \nabla F( \hat{w}_T )  \|^2 \big] \leq\varepsilon^2$.
Compared to the complexity $\Ocal(\varepsilon^{-4})$ of SGD, our randomized reshuffling variant is better than SGD if $n \leq \Ocal(\varepsilon^{-2})$.

\begin{cor}[\textbf{Diminishing learning rate}]\label{co:diminishing_LR_RR}
Let us choose a diminishing learning rate $\eta_t := \frac{\gamma n^{1/3}}{(t + \lambda)^{1/3}}$ for some $\gamma > 0$ and $\lambda \geq 0$ for $t\geq 1$ such that $\eta_1 := \frac{\gamma n^{1/3}}{(1 + \lambda)^{1/3}} \leq \frac{1}{L\sqrt{D}}$ in Algorithm~\ref{sgd_momentum_shuffling_01}.
Then, under the conditions of Theorem~\ref{thm_momentum_variance_RR}, after $T$ epochs with $T \geq 2$, we have
\begin{equation*} 
    \mathbb{E}\big[\| \nabla F( \hat{w}_T )  \|^2 \big]    \leq \frac{C_3 + C_4 \log (T - 1 +\lambda)}{n^{1/3} \left[ (T + \lambda)^{2/3} - (1 + \lambda)^{2/3}\right]},
\end{equation*}
where $C_3$ and $C_4$ are respectively given by
\begin{align*}
    C_3 &:=  \frac{4\left[F(  \tilde{w}_{0} )- F_*\right]}{(1-\beta)\gamma }  +  \frac{12 \sigma^2(5 -3\beta)L^2 \gamma^2 }{(1-\beta) (1+\lambda)} \quad \text{and} \\
    C_4 &:= \frac{6 \sigma^2(5 -3\beta)L^2 \gamma^2 }{ (1-\beta)}.
\end{align*}
Consequently, the convergence rate of Algorithm~\ref{sgd_momentum_shuffling_01} is $\mathcal{O}(n^{-1/3}T^{-2/3}\log(T))$ in epoch.
\end{cor}
%
%
%
\begin{remark}\label{re:exp_cos}
Algorithm~\ref{sgd_momentum_shuffling_01} under randomized reshuffling still works with exponential and cosine scheduled learning rates, and our analysis is similar to the one with general shuffling schemes.
However, we omit their analysis here.
\end{remark}

\section{Single Shuffling Variant}\label{sec:single_shuffling}
In this section, we modify the single-shuffling gradient method by directly incorporating a momentum term at each iteration.
We prove for the first time that this variant still achieves state-of-the-art convergence rate guarantee under the smoothness and bounded gradient (i.e. Assumption~\ref{as:A2}) assumptions as in existing shuffling methods.
This variant, though somewhat special, also covers an incremental gradient method with momentum as a special case.

Our new momentum algorithm using single shuffling strategy for solving \eqref{ERM_problem_01} is presented in Algorithm~\ref{sgd_momentum_shuffling2}.

\begin{algorithm}[hpt!]
 \caption{Single Shuffling Momentum Gradient}
 \label{sgd_momentum_shuffling2}
\begin{algorithmic}[1]
  \STATE {\bfseries Initialization:} Choose $\tilde{w}_0 \in \R^d$ and set $\tilde{m}_0 := \textbf{0}$;{\!\!}
  \STATE Generate a permutation $\pi$ of $[n]$; \label{alg:A2_step1}
  \FOR{$t := 1,2,\cdots,T $}
  \STATE Set $w_0^{(t)} := \tilde{w}_{t-1}$ and $m_0^{(t)} := \tilde{m}_{t-1}$; \label{alg:A2_step2}
  \FOR{$i = :0,\cdots,n-1$}
    \STATE Query $g_{i}^{(t)} := \nabla f ( w_{i}^{(t)} ; \pi( i + 1 ) )$; \label{alg:A2_step3}
    \STATE Choose $\eta^{(t)}_i := \frac{\eta_t}{n}$ and update \label{alg:A2_step4}
    \begin{equation*}
    \arraycolsep=0.2em
    \left\{\begin{array}{lcl}
        m_{i+1}^{(t)} &:=& \beta m_{i}^{(t)} + (1-\beta) g_{i}^{(t)} \\
        w_{i+1}^{(t)} &:= & w_{i}^{(t)} - \eta_i^{(t)} m_{i+1}^{(t)};
    \end{array}\right.
    \end{equation*}
  \ENDFOR
  \STATE Set $\tilde{w}_t := w_{n}^{(t)}$ and $\tilde{m}_t := m_{n}^{(t)}$; \label{alg:A2_step5}
  \ENDFOR
  \STATE\textbf{Output:}  
  Choose $\hat{w}_T \in \{ \tilde{w}_{0},\cdots,\tilde{w}_{T-1}\}$ at random with probability $\mathbb{P}[\hat{w}_T = \tilde{w}_{t-1}] = \frac{\eta_t}{\sum_{t=1}^{T} \eta_t}$.
  \label{alg:A2_step6}
\end{algorithmic}
\end{algorithm} 

Besides fixing the permutation $\pi$, Algorithm~\ref{sgd_momentum_shuffling2} is different from Algorithm~\ref{sgd_momentum_shuffling_01} at the point of updating $m_i^{(t)}$.
Here, $m_{i+1}^{(t)}$ is updated from $m_{i}^{(t)}$ instead of $m_0^{(t)}$ as in Algorithm \ref{sgd_momentum_shuffling_01}.
In addition, we update $\tilde{m}_t := m_{n}^{(t)}$ instead of the epoch gradient average.
Note that we can write the main update of Algorithm~\ref{sgd_momentum_shuffling2} as $$w_{i+1}^{(t)} := w_i^{(t)} - \frac{(1-\beta)\eta_t}{n}g_i^{(t)} + \beta(w^{(t)}_i - w_{i-1}^{(t)}),$$ which exactly reduces to existing momentum updates.

\textbf{\textit{Incremental gradient method with momentum.}}
If we choose $\pi := [n]$, then we obtain the well-known incremental gradient variant, but with momentum.
Hence, Algorithm~\ref{sgd_momentum_shuffling2} is still new compared to the standard  incremental gradient algorithm \citep{Bertsekas2011}. 
To prove convergence of Algorithm~\ref{sgd_momentum_shuffling2}, we replace Assumption~\ref{as:A1}(c) by the following:

\begin{ass}[\textbf{Bounded gradient}]\label{as:A2}
There exists $G>0$ such that $\norms{\nabla{f}(x; i)} \leq G$, $\forall x\in\dom{F}$ and $i \in [n]$.
\end{ass}
Now, we state the convergence of Algorithm~\ref{sgd_momentum_shuffling2} in the following theorem as our third main result.

\begin{thm}\label{thm_momentum}
Let $\{w_i^{(t)}\}_{t=1}^T$ be  generated by Algorithm \ref{sgd_momentum_shuffling2} with a LR $\eta_i^{(t)} := \frac{\eta_t}{n}$ and $0 < \eta_t \leq \frac{1}{L}$ for  $t\geq 1$.
Then, under Assumption \ref{as:A1}(a)-(b) and Assumption~\ref{as:A2}, we have
\begin{equation}\label{eq_thm_momentum}
\hspace{-0ex}
\arraycolsep=0.1em
\begin{array}{lcl}
\mathbb{E}\big[\| \nabla F( \hat{w}_T )  \|^2 \big]  & &\leq    \dfrac{\Delta_1 } {\big(\sum_{t=1}^T \eta_t \big)(1 - \beta^n)} \\
&& + {~} L^2 G^2\left(  \dfrac{\sum_{t=1}^T \xi_t^3}{\sum_{t=1}^T \eta_t} \right) + \dfrac{4\beta^n G^2}{1 - \beta^n},
\end{array}
\hspace{-1ex}
\end{equation}
where $\xi_t := \max (\eta_t, \eta_{t-1})$ for $t \geq 2$, $\xi_1 = \eta_1$, and
\begin{align*}
    \Delta_1 & := 2 [F( \tilde{w}_{0} ) - F_* ] + \left( \frac{1}{L} + \eta_1 \right) \norms{ \nabla F( \tilde{w}_{0} ) }^2  \\
    & \qquad + 2 L \eta_1^2 G^2. 
\end{align*}
\end{thm}
Compared to \eqref{eq_thm_momentum_variance_RR} in Theorem~\ref{thm_momentum_variance_RR}, \eqref{eq_thm_momentum} strongly depends on the weight $\beta$ as $ \frac{4\beta^n G^2}{1 - \beta^n}$ appears on the right-hand side of \eqref{eq_thm_momentum}.
Hence, $\beta$ must be chosen in a specific form to obtain desired convergence rate.

Theorem~\ref{thm_momentum} provides a key bound to derive concrete convergence rates in the following two corollaries.

\begin{cor}[\textbf{Constant learning rate}]\label{co:SSM_constant_LR}
In Algorithm~\ref{sgd_momentum_shuffling2}, let us fix $T \geq 1$ and choose the parameters:
\begin{compactitem}
    \item $\beta := \left(\frac{\nu}{T^{2/3}}\right)^{1/n}$ for some constant $\nu \geq 0$ such that $\beta \leq \left(\frac{R-1}{R}\right)^{1/n}$ for some $R \geq 1$, and
    \item $\eta_t := \frac{\gamma}{T^{1/3}}$ for some $\gamma > 0$ such that $\frac{\gamma}{T^{1/3}} \leq \frac{1}{L}$. 
\end{compactitem}
Then, under the conditions of Theorem~\ref{thm_momentum}, we have
\begin{equation*}
\arraycolsep=0.2em
\begin{array}{lcl}
\mathbb{E}\big[\| \nabla F( \hat{w}_T )  \|^2 \big]  & \leq &  \dfrac{D_0}{T^{2/3}} + \dfrac{R \norms{ \nabla F( \tilde{w}_{0} ) }^2}{T} + \dfrac{2 L G^2 R \gamma}{T^{4/3}} ,
\end{array}
\end{equation*}
where
\begin{align*}
    D_0 &:= \frac{2 R}{\gamma} [F( \tilde{w}_{0} ) - F_*] + \frac{R}{\gamma L} \norms{ \nabla F( \tilde{w}_{0} ) }^2  \\
    & \qquad + L^2 G^2 \gamma^2 + 4 \nu G^2 R. 
\end{align*}
Thus the convergence rate of Algorithm~\ref{sgd_momentum_shuffling2} is $\mathcal{O}(T^{-2/3})$. 
\end{cor}

With a constant learning rate as in Corollary \ref{co:SSM_constant_LR}, the convergence rate of Algorithm~\ref{sgd_momentum_shuffling2} is
\begin{equation*}
\Ocal\left(\frac{L[F( \tilde{w}_{0} ) - F_*] + \norms{\nabla{F}(\tilde{w}_0)}^2 + G^2}{T^{2/3}}\right),
\end{equation*}
which depends on $L[F( \tilde{w}_{0} ) - F_*]$, $\norms{\nabla{F}(\tilde{w}_0)}^2$, and $G^2$, and is slightly different from Corollary~\ref{co:constant_LR}.

\begin{cor}[\textbf{Diminishing learning rate}]\label{co:SSM_diminishing_LR}
In Algorithm~\ref{sgd_momentum_shuffling2}, let us choose the parameters:
\begin{compactitem}
    \item $\beta := \left(\frac{\nu}{T^{2/3}}\right)^{1/n}$ for some constant $\nu \geq 0$ such that $\beta \leq \left(\frac{R-1}{R}\right)^{1/n}$ for some $R \geq 1$, and
    \item a diminishing learning rate $\eta_t := \frac{\gamma}{(t + \lambda)^{1/3}}$ for all $t \in [T]$ for some $\gamma > 0$ and $\lambda \geq 0$ such that $\eta_1 = \frac{\gamma}{(1 + \lambda)^{1/3}} \leq \frac{1}{L}$. 
\end{compactitem}
Then, under the conditions of Theorem~\ref{thm_momentum}, we have
\begin{equation*}
\arraycolsep=0.2em
\begin{array}{lcl}
\mathbb{E}\big[\| \nabla F( \hat{w}_T )  \|^2 \big]  & \leq &  \dfrac{D_1}{[(T+\lambda)^{2/3} - (1+\lambda)^{2/3}]} + \dfrac{4\nu G^2 R}{T^{2/3}} \\
&& +  {~} \dfrac{L^2 G^2 \gamma^2 \log(T - 1 + \lambda)}{[(T+\lambda)^{2/3} - (1+\lambda)^{2/3}]}, 
\end{array}
\end{equation*}
for $T \geq 2$, where
\begin{equation*}
\arraycolsep=0.1em
\begin{array}{lcl}
D_1 &:= & \frac{2R}{\gamma} [F( \tilde{w}_{0} ) - F_*] + \big[ \frac{R}{\gamma L} + \frac{R}{(1+\lambda)^{1/3}} \big] \norms{ \nabla F( \tilde{w}_{0} ) }^2 \\
&& +  {~} \frac{2 R L \gamma G^2}{(1+\lambda)^{2/3}} + \frac{2}{1+\lambda} L^2 G^2 \gamma^2.
\end{array}
\end{equation*}
Thus the convergence rate of Algorithm~\ref{sgd_momentum_shuffling2} is $\mathcal{O}(\frac{\log(T)}{T^{2/3}})$.
\end{cor}



\begin{remark}\label{re:adaptive_LR}
Algorithm~\ref{sgd_momentum_shuffling2} still works with exponential and cosine scheduled LRs, and our analysis is similar to the one in Algorithm~\ref{sgd_momentum_shuffling_01}, which is 
deferred to the Supp. Doc. 
\end{remark}

\section{Numerical Experiments}\label{sec:experiments}
In order to examine our algorithms, we present two numerical experiments for different nonconvex problems and compare them with some state-of-the-art SGD-type and shuffling gradient methods. 

\subsection{Models and Datasets}


\textbf{Neural Networks.}
We test our Shuffling Momentum Gradient (SMG) algorithm using two standard network architectures: fully connected network (FCN) and convolutional neural network (CNN). 
For the fully connected setting, we train the classic LeNet-300-100 model \citep{MNIST} on the \texttt{Fashion-MNIST} dataset \citep{xiao2017fashion} with $60,000$ images. 

We also use the convolutional LeNet-5 \citep{MNIST} architecture to train the well-known \texttt{CIFAR-10} dataset \citep{CIFAR10} with $50,000$ samples. 
We repeatedly run the experiments for 10 random seeds and report the average results. 
All the algorithms are implemented and run in Python using the PyTorch package \citep{pytorch}.

\textbf{Nonconvex Logistic Regression.}
{Nonconvex regularizers are widely used in statistical learning such as approximating sparsity or gaining robustness.} We consider the following nonconvex binary classification problem:
\begin{align*}
   \min_{w \in \mathbb{R}^d} \Big \{ F(w) \! := \! \frac{1}{n} \sum_{i=1}^n \log(1 \! + \! \exp(- y_i x_i^\top w )) \!+ \!\lambda r(w)  \Big \}, 
\end{align*}
where $\sets{(x_i, y_i)}_{i=1}^n$ is a set of training samples; and $r(w) := \frac{1}{2} \sum_{j=1}^d \frac{w_j^2}{1 + w_j^2}$ is a nonconvex regularizer, and $\lambda := 0.01$ is a regularization parameter. {This example was also previously used in \citep{wang2019spiderboost,tran2019hybrid,nguyen2020unified}.}

We have conducted the experiments on two classification datasets \texttt{w8a} ($49,749$ samples) and \texttt{ijcnn1} ($91,701$ samples) from \texttt{LIBSVM} \citep{LIBSVM}. 
The experiment is repeated with random seeds 10 times and the average result is reported.

\subsection{Comparing SMG with Other Methods}
We compare our SMG algorithm with Stochastic Gradient Descent (SGD) and two other methods: SGD with Momentum (SGD-M) \citep{polyak1964some} and Adam \citep{Kingma2014}. 
For the latter two algorithms, we use the hyper-parameter settings recommended and widely used in practice (i.e. momentum: 0.9 for SGD-M, and two hyper-parameters $\beta_1 := 0.9$, $\beta_2 := 0.999$ for Adam). 
For our new SMG algorithm, we fixed the parameter $\beta := 0.5$ since it usually performs the best in our experiments.

To have a fair comparison, we apply the randomized reshuffling  scheme to all methods. 
Note that shuffling strategies are broadly used in practice and have been implemented in TensorFlow, PyTorch, and Keras \citep{chollet2015keras}. 
We tune each algorithm using constant learning rate and report the best result obtained. 

\begin{figure}[ht!] 
%
\includegraphics[width=0.238\textwidth]{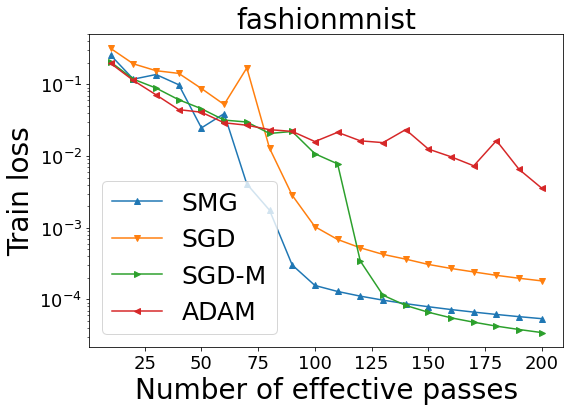}
\includegraphics[width=0.238\textwidth]{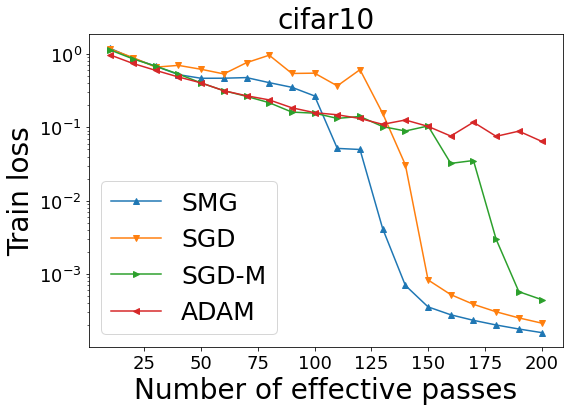}
\caption{The train loss produced by SMG, SGD, SGD-M, and Adam for \texttt{Fashion-MNIST} and \texttt{CIFAR-10}, respectively.}
\label{fig_NN_02}
  
\end{figure}

\textbf{Results.}
Our first experiment is presented in Figure~\ref{fig_NN_02}, where we depict the value of ``train loss'' (i.e. $F(w)$ in \eqref{ERM_problem_01}) on the $y$-axis and the ``number of effective passes'' (i.e. the number of epochs) on the $x$-axis.
It was observed that SGD-M and Adam work  well  for machine learning tasks (see, e.g., \citep{ruder2017overview}). 
Align with this fact, from Figure~\ref{fig_NN_02}, we also observe  that our SMG algorithm and SGD is slightly worse than SGD-M and Adam at the initial stage when training a neural network. 
However, SMG quickly catches up to Adam and demonstrates a good performance at the end of the training process. 

\begin{figure}[ht!] 
%
\includegraphics[width=0.238\textwidth]{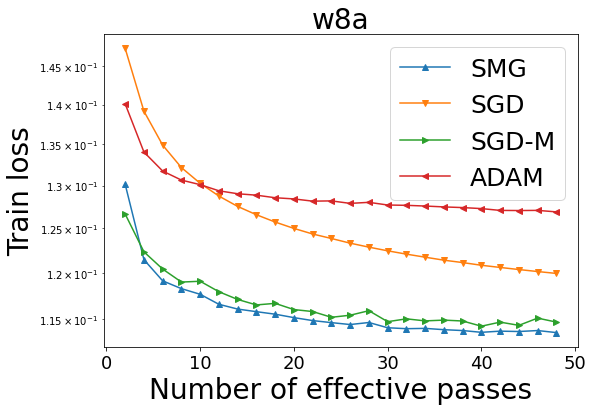}
\includegraphics[width=0.238\textwidth]{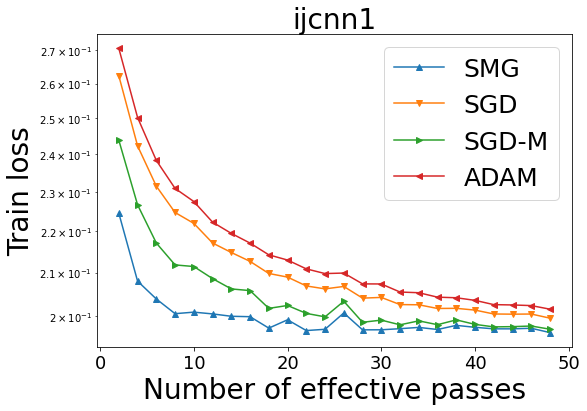}
\caption{The train loss produced by SMG, SGD, SGD-M, and Adam for the \texttt{w8a} and \texttt{ijcnn1} datasets, respectively.}
\label{fig_LR_02}
\end{figure}

For the nonconvex logistic regression problem,  our result is reported in Figure \ref{fig_LR_02}.
For two small datasets tested in our experiments, our algorithm performs significantly better than SGD and Adam, and slightly better than SGD with momentum.
\subsection{The Choice of Hyper-parameter $\beta$}
Since the hyper-parameter $\beta$ plays a critical role in the proposed SMG method, our next experiment is to investigate how this algorithm depends on $\beta$, while using the same constant learning rate. 

\begin{figure}[ht!] 
%
\includegraphics[width=0.238\textwidth]{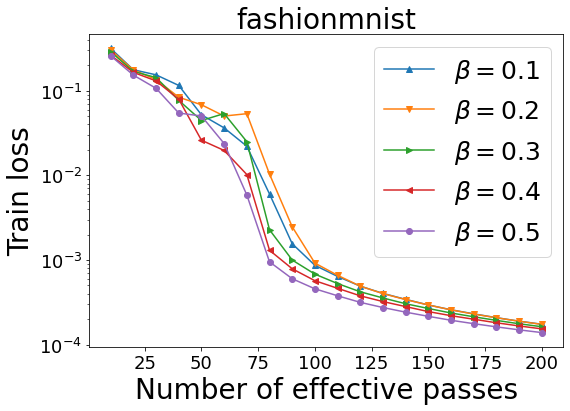}
\includegraphics[width=0.238\textwidth]{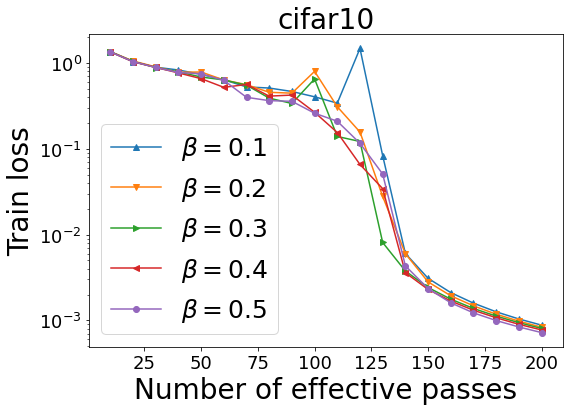}
\caption{The train loss reported by SMG with different $\beta$ on the \texttt{Fashion-MNIST} and \texttt{CIFAR-10} datasets, respectively.}
\label{fig_NN_03}
\end{figure}

\textbf{Results.}
Our result is presented in  Figure~\ref{fig_NN_03}.
We can observe from this figure that in the early stage of the training process, the choice $\beta := 0.5$ gives comparably good performance comparing to other smaller values. 
This choice also results in the best train loss in the end of the training process. 
However, the difference is not really significant, showing that SMG seems robust to the choice of the momentum weight $\beta$ in the range of $[0.1, 0.5]$.

\begin{figure}[ht!] 
%
\includegraphics[width=0.238\textwidth]{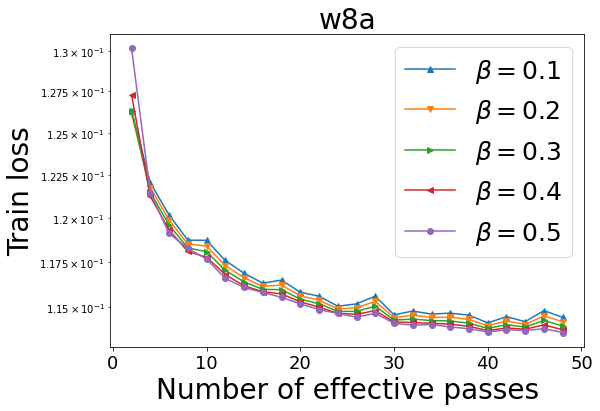}
\includegraphics[width=0.238\textwidth]{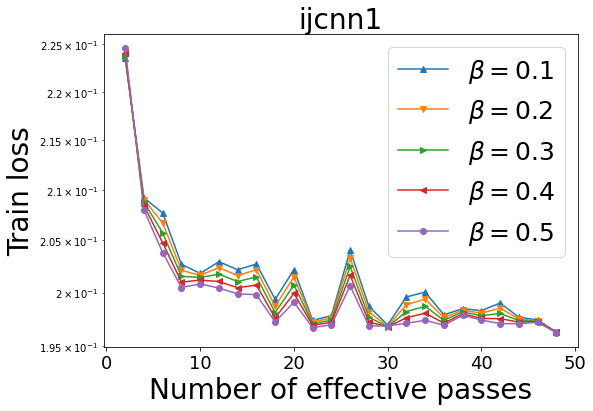}
\caption{The train loss produced by SMG under different values of $\beta$ on the \texttt{w8a} and \texttt{ijcnn1} datasets, respectively.}
\label{fig_LR_03}
\end{figure}

In the nonconvex logistic regression problem, the same choice of $\beta$ also yields similar outcome for two datasets: \texttt{w8a} and \texttt{ijcnn1}, as shown in Figure~\ref{fig_LR_03}. 
We have also experimented with different choices of $\beta \in \{0.6, 0.7, 0.8, 0.9\}$.
However, these choices of $\beta$ do not lead to good performance, and, therefore, we omit to report them here.
This finding raises an open question on the optimal choice of $\beta$.
Our empirical study here shows that the choice of $\beta$ in $[0.1, 0.5]$ works reasonably well, and $\beta := 0.5$ seems to be the best in our test.

\subsection{Different Learning Rate Schemes}
Our last experiment is to examine the effect of different learning rate variants on the performance of our SMG method, i.e. Algorithm~\ref{sgd_momentum_shuffling_01}.

\begin{figure}[ht!] 
%
\includegraphics[width=0.238\textwidth]{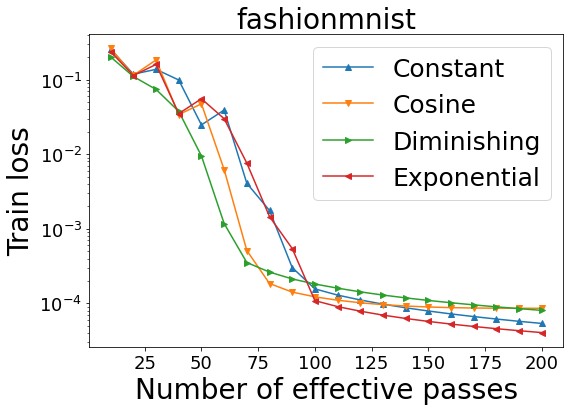}
\includegraphics[width=0.238\textwidth]{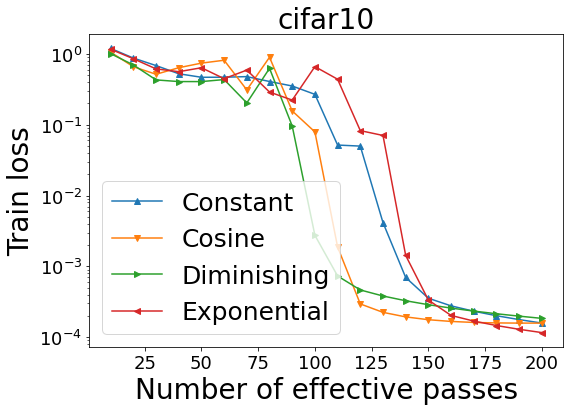}
\caption{The train loss produced by SMG under four different learning rate schemes on the \texttt{Fashion-MNIST} and \texttt{CIFAR-10} datasets, respectively.}
\label{fig_NN_04}
\end{figure}
\begin{figure}[ht!] 
%
\includegraphics[width=0.238\textwidth]{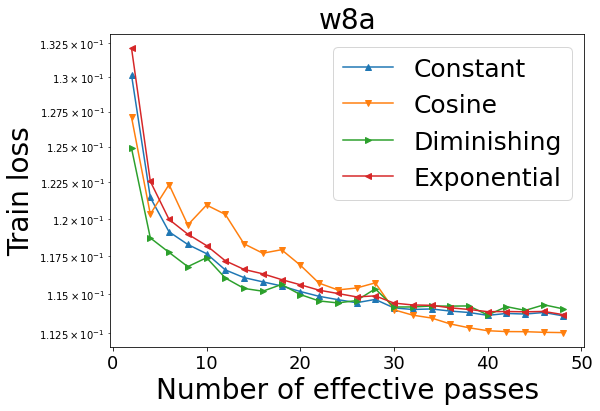}
\includegraphics[width=0.238\textwidth]{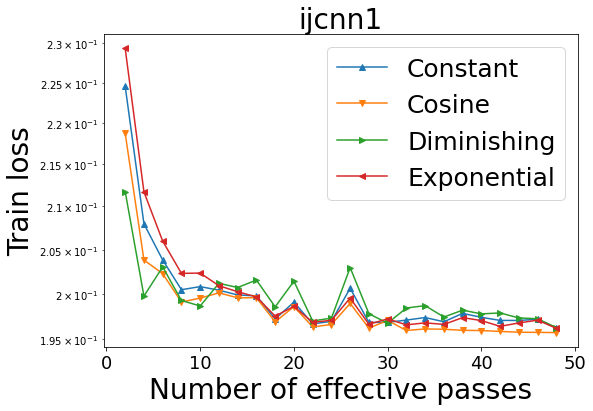}
\caption{The train loss produced by SMG using four different  learning rates on the \texttt{w8a} and \texttt{ijcnn1} datasets, respectively.}
\label{fig_LR_04}
\end{figure}

\textbf{Results.}
We conduct this test using four different learning rate variants: constant, diminishing, exponential decay, and cosine annealing learning rates.
Our results are reported in  Figure~\ref{fig_NN_04} and Figure~\ref{fig_LR_04}.
From Figure \ref{fig_NN_04}, we can observe that the cosine scheduled  and the diminishing learning rates converge relatively fast at the early stage. 
However, the exponential decay and the constant learning rates make faster progress in the last epochs and tend to give the best result at the end of the training process in our neural network training experiments. 

For the non-convex logistic regression, Figure~\ref{fig_LR_04} shows that the cosine learning rate has certain advantages compared to other choices after a few dozen numbers of epochs.

We note that the detailed settings and additional experiments in this section can be found in the Supp. Doc. 

\section{Conclusions and Future Work}\label{sec_conclusion}

We have developed two new shuffling gradient algorithms with momentum for solving nonconvex finite-sum minimization problems. 
{Our Algorithm~\ref{sgd_momentum_shuffling_01} is novel and can work with any shuffling strategy, while Algorithm~\ref{sgd_momentum_shuffling2} is similar to existing momentum methods using single shuffling strategy.}
Our methods achieve the state-of-the-art $\Ocal(1/T^{2/3})$ convergence rate under standard assumptions using different learning rates and shuffling strategies. 
When a randomized reshuffling strategy is exploited {for Algorithm~\ref{sgd_momentum_shuffling_01}}, we can further improve our rates by a fraction $n^{1/3}$ of the data size $n$, matching the best known results in the non-momentum shuffling method.
Our numerical results also show encouraging performance of the new methods.

We  believe  that  our  analysis  framework  can  be  extended  to  study  non-asymptotic rates for some recent adaptive SGD schemes such as Adam \citep{Kingma2014} and AdaGrad \citep{AdaGrad} as well as variance-reduced methods such as SARAH \citep{Nguyen2017sarah} under shuffling strategies.
It is also interesting to extend our framework to other problems such as minimax and federated learning.
We will further investigate these opportunities in the future.

\section*{Acknowledgements}

The authors would like to thank the reviewers for their useful comments and suggestions which helped to improve the exposition in this paper. The work of Q. Tran-Dinh has partly been supported by the Office of Naval Research under Grant No. ONR-N00014-20-1-2088 (2020-2023).

\newpage
\bibliographystyle{icml2021}
\bibliography{all_refs,reference}

\clearpage
\onecolumn
\begin{center}
\textsc{\large Supplementary Document:}

\textbf{\Large SMG: A Shuffling Gradient-Based Method with Momentum}
\end{center}
\vspace{2ex}

\section{Technical Proofs of Theorem~\ref{thm_momentum_variance} and Its Consequences in Section~\ref{sec:convergence_analysis}}
\label{apdx:sec:convergence_analysis}
This section provide the full proofs of Theorem~\ref{thm_momentum_variance} and its consequences  in Section~\ref{sec:convergence_analysis}.

\subsection{Auxiliary Notation and Common Expressions}\label{subsec:technicality}
\begin{remark}
We use the superscript ``$(t)$" for the epoch counter, and the subscript ``$i$" for the counter of the inner loop of the $t$-th epoch.
In addition, at each epoch $t$, we fix the learning rate $\eta_i^{(t)} := \frac{\eta_t}{n}$ for given $\eta_t > 0$. 
\end{remark}
In this Supp. Doc. we repeatedly use some common notations to analyze the analysis for three different shuffling strategies, including: Randomized Reshuffling, Shuffle Once and Incremental Gradient. We are ready to introduce the notations here. 

For any $t \geq 1$, $i=0,\dots, n$, $w_0^{(t)}\in\R^d$ generated by Algorithm \ref{sgd_momentum_shuffling_01}, and two permutations $\pi^{(t)}$ and $\pi^{(t-1)}$ of $[n]$ generated at epochs $t$ and $t-1$, we denote
\begin{align*}
    & A_i^{(t)} := \Big\Vert  \sum_{j=0}^{i-1}  g_{j}^{(t)} \Big\Vert ^2,  \tagthis \label{define_A} \\ 
    & B_i^{(t)} := \Big\Vert \sum_{j = i}^{n-1} g_{j}^{(t)}\Big\Vert^2,    \tagthis \label{define_B}\\
    & I_t := \sum_{i=0}^{n-1} \Big \Vert \nabla f ( w_{0}^{(t)} ; \pi^{(t)} ( i + 1 ) )  - g_{i}^{(t)} \Big\Vert^2,  \tagthis \label{define_I} \\ 
    & J_t := \sum_{i=0}^{n-1} \Big \Vert \nabla f ( w_{0}^{(t)} ; \pi^{(t-1)} ( i + 1 ) )  - g_{i}^{(t-1)} \Big\Vert^2, \quad  t \geq 2, \tagthis \label{define_J} \\
    & K_{t} := n (\Theta +1) \big\Vert \nabla F(w_{0}^{(t)}) \big\Vert ^2  + n \sigma^2, \tagthis \label{define_K} \vspace{1.5ex}\\
    & N_t =(n +\Theta) \big\Vert \nabla  F(w_{0}^{(t)})  \big\Vert ^2 + \sigma^2, \tagthis \label{define_N}
\end{align*}
and
\begin{equation}\label{define_xi}
\xi_t := \left\{\begin{array}{ll}
\max (\eta_t, \eta_{t-1})  &\text{if} \ t > 1,\\
\eta_1 &\text{if} \ t = 1.
\end{array}\right.
\end{equation} 
Note that we adopt the convention $\sum_{j=0}^{-1}  g_{j}^{(t)} = 0$, and $\sum_{j=n}^{n-1}  g_{j}^{(t)} = 0$ in the definitions of $A_i^{(t)} $ and $B_i^{(t)} $.

For each epoch $t=1, \cdots, T$, we denote  $\mathcal{F}_t $ by $ \sigma(w_0^{(1)},\cdots,w_0^{(t)})$, the $\sigma$-algebra generated by the iterates of Algorithm~\ref{sgd_momentum_shuffling_01} up to the beginning of the epoch $t$. From Step 4 of Algorithm~\ref{sgd_momentum_shuffling_01}, we observe that the permutation $\pi^{(t)}$ used at time $t$ is independent of the $\sigma$-algebra $\mathcal{F}_t$. This will be the key factor of our analysis for Randomized Reshuffling methods. We also denote $\E_t[\cdot] $ by $\E [\cdot \mid \mathcal{F}_t]$, the conditional expectation 
on the $\sigma$-algebra $\mathcal{F}_t$.

Now, we collect the following expressions derived from Algorithm~\ref{sgd_momentum_shuffling_01}, which are commonly used for every shuffling method.
First, from the update rule at Steps 1, 4, and 5 of Algorithm~\ref{sgd_momentum_shuffling_01}, for $t > 1$, we have
\begin{equation}\label{update_m_0^{T}}
    m_0^{(t)} = \tilde{m}_{t-1} = v_{n}^{(t-1)} \overset{\eqref{eq:main_update}}{=} v_{n-1}^{(t-1)} + \frac{1}{n} g_{n-1}^{(t-1)}
    = v_{0}^{(t-1)} + \frac{1}{n} \sum_{j=0}^{n-1} g_{j}^{(t-1)} = \frac{1}{n} \sum_{j=0}^{n-1} g_{j}^{(t-1)}.
\end{equation}
For the special case when $t = 1$, we set $m_0^{(1)} = \tilde{m}_{0} = \textbf{0} \in\R^d$.

Next, from the update $w_{i+1}^{(t)} := w_{i}^{(t)} - \eta_i^{(t)} m_{i+1}^{(t)}$  with $\eta_i^{(t)} := \frac{\eta_t}{n}$, for $i = 1, 2,\dots ,n-1$, at Step \ref{alg:A1_step4} of Algorithm~\ref{sgd_momentum_shuffling_01}, we can derive
\begin{align*}
    w_{i}^{(t)} &\overset{\eqref{eq:main_update}}{=} w_{i-1}^{(t)} - \frac{\eta_t}{n} m_{i}^{(t)} = w_{0}^{(t)} - \frac{\eta_t}{n} \sum_{j=0}^{i-1} m_{j+1}^{(t)}, \quad  t \geq 1, \tagthis \label{update_w_i^{T}}\\
    w_{0}^{(t)} = w_{n}^{(t-1)} &\overset{\eqref{eq:main_update}}{=} w_{n-1}^{(t-1)} - \frac{\eta_{t-1}}{n} m_{n}^{(t-1)} = w_{i}^{(t-1)} - \frac{\eta_{t-1}}{n} \sum_{j = i}^{n-1} m_{j+1}^{(t-1)}, \quad t \geq 2. \tagthis \label{update_w_i^{T}-1}
\end{align*}
Note that $\sum_{j=0}^{-1} m_{j+1}^{(t)}$ and $\sum_{j = n}^{n-1} m_{j+1}^{(t-1)} = 0$ by convention, these equations also hold true for $i=0$ and $i=n$.\\


%

%


\subsection{The Proof Sketch of Theorem~\ref{thm_momentum_variance}}\label{subsec:proof_TH1}
Due to the technicality of our proofs, we first provide here a proof sketch of our first main result, Theorem~\ref{thm_momentum_variance}, while the full proof is given in the next subsections.

The proof of Theorem~\ref{thm_momentum_variance} is divided into the following steps.
\begin{compactitem}
\item First, from the $L$-smoothness of $f_i$ and the update of $w^{(t)}_i$, we can bound 
\begin{equation*}
F(\tilde{w}_{t}) \leq  F(\tilde{w}_{t-1}) - \frac{\eta_t}{2}\norms{\nabla{F}(\tilde{w}_{t-1})}^2 + \mathcal{T}_{[1]}, 
\end{equation*}
where 
\begin{equation*}
\mathcal{T}_{[1]} := \frac{\eta_t}{2} \big\Vert \nabla{F}( \tilde{w}_{t-1})  - \frac{1}{n} \sum_{j=0}^{n-1} \big( \beta g_{j}^{(t-1)} + (1-\beta)   g_{j}^{(t)} \big) \big\Vert^2.
\end{equation*}
Hence, the key step here is to upper bound $\mathcal{T}_{[1]}$ by $\sum_{t=1}^T\eta_t\norms{\nabla{F}(\tilde{w}_{t-1})}^2$.
\item Next, we upper bound $\mathcal{T}_{[1]} \leq \frac{\eta_t}{2n}(\beta I_t + (1-\beta)J_t) $, where $I_t$ and $J_t$ are defined by \eqref{define_I} and \eqref{define_J}, respectively, which bound the sum of square errors between $g_i^{(t)}$ and the gradient $\nabla{f}(w_0^{(t)}, \pi^{(t)}(i+1))$ at $t$ and $t-1$, respectively.

\item Then, $\beta I_t + (1-\beta)J_t$ can be upper bounded as
\begin{equation*}
\beta I_t + (1-\beta)J_t \leq \Ocal\Big(\xi_t^2\sum_{j=0}^t\beta^{t-j}\norms{\nabla{F}(w_0^{(j)})}^2 + \xi^2_t\sigma^2 \Big),
\end{equation*}
as shown in Lemma~\ref{lem_momentum_variance_04}, where $\xi_t := \max\{\eta_{t-1}, \eta_t\}$ defined in \eqref{define_xi}.
\item We further upper bound the sum of the right-hand side of the last estimate as in Lemma~\ref{lem_momentum_variance_05} in terms $\sum_{t=0}^T\eta_t\norms{\nabla{F}(\tilde{w}_{t-1})}^2$.
\end{compactitem}
Combining these steps together, and using some simplification, we obtain \eqref{eq_thm_momentum_variance} in Theorem~\ref{thm_momentum_variance}.

\subsection{Technical Lemmas}
The proof of Theorem~\ref{thm_momentum_variance} requires the following lemmas as intermediate steps of its proof.
Let us first upper bound $\big\Vert w_{i}^{(t)} - w_{0}^{(t)} \big\Vert^2$ and $\big\Vert w_{i}^{(t-1)} - w_{0}^{(t)} \big\Vert^2$ in the following lemma.

\begin{lem} \label{lem_momentum_variance_01}
Let $\{w_i^{(t)}\}$ be generated by Algorithm \ref{sgd_momentum_shuffling_01} with $0\leq \beta < 1$ and $\eta_i^{(t)} := \frac{\eta_t}{n}$ for every $t \geq 1$.
Then, for $i = 0, 1, \dots, n$, we have 
\begin{align*}
    (\textbf{a}) \ &\big\Vert w_{i}^{(t)} - w_{0}^{(t)} \big\Vert^2  \leq
    \left\{\begin{array}{ll}
    \frac{\xi_t^2}{n^2} \left[ \beta \frac{i^2}{n^2} A_n^{(t-1)} + (1-\beta) A_i^{(t)} \right], &\text{if}~t > 1, \\ 
    \frac{\xi_t^2}{n^2}  (1-\beta)^2 A_i^{(t)}, &\text{if}~ t = 1, 
    \end{array}\right. \tagthis \label{update_distance_t} \\
    (\textbf{b}) \ &\big\Vert w_{i}^{(t-1)} - w_{0}^{(t)} \big\Vert^2 \leq 
    \left\{\begin{array}{ll}
    \frac{\xi_t^2}{n^2}  \left[ \beta \frac{(n-i)^2}{n^2} A_n^{(t-2)} + (1-\beta) B_i^{(t-1)} \right] , &\text{if}~t > 2, \\
    \frac{\xi_t^2}{n^2}   (1-\beta)^2 B_i^{(t-1)}, &\text{if}~ t = 2.
    \end{array}\right.
    \tagthis \label{update_distance_t-1}
\end{align*}
\end{lem}

Next, we upper bound the quantities $A_n^{(t)}$, $\sum_{i=0}^{n-1}A_i^{(t)}$, and $\sum_{i=0}^{n-1}B^{(t)}_i$ defined above in terms of $I_t$ and $K_t$.
\begin{lem} \label{lem_momentum_variance_02}
Under the same setting as of Lemma~\ref{lem_momentum_variance_01} and Assumption~\ref{as:A1}(c), for $t\geq 1$, it holds that
\begin{align*}
    (\textbf{a}) \qquad & A_n^{(t)} \leq 3 n  \Big( I_t + K_{t} \Big)
    \quad \text{and} \quad \sum_{i=0}^{n-1} A_i^{(t)} \leq \frac{3n^2}{2} \Big( I_t + K_t\Big), \tagthis \label{bound_A^{T}}\\
    (\textbf{b}) \qquad & \sum_{i=0}^{n-1} B_i^{(t)} \leq 3n^2 \Big(I_t + K_t \Big). \tagthis \label{bound_B^{T}}
\end{align*}
\end{lem}

The results of Lemma~\ref{lem_momentum_variance_03} below are direct consequences of Lemma \ref{lem_momentum_variance_01}, Lemma~\ref{lem_momentum_variance_02}, and the fact that $L^2\xi_t^2 \leq \frac{2}{5}$.
\begin{lem} \label{lem_momentum_variance_03}
Under the same setting as of Lemma~\ref{lem_momentum_variance_02} and Assumption~\ref{as:A1}(b), it holds that
\begin{align*}
    (\textbf{a}) \ 
    &I_t \leq \left\{\begin{array}{ll} 
    L^2\xi_t^2 \Big[\beta (  I_{t-1} + K_{t-1} )  + \frac{3}{2} (1-\beta) (I_t + K_t ) \Big], &\text{if}~t > 1,  \\
    \frac{3}{2} L^2\xi_t^2  (1-\beta)^2 \Big(I_t + K_t \Big), &\text{if}~t = 1,
    \end{array}\right. \\
    (\textbf{b}) \ 
    &J_t \leq 
    \left\{\begin{array}{ll}
    3L^2\xi_t^2 \Big[ \beta  (  I_{t-2} + K_{t-2} )  + (1-\beta)  (I_{t-1} + K_{t-1} ) \Big], &\text{if}~ t > 2, \\
    3L^2\xi_t^2  \Big(I_{t-1} + K_{t-1} \Big),  &\text{if}~ t =2,
    \end{array}\right.
    \\
    (\textbf{c}) \ &\text{If } L^2\xi_t^2 \leq \frac{2}{5}, \text{ then }  I_t + K_t \leq 
    \left\{\begin{array}{ll}
    \beta \Big(I_{t-1} + K_{t-1}\Big) + \frac{5 -3\beta}{2} K_t, &\text{if}~ t > 1,   \\ 
    \frac{5 -3\beta}{2} K_1, &\text{if}~t = 1.  
    \end{array}\right.
\end{align*}
\end{lem}

One of our key estimates is to upper bound the quantity $\beta J_t + (1-\beta) I_t$, which is derived in the following lemma.

\begin{lem} \label{lem_momentum_variance_04}
Under the same conditions as in Lemma~\ref{lem_momentum_variance_03}, for any $t\geq 2$, we have
\begin{align*}
    \beta J_t + (1-\beta) I_t & \leq  \frac{9(5 -3\beta) L^2\xi_t^2}{2} \left( n (\Theta +1) \sum_{j=1}^{t}  \beta^{t-j} \big\Vert\nabla  F(w_{0}^{(j)}) \big\Vert ^2 +  \frac{n \sigma^2}{1-\beta} \right). 
    \tagthis \label{bound_beta_IJ}
\end{align*}
\end{lem}

\begin{lem} \label{lem_momentum_variance_05}
Under the same settings as in Theorem \ref{thm_momentum_variance} (or Theorem \ref{thm_momentum_variance_RR}, respectively) and $\eta_t \geq \eta_{t+1}$ for $t\geq 1$, we have 
\begin{align*}
    \sum_{t=1}^{T} \eta_t \sum_{j=1}^{t}  \beta^{t-j} \big\Vert \nabla F(\tilde{w}_{j-1}) \big\Vert ^2 \leq \frac{1}{1-\beta}  \sum_{t=1}^{T} \eta_t \big\Vert \nabla F(\tilde{w}_{t-1}) \big\Vert ^2.
\end{align*}
\end{lem}

\begin{lem} \label{lem_momentum_variance_06}
Under the same setting as in Theorem \ref{thm_momentum_variance}, we have
\begin{align*} 
    \frac{\eta_1}{2} \norms{ \nabla F(\tilde{w}_{0} )}^2  \leq \frac{ F( \tilde{w}_{0} ) - F( \tilde{w}_{1} )}{(1-\beta)}
    + \frac{(1-\beta)\eta_1 }{4} \big\Vert \nabla F(\tilde{w}_{0}) \big\Vert ^2  +  \frac{9\sigma^2(5 - 3\beta)L^2}{4(1-\beta)}\cdot \xi_1^3. \tagthis \label{bound_F(w)_1_tilde}
\end{align*}

\end{lem}

\subsection{The Proof of Theorem \ref{thm_momentum_variance}: Key Estimate for Algorithm~\ref{sgd_momentum_shuffling_01}}
First, from the assumption $0 < \eta_t \leq \frac{1}{L\sqrt{K}}$, $t \geq 1$, we have $0 < \eta_t^2 \leq \frac{1}{K L^2}$. 
Next, from \eqref{define_xi}, we have $\xi_t = \max (\eta_t; \eta_{t-1})$ for $t > 1$ and $\xi_1 = \eta_1$, which lead to $0 < \xi_t^2 \leq \frac{1}{K L^2}$ for $t \geq 1$. 
Moreover, from the definition of $K = \max \left(\frac{5}{2}, \frac{9(5 -3\beta) (\Theta +1)}{1-\beta}\right)$ in Theorem~\ref{thm_momentum_variance}, we have $L^2\xi_t^2 \leq \frac{2}{5}$ and $9 L^2 \xi_t^2 (5 -3\beta) (\Theta +1) \leq 1-\beta$ for $t\geq 1.$




Now, letting $i=n$ in Equation \eqref{update_w_i^{T}}, for all $t\geq 1$, we have
\begin{align*}
    w_{n}^{(t)} - w_{0}^{(t)} &= - \frac{\eta_t}{n} \sum_{j=0}^{n-1} m_{j+1}^{(t)} 
    \overset{\eqref{eq:main_update}}{=}  - \frac{\eta_t}{n} \sum_{j=0}^{n-1} \big( \beta m_{0}^{(t)} + (1-\beta) g_{j}^{(t)}\big) 
    =  - \frac{\eta_t}{n} \Big( n \beta m_{0}^{(t)} +  (1-\beta)  \sum_{j=0}^{n-1}  g_{j}^{(t)}  \Big).
\end{align*}
Since $m_0^{(t)} = \frac{1}{n} \sum_{j=0}^{n-1} g_{j}^{(t-1)}$ for $t \geq 2$ (due to \eqref{update_m_0^{T}}), we have the following update
\begin{align*}
    w_{n}^{(t)} - w_{0}^{(t)} &= - \frac{\eta_t}{n}  \sum_{j=0}^{n-1} \Big( \beta g_{j}^{(t-1)} + (1-\beta)   g_{j}^{(t)} \Big).
    \tagthis \label{update_epoch_02}
\end{align*}
From the $L$-smoothness of $F$ in Assumption \eqref{as:A1}(b), \eqref{update_epoch_02}, and the fact that $w_0^{(t+1)} := w_{n}^{(t)}$, for every epoch $t\geq 2$, we can derive
\begin{align*} 
    F( w_0^{(t+1)} )  & \overset{\eqref{eq:Lsmooth}}{\leq} F( w_0^{(t)} ) + \nabla F( w_0^{(t)} )^{\top}(w_0^{(t+1)} - w_0^{(t)}) + \frac{L}{2}\norms{w_0^{(t+1)} - w_0^{(t)}}^2  \\
    &\overset{\eqref{update_epoch_02}}{=} F( w_0^{(t)} ) - \eta_t \nabla F( w_0^{(t)} )^{\top} \left( \frac{1}{n} \sum_{j=0}^{n-1} \Big( \beta g_{j}^{(t-1)} + (1-\beta)   g_{j}^{(t)} \Big) \right)  \\
    & \quad + \frac{L \eta_t^2}{2} \Big \Vert \frac{1}{n} \sum_{j=0}^{n-1} \Big( \beta g_{j}^{(t-1)} + (1-\beta)   g_{j}^{(t)} \Big) \Big\Vert ^2  \\
    & \overset{\tiny(a)}{=} F( w_0^{(t)} ) - \frac{\eta_t}{2} \norms{ \nabla F( w_0^{(t)} )}^2 + \frac{\eta_t}{2} \Big\Vert \nabla F( w_0^{(t)} )  - \frac{1}{n} \sum_{j=0}^{n-1} \left( \beta g_{j}^{(t-1)} + (1-\beta)   g_{j}^{(t)} \right)\Big\Vert^2  \\
    & \quad - \frac{\eta_t}{2} \left( 1 - L\eta_t  \right) \Big\Vert \frac{1}{n} \sum_{j=0}^{n-1} \left( \beta g_{j}^{(t-1)} + (1-\beta)   g_{j}^{(t)} \right) \Big\Vert^2   \\
    & \leq F( w_0^{(t)} ) - \frac{\eta_t}{2} \norms{ \nabla F( w_0^{(t)} )}^2 + \frac{\eta_t}{2} \Big\Vert \nabla F( w_0^{(t)} )  - \frac{1}{n} \sum_{j=0}^{n-1} \left( \beta g_{j}^{(t-1)} + (1-\beta)   g_{j}^{(t)} \right)\Big\Vert^2, \tagthis \label{bound_F(w)}
\end{align*}
where (\textit{a}) follows from the elementary equality $u^{\top}v = \frac{1}{2}\left(\norms{u}^2 + \norms{v}^2 - \norms{u - v}^2\right)$ and the last equality comes from the fact that $\eta_t \leq \frac{\sqrt{2}}{\sqrt{5} L} \leq \frac{1}{L}$, due to the choice of our learning rate $\eta_t$. 

Next, we use an interesting fact of the permutations $\pi^{(t)}$ and $\pi^{(t-1)}$ to rewrite the full gradient as
\begin{align*}
    \nabla F( w_0^{(t)} )  &= \beta \nabla F( w_0^{(t)} ) + (1- \beta) \nabla F( w_0^{(t)} )\\
    &=  \frac{\beta}{n} \sum_{j=0}^{n-1} \nabla f ( w_{0}^{(t)} ; \pi^{(t-1)} ( j + 1 ) ) +  \frac{(1-\beta)}{n} \sum_{j=0}^{n-1} \nabla f ( w_{0}^{(t)} ; \pi^{(t)} ( j + 1 ) ) \\
    &= \frac{1}{n} \sum_{j=0}^{n-1} \left(\beta \nabla f ( w_{0}^{(t)} ; \pi^{(t-1)} ( j + 1 ) ) + (1-\beta) \nabla f ( w_{0}^{(t)} ; \pi^{(t)} ( j + 1 ) ) \right).
\end{align*}
Let us upper bound the last term of \eqref{bound_F(w)} as follows:
\begin{align*} 
    \mathcal{T}_{[1]} &:= \frac{\eta_t}{2} \Big\Vert \nabla F( w_0^{(t)} )  - \frac{1}{n} \sum_{j=0}^{n-1} \left( \beta g_{j}^{(t-1)} + (1-\beta)   g_{j}^{(t)} \right)\Big\Vert^2  \\
    & = \frac{\eta_t}{2} \Big\Vert \frac{1}{n} \sum_{j=0}^{n-1} \left[ \beta \Big(\nabla f ( w_{0}^{(t)} ; \pi^{(t-1)} ( j + 1 ) ) - g_{j}^{(t-1)} \Big) + (1-\beta) \Big( \nabla f ( w_{0}^{(t)} ; \pi^{(t)} ( j + 1 ) )  - g_{j}^{(t)} \Big) \right]\Big\Vert^2  \\
    & \overset{(b)}{\leq} \frac{\eta_t}{2n} \sum_{j=0}^{n-1} \Big\Vert 
    \beta \Big(\nabla f ( w_{0}^{(t)} ; \pi^{(t-1)} ( j + 1 ) ) - g_{j}^{(t-1)} \Big) + (1-\beta) \Big( \nabla f ( w_{0}^{(t)} ; \pi^{(t)} ( j + 1 ) )  - g_{j}^{(t)} \Big) \Big\Vert^2 \\
    & \overset{(c)}{\leq} \frac{\eta_t}{2n} \sum_{j=0}^{n-1} \left[ \beta \Big\Vert 
    \nabla f ( w_{0}^{(t)} ; \pi^{(t-1)} ( j + 1 ) ) - g_{j}^{(t-1)} \Big\Vert^2  + (1-\beta) \Big \Vert \nabla f ( w_{0}^{(t)} ; \pi^{(t)} ( j + 1 ) )  - g_{j}^{(t)} \Big\Vert^2 \right] \\
    &\overset{\eqref{define_I},\eqref{define_J}}{\leq} \frac{\eta_t}{2n} \Big[ \beta J_t + (1-\beta) I_t \Big] \qquad \textrm{(by the definition of $I_t$ and $J_t$)},
\end{align*}
where (\textit{b}) is from the Cauchy-Schwarz inequality, (\textit{c}) follows from the convexity of $\norms{\cdot}^2$ for $0\leq \beta < 1$.



Applying this into \eqref{bound_F(w)}, we get the following result:
\begin{align*} 
    F( w_0^{(t+1)} )& \leq F( w_0^{(t)} ) - \frac{\eta_t}{2} \norms{ \nabla F( w_0^{(t)} )}^2 + \frac{\eta_t}{2n} \Big[ \beta J_t + (1-\beta) I_t \Big]. \tagthis \label{estimate_alg_1}
\end{align*}

\begin{remark}
    Note that the derived estimate \eqref{estimate_alg_1} is true for all shuffling strategies, including Random Reshuffling, Shuffle Once, and Incremental Gradient (under the assumptions of Theorem \ref{thm_momentum_variance} and \ref{thm_momentum_variance_RR}, respectively).
    
\end{remark}

Applying the result of Lemma \ref{lem_momentum_variance_04} for $t \geq 2$, then we obtain
\begin{align*} 
    F( w_0^{(t+1)} ) & \leq F( w_0^{(t)} ) - \frac{\eta_t}{2} \norms{ \nabla F( w_0^{(t)} )}^2 + \frac{9(5 -3\beta) \eta_t L^2\xi_t^2}{4n} \Bigg[ n (\Theta+1) \sum_{j=1}^{t}  \beta^{t-j} \big\Vert \nabla F(w_{0}^{(j)}) \big\Vert ^2  +  \frac{n \sigma^2}{1-\beta} \Bigg],
\end{align*}
which leads to 
\begin{align*} 
	F( w_0^{(t+1)} ) & \leq F( w_0^{(t)} ) - \frac{\eta_t}{2} \norms{ \nabla F( w_0^{(t)} )}^2 + \frac{9(5 -3\beta)\eta_t L^2 \xi_t^2(\Theta+1)}{4} \sum_{j=1}^{t}  \beta^{t-j} \big\Vert \nabla F(w_{0}^{(j)}) \big\Vert ^2  + \frac{9 \sigma^2 (5 -3\beta)L^2 \xi_t^3 }{4(1-\beta)}. 
\end{align*}
Since $\xi_t$ satisfies $9 L^2 \xi_t^2 (5 -3\beta) (\Theta+1) \leq 1-\beta$ as proved above, we can deduce from the last estimate that
\begin{align*} 
    F( w_0^{(t+1)} )
    & \leq F( w_0^{(t)} ) - \frac{\eta_t}{2} \norms{ \nabla F( w_0^{(t)} )}^2 + \frac{(1-\beta)\eta_t }{4} \sum_{j=1}^{t}  \beta^{t-j} \big\Vert \nabla F(w_{0}^{(j)}) \big\Vert ^2  + \frac{9\sigma^2(5 -3\beta)L^2 }{4(1-\beta)} \cdot \xi_t^3.
\end{align*}
Rearranging this inequality and noting that $\tilde{w}_{t-1} = w_0^{(t)}$, for $t \geq 2$, we obtain the following estimate: 
\begin{equation*}
\arraycolsep=0.2em
\begin{array}{lcl} 
    \frac{\eta_t}{2} \norms{ \nabla F( \tilde{w}_{t-1} )}^2     & \leq & F( \tilde{w}_{t-1} ) - F( \tilde{w}_t ) +  \frac{(1-\beta)\eta_t}{4} \sum_{j=1}^{t}  \beta^{t-j} \big\Vert \nabla F(\tilde{w}_{j-1}) \big\Vert ^2  +  \frac{9\sigma^2(5 -3\beta)L^2}{4(1-\beta)} \cdot \xi_t^3
     \\
    & \overset{1-\beta <1}{\leq} & \frac{F( \tilde{w}_{t-1} ) - F( \tilde{w}_t )}{(1-\beta)} +  \frac{(1-\beta)\eta_t}{4} \sum_{j=1}^{t}  \beta^{t-j} \big\Vert \nabla F(\tilde{w}_{j-1}) \big\Vert ^2  +  \frac{9\sigma^2(5 -3\beta)L^2}{4(1-\beta)} \cdot \xi_t^3. 
\end{array}
\end{equation*}
For $t = 1$, since $\xi_1 = \eta_1$ as previously defined in \eqref{define_xi}, from the result of Lemma \ref{lem_momentum_variance_06} we have 
\begin{align*} 
    \frac{\eta_1}{2} \norms{ \nabla F(\tilde{w}_{0} )}^2  \leq \frac{ F( \tilde{w}_{0} ) - F( \tilde{w}_{1} )}{(1-\beta)}
    + \frac{(1-\beta)\eta_1 }{4} \sum_{j=1}^{1}  \beta^{1-j} \big\Vert \nabla F(\tilde{w}_{j-1}) \big\Vert ^2  +  \frac{9\sigma^2(5 - 3\beta)L^2}{4(1-\beta)}\cdot \xi_1^3.
\end{align*}

Summing the previous estimate for $t := 2$ to $t := T$, and then using the last one, we obtain
\begin{align*} 
    \frac{1}{2} \sum_{t=1}^{T} \eta_t \norms{ \nabla F( \tilde{w}_{t-1} )}^2 
    & \leq \frac{ F( \tilde{w}_0 ) - F_*}{(1-\beta)} +  \frac{(1-\beta)}{4} \sum_{t=1}^{T}  \eta_t \sum_{j=1}^{t}  \beta^{t-j} \big\Vert \nabla F(\tilde{w}_{j-1}) \big\Vert ^2  +  \frac{9 \sigma^2 L^2 (5 -3\beta)}{4(1-\beta)} \sum_{t=1}^{T} \xi_t^3. 
\end{align*}

Applying Lemma \ref{lem_momentum_variance_05} to the last estimate, we get
\begin{align*} 
    \frac{1}{2} \sum_{t=1}^{T} \eta_t \norms{ \nabla F( \tilde{w}_{t-1} )}^2 
    & \leq \frac{ F( \tilde{w}_0) - F_* }{(1-\beta)} +  \frac{(1-\beta)}{4} \frac{1}{(1-\beta)} \sum_{t=1}^{T} \eta_t \norms{ \nabla F( \tilde{w}_{t-1} )}^2   +  \frac{9 \sigma^2 L^2 (5 -3\beta)}{4(1-\beta)}\sum_{t=1}^{T} \xi_t^3, 
\end{align*}
which is equivalent to 
\begin{align*} 
    \sum_{t=1}^{T} \eta_t \norms{ \nabla F( \tilde{w}_{t-1} )}^2 
    & \leq \frac{ 4[F( \tilde{w}_0 ) - F_* ]}{(1-\beta)} +  \frac{9 \sigma^2 L^2 (5 -3\beta)}{(1-\beta)} \sum_{t=1}^{T} \xi_t^3. 
    \tagthis \label{eq_thm_01}
\end{align*}
Dividing both sides of the resulting estimate by $\sum_{t=1}^{T} \eta_t$, we obtain \eqref{eq_thm_momentum_variance}.
Finally, due to the choice of $\hat{w}_T$ at Step~\ref{alg:A1_step6} in Algorithm~\ref{sgd_momentum_shuffling_01}, we have $\mathbb{E}\big[\| \nabla F( \hat{w}_T )  \|^2 \big]  = \frac{1}{\sum_{t=1}^T \eta_t} \sum_{t=1}^T \eta_t  \norms{ \nabla F( \tilde{w}_{t-1} )}^2 $.
\Eproof

\subsection{The Proof of Corollaries~\ref{co:constant_LR} and \ref{co:diminishing_LR}: Constant and Diminishing Learning Rates}\label{apdx:subsec:main_result}

\begin{proof}[The proof of Corollary~\ref{co:constant_LR}]
Since $T \geq 1$ and $\eta_t := \frac{\gamma}{T^{1/3}}$, we have $\eta_t \geq \eta_{t+1}$. We also have $0 < \eta_t \leq \frac{1}{L\sqrt{K}}$ for all $t \geq 1$, and $\sum_{t=1}^T\eta_t = \gamma T^{2/3}$ and $\sum_{t=1}^T\eta_{t-1}^3 = \gamma^3$.
Substituting these expressions into \eqref{eq_thm_momentum_variance}, we obtain
\begin{align*} 
    \mathbb{E}\big[\| \nabla F( \hat{w}_T )  \|^2 \big] 
    & \leq \frac{ 4[F( \tilde{w}_0 ) - F_* ]}{(1-\beta)\gamma T^{2/3}} +  \frac{9 \sigma^2 L^2 (5 -3\beta)}{(1-\beta)} \cdot \frac{\gamma^2}{T^{2/3}} .
\end{align*}
which is our desired result. 
\end{proof}


\begin{proof}[The proof of Corollary~\ref{co:diminishing_LR}]
For $t \geq 1$, since $\eta_t = \frac{\gamma}{(t+\lambda)^{1/3}}$, we have $\eta_t \geq \eta_{t+1}$. We also have $0 < \eta_t \leq \frac{1}{L\sqrt{K}}$ for all $t \geq 1$, and
\begin{equation*}
\begin{array}{lll}
& \sum_{t=1}^T\eta_t  & = \gamma \sum_{t=1}^T\frac{1}{(t+\lambda)^{1/3}} \geq \gamma \int_1^T\frac{d\tau}{(\tau + \lambda)^{1/3}} \geq \gamma \left[ (T + \lambda)^{2/3} - (1 + \lambda)^{2/3} \right],  \\
& \sum_{t=1}^T\eta_{t-1}^3 & = 2\eta_1^3 + \gamma^3 \sum_{t=3}^{T} \frac{1}{t - 1 + \lambda} 
\leq \frac{2\gamma^3}{(1+\lambda)} + \gamma^3 \int_{t=2}^T\frac{d\tau}{\tau - 1 + \lambda} \leq \gamma^3 \left[ \frac{2}{(1+\lambda)} + \log (T - 1 +\lambda) \right].
\end{array}
\end{equation*}
Substituting these expressions into \eqref{eq_thm_momentum_variance}, we obtain
\begin{align*} 
    \mathbb{E}\big[\| \nabla F( \hat{w}_T )  \|^2 \big] 
    & \leq \frac{ 4[F( \tilde{w}_0 ) - F_* ]}{(1-\beta)\gamma \left[ (T + \lambda)^{2/3} - (1 + \lambda)^{2/3} \right]} +  \frac{9 \sigma^2 L^2 (5 -3\beta)}{(1-\beta)} \cdot \frac{\gamma^2 \left[ \frac{2}{(1+\lambda)} + \log (T - 1 +\lambda) \right]}{\left[ (T + \lambda)^{2/3} - (1 + \lambda)^{2/3} \right]}. 
\end{align*}
Let us define $C_1$ and $C_2$ as follows:
\begin{align*}
    C_1 := \frac{4\left[F(  \tilde{w}_{0} )- F_*\right]}{(1-\beta)\gamma }  + \frac{18 \sigma^2 L^2 (5 -3\beta)\gamma^2}{(1-\beta)(1+\lambda)} , \text{ and } C_2:= \frac{9 \sigma^2 L^2 (5 -3\beta)\gamma^2}{(1-\beta)}.
\end{align*}
Then, we obtain from the last estimate that
\begin{equation*} 
    \mathbb{E}\big[\| \nabla F( \hat{w}_T )  \|^2 \big]    \leq \frac{C_1 + C_2 \log (T - 1 +\lambda)}{ (T + \lambda)^{2/3} - (1 + \lambda)^{2/3}},
\end{equation*}
which completes our proof.
\end{proof}

\subsection{The Proof of Corollaries~\ref{co:exponential_LR} and \ref{co:cos_LR}: Scheduled Learning Rates}\label{apdx:subsec:scheduled_LRs}

\begin{proof}[The proof of Corollary \ref{co:exponential_LR}: Exponential LR]
For $t \geq 1$, since $\eta_t =  \frac{\gamma \alpha^t}{T^{1/3}}$ where $\alpha := \rho^{1/T} \in (0, 1)$, we have $\eta_t \geq \eta_{t+1}$. We also have $0 < \eta_t \leq \frac{1}{L\sqrt{K}}$ for all $t \geq 1$ and
\begin{equation*}
\begin{array}{lll}
& \sum_{t=1}^T\eta_t  & = \frac{\gamma}{T^{1/3}} \sum_{t=1}^{T} \alpha^t \geq \frac{\gamma}{T^{1/3}} \sum_{t=1}^{T} \alpha^T = \gamma \rho T^{2/3},  \\
& \sum_{t=1}^T\eta_{t-1}^3 & \leq \sum_{t=1}^T\eta_{1}^3 =  \gamma^3  \alpha^{3} \leq \gamma^3.
\end{array}
\end{equation*}
Substituting these expressions into \eqref{eq_thm_momentum_variance}, we obtain 
\begin{align*} 
    \mathbb{E}\big[\| \nabla F( \hat{w}_T )  \|^2 \big] 
    & \leq \frac{ 4[F( \tilde{w}_0 ) - F_* ]}{(1-\beta)\gamma \rho T^{2/3}} +  \frac{9 \sigma^2 L^2 (5 -3\beta)}{(1-\beta)}\cdot \frac{\gamma^2}{\rho T^{2/3}}. 
\end{align*}
which is our desired result.
\end{proof}

\begin{proof}[The proof of Corollary \ref{co:cos_LR}: Cosine LR]
First, we would like to show that $\sum_{t=1}^T \cos \frac{t \boldsymbol{\pi}}{T}  = -1$. Let us denote $A := \sum_{t=1}^T \cos \frac{t \boldsymbol{\pi}}{T}$. 
Multiplying this sum  by $\sin \frac{\boldsymbol{\pi}}{2T}$, we get
\begin{align*}
    A \cdot \sin \frac{\boldsymbol{\pi}}{2T} &= \sum_{t=1}^T \cos \frac{t \boldsymbol{\pi}}{T} \sin \frac{\boldsymbol{\pi}}{2T} 
    \overset{(a)}{=} \frac{1}{2} \sum_{t=1}^T \Big( \sin \frac{(2t+1) \boldsymbol{\pi}}{2T} - \sin \frac{ (2t-1) \boldsymbol{\pi}}{2T} \Big) \\
    &=  \frac{1}{2} \Big( \sin \frac{(2T+1) \boldsymbol{\pi}}{2T} - \sin \frac{ \boldsymbol{\pi}}{2T} \Big) 
    = \frac{1}{2} \Big( - \sin \frac{\boldsymbol{\pi}}{2T} - \sin \frac{\boldsymbol{\pi}}{2T} \Big) = - \sin \frac{\boldsymbol{\pi}}{2T},
\end{align*}
where (\textit{a}) comes from the identity
$\cos a \cdot \sin b = \frac{1}{2} (\sin (a+b) - \sin (a-b))$. 
Since $\sin \frac{\boldsymbol{\pi}}{2T}$ is nonzero, we obtain $A = \sum_{t=1}^T \cos \frac{t \boldsymbol{\pi}}{T}  = -1.$ 

For $1 \leq t \leq T$, since $\eta_t = \frac{\gamma}{T^{1/3}}\left(1+ \cos \frac{t \boldsymbol{\pi}}{T} \right)$, we have $\eta_t \geq \eta_{t+1}$. We also have $0 < \eta_t \leq \frac{1}{L\sqrt{K}}$ for all $t \geq 1$ and
\begin{equation*}
\begin{array}{lll}
& \sum_{t=1}^T\eta_t  & = \frac{\gamma}{T^{1/3}} \sum_{t=1}^{T} \left(1+ \cos \frac{t \boldsymbol{\pi}}{T} \right)  \overset{(b)}{=} \frac{\gamma}{T^{1/3}} (T-1) \geq \frac{\gamma}{T^{1/3}} \cdot \frac{T}{2} = \frac{\gamma T^{2/3}}{2}, \ \text{for} \ T \geq 2,   \\
& \sum_{t=1}^T \eta_{t-1}^3 & \leq  \frac{\gamma^3}{T} \sum_{t=1}^{T} 2^3  = 8 \gamma^3 \text{, since } \left(1+ \cos \frac{t \boldsymbol{\pi}}{T} \right)^3 \leq 2^3 \text{ for all } t \geq 1,
\end{array}
\end{equation*}
where $(b)$ follows since $\sum_{t=1}^T \cos \frac{t \boldsymbol{\pi}}{T}  = -1$ as shown above. 
Substituting these expressions into \eqref{eq_thm_momentum_variance}, we obtain
\begin{align*} 
    \mathbb{E}\big[\| \nabla F( \hat{w}_T )  \|^2 \big] 
    & \leq \frac{8[F( \tilde{w}_0 ) - F_* ]}{(1-\beta)\gamma T^{2/3}} +  \frac{9 \sigma^2 L^2 (5 -3\beta)}{(1-\beta)} \cdot \frac{16 \gamma^2}{ T^{2/3}}. 
\end{align*}
which is our desired estimate.
\end{proof}

\subsection{The Proof of Technical Lemmas}\label{appdx:subec:Tech_lemmas_proof}
We now provide the proof of six lemmas that serve for the proof of Theorem~\ref{thm_momentum_variance} above.

\subsubsection{Proof of Lemma~\ref{lem_momentum_variance_01}: Upper Bounding Two Terms $\Vert w_i^{(t)} - w_0^{(t)}\Vert^2$ and $\Vert w_i^{(t-1)} - w_0^{(t)}\Vert^2$}

\begin{remark}
    Note that the result of Lemma~\ref{lem_momentum_variance_01} is true for all shuffling strategies including Random Reshuffling, Shuffle Once and Incremental Gradient (under the assumptions of Theorem \ref{thm_momentum_variance} and \ref{thm_momentum_variance_RR}, respectively).
\end{remark}

(\textbf{a}) 
From Equation \eqref{update_w_i^{T}}, for $i = 0, 1, \dots, n$ and $t \geq 1$, we have
\begin{align*}
    w_{i}^{(t)} - w_{0}^{(t)} &= - \frac{\eta_t}{n} \sum_{j=0}^{i-1} m_{j+1}^{(t)} 
    \overset{\eqref{eq:main_update}}{=}  - \frac{\eta_t}{n} \sum_{j=0}^{i-1} \big( \beta m_{0}^{(t)} + (1-\beta) g_{j}^{(t)}\big) 
    =  - \frac{\eta_t}{n} \Big( i \beta m_{0}^{(t)} +  (1-\beta)  \sum_{j=0}^{i-1}  g_{j}^{(t)}  \Big).
\end{align*}
Moreover, by \eqref{define_xi}, we have $\eta_{t} \leq \xi_t$ for all $t\geq 1$.
Therefore, for $t > 1$ and $i = 0, 1, \dots, n$, we can derive that
\begin{align*}
    \big\Vert w_{i}^{(t)} - w_{0}^{(t)} \big\Vert^2 
    &\leq \frac{ \xi_t^2 }{n^2}  \Big\Vert i \beta m_{0}^{(t)} +  (1-\beta)  \sum_{j=0}^{i-1}  g_{j}^{(t)}  \Big\Vert ^2 \\
    & \overset{\eqref{update_m_0^{T}}}{=}\frac{ \xi_t^2 }{n^2}  \Big\Vert \beta \frac{i}{n}  \sum_{j=0}^{n-1} g_{j}^{(t-1)} + (1-\beta)  \sum_{j=0}^{i-1}  g_{j}^{(t)} \Big\Vert ^2 \\
    & \overset{(a)}{\leq} \frac{\xi_t^2}{n^2} \left( \beta \Big\Vert  \frac{i}{n}  \sum_{j=0}^{n-1} g_{j}^{(t-1)}\Big\Vert ^2 + (1-\beta) \Big\Vert \sum_{j=0}^{i-1}  g_{j}^{(t)} \Big\Vert ^2 \right) \\
    & \overset{\eqref{define_A}}{=} \frac{\xi_t^2}{n^2} \left( \beta \frac{i^2}{n^2} A_n^{(t-1)} + (1-\beta) A_i^{(t)} \right),
\end{align*}
where (\textit{a}) follows from the convexity of $\norms{\cdot}^2$ and $0\leq \beta < 1$.

For $t=1$ and $i = 0, 1, \dots, n$, we have $m_0^{(t)} = \textbf{0}$, which leads to
\begin{align*}
    \big\Vert w_{i}^{(t)} - w_{0}^{(t)} \big\Vert^2 
    &\leq \frac{ \xi_t^2 }{n^2}  \Big\Vert i \beta m_{0}^{(t)} +  (1-\beta)  \sum_{j=0}^{i-1}  g_{j}^{(t)}  \Big\Vert ^2 
    = \frac{ \xi_t^2 }{n^2}  \Big\Vert (1-\beta)  \sum_{j=0}^{i-1}  g_{j}^{(t)} \Big\Vert ^2 
    =  \frac{\xi_t^2}{n^2}  (1-\beta)^2 A_i^{(t)}.
\end{align*}
Combining these two cases, we obtain \eqref{update_distance_t}.

(\textbf{b}) 
Using similar argument for the second equation \eqref{update_w_i^{T}-1}, for $t \geq 2$ and $i = 0, 1, \dots, n$; we can derive that
\begin{align*}
    w_{i}^{(t-1)} - w_{0}^{(t)} &= \frac{\eta_{t-1}}{n} \sum_{j = i}^{n-1} m_{j+1}^{(t-1)} 
    \overset{\eqref{eq:main_update}}{=}\frac{\eta_{t-1}}{n} \sum_{j = i}^{n-1}\left[ \beta m_{0}^{(t-1)} + (1-\beta) g_{j}^{(t-1)} \right]\\
    &= \frac{\eta_{t-1}}{n} \left[ (n-i)\beta m_{0}^{(t-1)} + (1-\beta) \sum_{j = i}^{n-1} g_{j}^{(t-1)} \right].
\end{align*}
Note again that $\eta_{t-1} \leq \xi_t$, for $t > 2$ and $i = 0, 1, \dots, n$; we can show that
\begin{align*}
    \big\Vert w_{i}^{(t-1)} - w_{0}^{(t)} \big\Vert^2 
    &\leq \frac{\xi_t^2}{n^2}  \Big\Vert (n-i)\beta m_{0}^{(t-1)} + (1-\beta) \sum_{j = i}^{n-1} g_{j}^{(t-1)}\Big\Vert^2 \\
    &\overset{\eqref{update_m_0^{T}}}{=} \frac{\xi_t^2}{n^2}  \Big\Vert \beta \frac{n-i}{n} \sum_{j=0}^{n-1} g_{j}^{(t-2)} + (1-\beta) \sum_{j = i}^{n-1} g_{j}^{(t-1)}\Big\Vert^2 \\
    &\overset{(b)}{\leq} \frac{\xi_t^2}{n^2}  \left[ \beta \Big\Vert  \frac{n-i}{n} \sum_{j=0}^{n-1} g_{j}^{(t-2)}\Big\Vert^2 + (1-\beta) \Big\Vert \sum_{j = i}^{n-1} g_{j}^{(t-1)}\Big\Vert^2 \right] \\
    & \overset{\eqref{define_A}, \eqref{define_B}}{=} \frac{\xi_t^2}{n^2}  \left[ \beta   \frac{(n-i)^2}{n^2} A_n^{(t-2)} + (1-\beta) B_i^{(t-1)} \right],
\end{align*} 
where in (\textit{b}) we use again the convexity of $\norms{\cdot}^2$ and $0\leq \beta < 1$.

For $t=2$, we easily get
\begin{align*}
    \big\Vert w_{i}^{(t-1)} - w_{0}^{(t)} \big\Vert^2 
    &\leq \frac{\xi_t^2}{n^2}  \Big\Vert (n-i)\beta m_{0}^{(t-1)} + (1-\beta) \sum_{j = i}^{n-1} g_{j}^{(t-1)}\Big\Vert^2 
    \leq \frac{\xi_t^2}{n^2}  \Big\Vert (1-\beta) \sum_{j = i}^{n-1} g_{j}^{(t-1)}\Big\Vert^2 
    = \frac{\xi_t^2}{n^2}   (1-\beta)^2 B_i^{(t-1)}.
\end{align*} 
Combining the two cases above, we obtain \eqref{update_distance_t-1}.
\Eproof

\subsubsection{Proof of Lemma~\ref{lem_momentum_variance_02}: Upper Bounding  The Terms $A_n^{(t)}$, $\sum_{i=0}^{n-1}A_i^{(t)}$, and $\sum_{i=0}^{n-1}B^{(t)}_i$}
(\textbf{a}) Let first upper bound the term $A_i^{(t)}$ defined by \eqref{define_A}.
Indeed, for $i = 0, \dots, n-1$ and $t \geq 1$, we have
\allowdisplaybreaks
\begin{align*}
    A_i^{(t)} &= \Big\Vert  \sum_{j=0}^{i-1}  g_{j}^{(t)} \Big\Vert ^2 \\
    & \overset{(a)}{\leq} 3\Big\Vert   \sum_{j=0}^{i-1}  \Big(g_{j}^{(t)} - \nabla f ( w_{0}^{(t)} ; \pi^{(t)} ( j + 1 ) ) \Big) \Big\Vert ^2  
    + 3\Big\Vert  \sum_{j=0}^{i-1}  \Big(\nabla f ( w_{0}^{(t)} ; \pi^{(t)} ( j + 1 ) ) - \nabla F(w_{0}^{(t)}) \Big) \Big\Vert ^2 
    + 3\Big\Vert   \sum_{j=0}^{i-1}   \nabla F(w_{0}^{(t)})  \Big\Vert ^2 \\
    &\overset{(a)}{\leq} 
    3 i  \sum_{j=0}^{i-1} \Big\Vert g_{j}^{(t)} - \nabla f ( w_{0}^{(t)} ; \pi^{(t)} ( j + 1 ) )  \Big\Vert ^2 
    + 3  i  \sum_{j=0}^{i-1} \Big\Vert  \nabla f ( w_{0}^{(t)} ; \pi^{(t)} ( j + 1 ) ) - \nabla F(w_{0}^{(t)})  \Big\Vert ^2 
    + 3 i^2 \big\Vert  \nabla F(w_{0}^{(t)})  \big\Vert ^2 \\
    &\leq 
    3 i  \sum_{j=0}^{n-1} \Big\Vert g_{j}^{(t)} - \nabla f ( w_{0}^{(t)} ; \pi^{(t)} ( j + 1 ) )  \Big\Vert ^2 
    + 3  i \sum_{j=0}^{n-1} \Big\Vert  \nabla f ( w_{0}^{(t)} ; \pi^{(t)} ( j + 1 ) ) - \nabla F(w_{0}^{(t)})  \Big\Vert ^2 
    + 3i^2 \big\Vert  \nabla F(w_{0}^{(t)})  \big\Vert ^2 \\
    & \overset{(b)}{\leq} 3 i  I_t
    + 3  i  n \Big[ \Theta \big\Vert \nabla  F(w_{0}^{(t)})  \big\Vert ^2 + \sigma^2 \Big]
    + 3 i^2\big\Vert \nabla  F(w_{0}^{(t)})  \big\Vert ^2 \\
    &\leq 3 i  I_t + 3 (i^2 + in \Theta) \big\Vert \nabla  F(w_{0}^{(t)})  \big\Vert ^2  +3in \sigma^2,
\end{align*}
where we use the Cauchy-Schwarz inequality in (\textit{a}), and (\textit{b}) follows from Assumption~\ref{as:A1}(c).


Letting $i = n$ in the above estimate, we obtain the first estimate of Lemma \ref{lem_momentum_variance_02}, i.e.:
\begin{align*}
    A_n^{(t)} &\leq 3 n  I_t + 3 n^2 (\Theta+1) \big\Vert \nabla  F(w_{0}^{(t)})  \big\Vert ^2  +3n^2 \sigma^2 \leq 3 n  \Big( I_t + K_{t} \Big), \ t \geq 1. 
\end{align*}
Now we are ready to calculate the sum $\sum_{i=0}^{n-1} A_i^{(t)}$ as
\begin{align*}
    \sum_{i=0}^{n-1} A_i^{(t)}  &\leq 3I_t \sum_{i=0}^{n-1} i + 3  \sum_{i=0}^{n-1}(i^2 + in \Theta) \big\Vert \nabla  F(w_{0}^{(t)})  \big\Vert ^2  +3n \sigma^2 \sum_{i=0}^{n-1} i\\
    &\leq \frac{3n^2}{2} I_t +    \frac{3 n^3}{2} (\Theta + 1) \big\Vert \nabla  F(w_{0}^{(t)})  \big\Vert ^2  + \frac{3n^3 \sigma^2}{2}
    \leq \frac{3n^2}{2} \Big( I_t + K_t\Big), \  t \geq 1, 
\end{align*}
where the last line follows from $\sum_{i=0}^{n-1} i \leq \frac{n^2}{2}$ and $\sum_{i=0}^{n-1} i^2 \leq \frac{n^3}{2}$. 
This proves the second estimate of Lemma \ref{lem_momentum_variance_02}.

(\textbf{b})
Using similar argument as above, for $i =0, 1, \dots, n-1$ and $t \geq 1$, we can derive
\begin{align*}
    B_i^{(t)} &= \Big\Vert \sum_{j = i}^{n-1} g_{j}^{(t)}\Big\Vert^2 
    \overset{(a)}{\leq}
    3\Big\Vert   \sum_{j = i}^{n-1}  \Big(g_{j}^{(t)} - \nabla f ( w_{0}^{(t)} ; \pi^{(t)} ( j + 1 ) ) \Big) \Big\Vert ^2 \\
    &\quad + 3\Big\Vert  \sum_{j=i}^{n-1}  \Big(\nabla f ( w_{0}^{(t)} ; \pi^{(t)} ( j + 1 ) ) - \nabla F(w_{0}^{(t)}) \Big) \Big\Vert ^2 
    + 3\Big\Vert   \sum_{j=i}^{n-1}   \nabla F(w_{0}^{(t)})  \Big\Vert ^2 \\
    &\overset{(a)}{\leq}
    3 (n-i)  \sum_{j=i}^{n-1} \Big\Vert g_{j}^{(t)} - \nabla f ( w_{0}^{(t)} ; \pi^{(t)} ( j + 1 ) )  \Big\Vert ^2 \\
    &\quad + 3  (n-i)  \sum_{j=i}^{n-1} \Big\Vert  \nabla f ( w_{0}^{(t)} ; \pi^{(t)} ( j + 1 ) ) - \nabla F(w_{0}^{(t)})  \Big\Vert ^2 
    + 3 (n-i)^2 \Big\Vert \nabla  F(w_{0}^{(t)})  \Big\Vert ^2 \\
    & \overset{\eqref{define_I}}{\leq} 
    3 (n-i)  I_t 
    + 3 (n-i) \sum_{j=0}^{n-1} \Big\Vert  \nabla f ( w_{0}^{(t)} ; \pi^{(t)} ( j + 1 ) ) -\nabla  F(w_{0}^{(t)})  \Big\Vert ^2 
    + 3(n-i)n \big\Vert \nabla  F(w_{0}^{(t)})  \big\Vert ^2 \\
    &\overset{(b)}{\leq} 3 (n-i)  I_t
    + 3  (n-i)  n \Big[ \Theta \big\Vert \nabla  F(w_{0}^{(t)})  \big\Vert ^2 + \sigma^2 \Big]
    + 3 (n -i)n\big\Vert \nabla  F(w_{0}^{(t)})  \big\Vert ^2 \\
    &\leq 3 (n-i)  I_t + 3 (n-i)n(\Theta+1) \big\Vert \nabla  F(w_{0}^{(t)})  \big\Vert ^2  +3(n-i)n \sigma^2,
\end{align*}
where we use again the Cauchy-Schwarz inequality in (\textit{a}), and (\textit{b}) follows from Assumption~\ref{as:A1}(c). Finally, for all $t\geq 1$, we calculate the sum $\sum_{i=0}^{n-1} B_i^{(t)}$ as follows:
\begin{align*}
    \sum_{i=0}^{n-1} B_i^{(t)}  &\leq 3I_t \sum_{i=0}^{n-1} (n-i) + 3  \sum_{i=0}^{n-1}(n-i)n(\Theta+1) \Big\Vert \nabla  F(w_{0}^{(t)})  \Big\Vert ^2  +3n \sigma^2 \sum_{i=0}^{n-1} (n-i)\\
    &\leq 3n^2I_t + 3  n^3(\Theta+1) \Big\Vert \nabla  F(w_{0}^{(t)})  \Big\Vert ^2  +3n^3 \sigma^2 \\
    & \leq 3n^2 \Big(I_t + K_t \Big),
\end{align*}
where the last line follows from the fact that $\sum_{i=0}^{n-1} (n-i) \leq n^2$.
This proves \eqref{bound_B^{T}}.
\Eproof

\subsubsection{Proof of Lemma~\ref{lem_momentum_variance_03}: Upper Bounding The Terms $I_t$, $J_t$, and $I_t + K_t$}
(\textbf{a})
First, for every $t \geq 1$ and $g_{i}^{(t)} := \nabla f ( w_{i}^{(t)} ; \pi^{(t)} ( i + 1 ) )$, we have 
\begin{align*}
    I_t &= \sum_{i=0}^{n-1} \Big \Vert \nabla f ( w_{0}^{(t)} ; \pi^{(t)} ( i + 1 ) )  - g_{i}^{(t)} \Big\Vert^2 \\
    &\overset{\tiny \textrm{Step~\ref{alg:A1_step3}}}{=} \sum_{i=0}^{n-1} \Big \Vert \nabla f ( w_{0}^{(t)} ; \pi^{(t)} ( i + 1 ) )  - \nabla f ( w_{i}^{(t)} ; \pi^{(t)} ( i + 1 ) )  \Big\Vert^2 \\
    &\leq L^2 \sum_{i=0}^{n-1} \big\Vert w_{i}^{(t)} - w_{0}^{(t)} \big\Vert^2, 
    \tagthis \label{bound_I_t_02}
\end{align*}
where the last inequality follows from Assumption~\ref{as:A1}(b).

Applying Lemma \ref{lem_momentum_variance_01} to \eqref{bound_I_t_02}, for $t> 1$, we have
\begin{align*}
    I_t &  \leq \frac{L^2\xi_t^2}{n^2} \sum_{i=0}^{n-1}\left[ \beta \frac{i^2}{n^2} A_n^{(t-1)} + (1-\beta) A_i^{(t)} \right] \\
    &\leq  \frac{L^2\xi_t^2}{n^2}  \left[ \beta \frac{ A_n^{(t-1)} }{n^2}  \sum_{i=0}^{n-1} i^2 +  (1-\beta) \sum_{i=0}^{n-1} A_i^{(t)} \right]\\
    &\overset{(a)}{\leq}  \frac{L^2\xi_t^2}{n^2} \left[ \beta \frac{n}{3} 3 n  \Big( I_{t-1} + K_{t-1} \Big)  +  (1-\beta)  \frac{3n^2}{2} \Big( I_t + K_t\Big) \right]\\
    &=  L^2\xi_t^2 \left[ \beta \big(  I_{t-1} + K_{t-1}\big) +   \frac{3}{2} (1-\beta) \Big(I_t + K_t \Big) \right],
\end{align*}
where (\textit{a}) follows from the result \eqref{bound_A^{T}} in Lemma \ref{lem_momentum_variance_02} and the fact that $\sum_{i=0}^{n-1} i^2 \leq \frac{n^3}{3}$.

For $t = 1$, by applying Lemma \ref{lem_momentum_variance_01} and \ref{lem_momentum_variance_02} consecutively to \eqref{bound_I_t_02} we have
\begin{align*}
    I_t & \leq L^2 
    \sum_{i=0}^{n-1} \big\Vert w_{i}^{(t)} - w_{0}^{(t)} \big\Vert^2 \overset{\eqref{update_distance_t}}{\leq} \frac{L^2\xi_t^2}{n^2} \sum_{i=0}^{n-1}\left( (1-\beta)^2 A_i^{(t)} \right) 
    \leq  \frac{L^2\xi_t^2 (1-\beta)^2}{n^2} \sum_{i=0}^{n-1} A_i^{(t)} \\
    &\overset{\eqref{bound_A^{T}}}{\leq}  \frac{L^2\xi_t^2 (1-\beta)^2}{n^2} \frac{3n^2}{2} \Big( I_t + K_t\Big) 
    = \frac{3}{2} L^2\xi_t^2  (1-\beta)^2 \Big(I_t + K_t \Big). \tagthis \label{bound_I=1}
\end{align*}

(\textbf{b})
Using the similar argument as above we can derive that
\begin{align*}
    J_t &= \sum_{i=0}^{n-1} \Big \Vert \nabla f ( w_{0}^{(t)} ; \pi^{(t-1)} ( i + 1 ) )  - g_{i}^{(t-1)} \Big\Vert^2 \\
    &\overset{\tiny\textrm{Step~\ref{alg:A1_step3}}}{=} \sum_{i=0}^{n-1} \Big \Vert \nabla f ( w_{0}^{(t)} ; \pi^{(t-1)} ( i + 1 ) )  - \nabla f ( w_{i}^{(t-1)} ; \pi^{(t-1)} ( i + 1 ) )  \Big\Vert^2 \\
    & \overset{\eqref{eq:Lsmooth_basic}}{\leq} L^2 
    \sum_{i=0}^{n-1} \big\Vert w_{i}^{(t-1)} - w_{0}^{(t)} \big\Vert^2.
\end{align*}
Applying Lemma \ref{lem_momentum_variance_01} for $t > 2$ we get to the above estimate, we get
\begin{align*}
    J_t &\leq \frac{L^2\xi_t^2}{n^2} \sum_{i=0}^{n-1}  \left[ \beta   \frac{(n-i)^2}{n^2} A_n^{(t-2)} + (1-\beta) B_i^{(t-1)} \right] \\
    &\leq \frac{L^2\xi_t^2}{n^2} \left[ \beta \frac{ A_n^{(t-2)} }{n^2}  \sum_{i=0}^{n-1} (n-i)^2 +  (1-\beta) \sum_{i=0}^{n-1} B_i^{(t-1)} \right] \\
    &\overset{(b)}{\leq}  \frac{L^2\xi_t^2}{n^2} \left[ \beta n \cdot 3n  \Big( I_{t-2} + K_{t-2} \Big)  + (1-\beta) \cdot 3n^2 \Big(I_{t-1} + K_{t-1}\Big) \right]\\
    & =  3L^2\xi_t^2 \left[ \beta  \Big(  I_{t-2} + K_{t-2} \Big) + (1-\beta)  \Big(I_{t-1} + K_{t-1} \Big)\right],
\end{align*}
where (\textit{b}) follows from the results \eqref{bound_A^{T}}, \eqref{bound_B^{T}} in Lemma \ref{lem_momentum_variance_02} and the fact that $\sum_{i=0}^{n-1} (n-i)^2 \leq n^3$.

Now applying Lemma \ref{lem_momentum_variance_01} for the special case $t = 2$, we obtain
\begin{align*}
    J_t &\leq \frac{L^2\xi_t^2}{n^2} \sum_{i=0}^{n-1} (1-\beta)^2 B_i^{(t-1)} 
    \leq \frac{L^2\xi_t^2}{n^2}  (1-\beta)^2 \cdot 3n^2 \Big(I_{t-1} + K_{t-1}\Big)
    \leq  3L^2\xi_t^2 \Big(I_{t-1} + K_{t-1} \Big),
\end{align*}
where we use the result \eqref{bound_B^{T}} from Lemma \ref{lem_momentum_variance_02} and the fact that $1 - \beta \leq 1$. 

(\textbf{c}) 
First, from the result of Part \textbf{(a)}, for $t > 1$, we have
\begin{align*}
    &I_t \leq L^2\xi_t^2 \left[ \beta \big(  I_{t-1} + K_{t-1}\big) +   \frac{3}{2} (1-\beta) K_t \right] + \frac{3}{2} (1-\beta) L^2\xi_t^2 I_t.
\end{align*}
Since $L^2\xi_t^2 \leq \frac{2}{5}$ (due to the choice of our learning rate) and $1-\beta \leq 1$, we further get
\begin{align*}
    I_t &\leq   \frac{2}{5} \left[ \beta \Big(  I_{t-1} + K_{t-1}\Big) + \frac{3}{2} (1-\beta) K_t \right] + \frac{3}{5} I_t,
\end{align*}
which is equivalent to 
\begin{align*}
    I_t \leq \beta \Big(  I_{t-1} + K_{t-1}\Big) +   \frac{3}{2} (1-\beta) K_t.
\end{align*}
Adding $K_t$ to both sides of this inequality, for $t > 1$, we get
\begin{align*}
    I_t + K_t \leq \beta \Big(  I_{t-1} + K_{t-1}\Big) +   \frac{5 -3\beta}{2} K_t. \tagthis \label{bound_sum_I_K>1}
\end{align*}
Finally, for $t = 1$, since $L^2\xi_t^2 \leq \frac{2}{5}$ and $1-\beta \leq 1$, we have
\begin{align*}
    I_t & \leq \frac{3}{2} L^2\xi_t^2  (1-\beta)^2 \Big(I_t + K_t \Big) 
    \leq \frac{3}{5} (1-\beta)^2 I_t + \frac{3}{5} (1-\beta)^2 K_t \leq \frac{3}{5} I_t + \frac{3}{5} (1-\beta) K_t,
\end{align*}
which is equivalent to 
\begin{align*}
    \frac{2}{5} I_t \leq \frac{3}{5} (1-\beta) K_t.
    \tagthis \label{bound_I_1}
\end{align*}
This leads to our desired result for $t=1$ as
\begin{align*}
    I_1 +K_1 \leq \frac{3}{2} (1-\beta) K_1 + K_1 = \frac{5 - 3\beta}{2} K_1. \tagthis \label{bound_sum_I_K=1}
\end{align*}
Combining two cases, we obtain the desired result in Part (\textbf{c}).
\Eproof

\subsubsection{Proof of Lemma~\ref{lem_momentum_variance_04}: Upper Bounding The Key Quantity $\beta J_t + (1-\beta) I_t$}
First, we analyze the case $t>2$ using the results of Lemma \ref{lem_momentum_variance_03} as follows:
\begin{align*}
    \beta J_t + (1-\beta) I_t & \leq 
    3\beta L^2\xi_t^2 \left[ \beta  \Big(  I_{t-2} + K_{t-2} \Big) + (1-\beta)  \Big(I_{t-1} + K_{t-1} \Big)\right] \\ 
    & \quad + (1-\beta) L^2\xi_t^2 \left[ \beta \big(  I_{t-1} + K_{t-1}\big) +   \frac{3}{2} (1-\beta) \Big(I_t + K_t \Big) \right] \\
    & \leq L^2\xi_t^2 \left[  2 \big(I_t + K_t \big)  + 4\beta \big(I_{t-1} + K_{t-1} \big)
    + 3\beta^2  \big(I_{t-2} + K_{t-2} \big) \right] \\
    & =: L^2\xi_t^2 P_t, \tagthis \label{define_P}
\end{align*}
where the last two lines follow since $1 - \beta \leq 1$ and $P_t := 2(I_t + K_t ) + 4\beta (I_{t-1} + K_{t-1} ) + 3\beta^2 (I_{t-2} + K_{t-2})$ for $t>2$.

Next, we bound the term $P_t$ in \eqref{define_P} as
\begin{align*}
    P_t &= 2 \big(I_t + K_t \big) 
    + 4\beta \big(I_{t-1} + K_{t-1} \big)
    + 3\beta^2  \big(I_{t-2} + K_{t-2} \big)\\
    &\leq 2  \Big[ \beta \Big(I_{t-1} + K_{t-1}\Big) + \frac{5 -3\beta}{2} K_t \Big]
    + 4\beta \Big(I_{t-1} + K_{t-1} \Big)
    + 3\beta^2  \Big(I_{t-2} + K_{t-2} \Big) &\text{apply \eqref{bound_sum_I_K>1} for} \ t\\
    &= (5 -3\beta) K_t  
    + 6\beta \Big(I_{t-1} + K_{t-1} \Big)
    + 3\beta^2  \Big(I_{t-2} + K_{t-2} \Big) \\
    &\leq (5 -3\beta) K_t + 6\beta \Big[ \beta \big(I_{t-2} + K_{t-2}\big) + \frac{5 -3\beta}{2} K_{t-1} \Big]
    + 3\beta^2  \Big(I_{t-2} + K_{t-2} \Big) &\text{apply \eqref{bound_sum_I_K>1} for} \ t-1\\
    &= (5 -3\beta) K_t + 3 \beta (5 -3\beta)K_{t-1}
    + 9\beta^2  \Big(I_{t-2} + K_{t-2} \Big).
\end{align*}
We consider the last term, which can be bounded as
\begin{align*}
    \beta^2  \Big(I_{t-2} + K_{t-2} \Big)
    &= \beta^3 \big(I_{t-3} + K_{t-3}\big) + \frac{5 -3\beta}{2} \beta^2 K_{t-2} &\text{apply \eqref{bound_sum_I_K>1} recursively for} \ t-2 > 1\\
    &= \beta^{t-1} \big(I_1 + K_1\big) +\frac{5 -3\beta}{2}   \sum_{j=2}^{t-2}  \beta^{t-j} K_j  \\
    &\leq \beta^{t-1} \frac{5 -3\beta}{2} K_1 +\frac{5 -3\beta}{2}   \sum_{j=2}^{t-2}  \beta^{t-j} K_j 
    &\text{apply \eqref{bound_sum_I_K=1}} \\
    &\leq \frac{5 -3\beta}{2}   \sum_{j=1}^{t-2}  \beta^{t-j} K_j.
\end{align*}

Note that this bound is also true for the case $t-2 = 1$:
\begin{align*}
    \beta^2  \Big(I_{t-2} + K_{t-2} \Big) &= \beta^2  \Big(I_{1} + K_{1} \Big) \leq \beta^2  \frac{5 -3\beta}{2} K_1  &\text{apply \eqref{bound_sum_I_K=1}} 
\end{align*}

Substituting this inequality into $P_t$, for $t > 2$, we get
\begin{align*}
    P_t &\leq (5 -3\beta) K_t + 3 \beta (5 -3\beta)K_{t-1}
    + 9\beta^2  \Big(I_{t-2} + K_{t-2} \Big) \\
    &\leq \frac{9}{2}(5 -3\beta)  K_t + \frac{9}{2} (5 -3\beta)  \beta K_{t-1}
    +  \frac{9}{2} (5 -3\beta)  \sum_{j=1}^{t-2}  \beta^{t-j} K_j \\
    & \leq \frac{9}{2} (5 -3\beta)   \sum_{j=1}^{t}  \beta^{t-j} K_j.
\end{align*}
Now we analyze similarly for the case $t=2$ as follows:
\begin{align*}
    \beta J_2 + (1-\beta) I_2 & \leq 
    3\beta L^2\xi_2^2  \Big(I_{1} + K_{1} \Big)
    + (1-\beta) L^2\xi_2^2 \left[ \beta \big(  I_{1} + K_{1}\big) +   \frac{3}{2} (1-\beta) \Big(I_2 + K_2 \Big) \right]\\
    & \leq L^2\xi_2^2 \left[  2 \Big(I_2 + K_2 \Big)  + 4\beta \Big(I_{1} + K_{1} \Big) \right] \\
    & =: L^2\xi_2^2 P_2, 
\end{align*}
where the last line follows since $1 - \beta \leq 1$ and $P_2 := 2 \Big(I_2 + K_2 \Big) + 4\beta \Big(I_{1} + K_{1} \Big)$. 

Next, we bound the term $P_2$ as follows:
\begin{align*}
    P_2 &= 2 \Big(I_2 + K_2 \Big)  + 4\beta \Big(I_1 + K_1 \Big) \\
    & \leq 2  \left(\beta \Big(I_1 + K_1\Big) + \frac{5 -3\beta}{2} K_2 \right)  + 4\beta \Big(I_1 + K_1 \Big) &\text{apply \eqref{bound_sum_I_K>1} for} \ t=2\\
    &= (5 -3\beta) K_2   + 6\beta \Big(I_1 + K_1 \Big) \\
    & \leq (5 -3\beta) K_2 + 6\beta  \frac{5 -3\beta}{2} K_1  &\text{apply \eqref{bound_sum_I_K=1}}\\
    & = (5 -3\beta) K_2 + 3 \beta (5 -3\beta)K_1 \\
    & \leq \frac{9}{2} (5 -3\beta)   \sum_{j=1}^{2}  \beta^{2-j} K_j.
\end{align*}
Hence, the statements $\beta J_t + (1-\beta)I_t \leq L^2\xi_t^2 P_t$ and $P_t \leq \frac{9}{2} (5 -3\beta)   \sum_{j=1}^{t}  \beta^{t-j} K_j $ are true for every $t \geq 2$. 

Combining these two cases, we have the following estimates:
\begin{align*}
    \beta J_t + (1-\beta)I_t \leq L^2\xi_t^2 P_t
    &\leq L^2\xi_t^2 \cdot \frac{9}{2} (5 -3\beta)   \sum_{j=1}^{t}  \beta^{t-j} K_j\\
    &\overset{\eqref{define_K}}{\leq} L^2\xi_t^2 \cdot \frac{9}{2} (5 -3\beta)   \sum_{j=1}^{t}  \beta^{t-j} \Big[ n (\Theta+1) \big\Vert \nabla  F(w_{0}^{(j)}) \big\Vert ^2  + n \sigma^2 \Big] \\
    &\leq L^2\xi_t^2 \cdot \frac{9}{2} (5 -3\beta)\Big[ n (\Theta+1) \sum_{j=1}^{t}  \beta^{t-j} \big\Vert \nabla F(w_{0}^{(j)}) \big\Vert ^2  +  \sum_{j=1}^{t}  \beta^{t-j} n \sigma^2 \Big] \\
    &\leq \frac{9(5 -3\beta)L^2\xi_t^2}{2} \Big[ n (\Theta+1) \sum_{j=1}^{t}  \beta^{t-j} \big\Vert \nabla F(w_{0}^{(j)}) \big\Vert ^2  +  \frac{n \sigma^2}{1-\beta}\Big],  
\end{align*}
where the last line follows since $\sum_{j=1}^{t}  \beta^{t-j} \leq \frac{1}{1-\beta}$ for every $t\geq2$.
Hence, we have proved \eqref{bound_beta_IJ}.
\Eproof

\subsubsection{The Proof of Lemma~\ref{lem_momentum_variance_05}}
Using the assumption that $\eta_1 \geq \eta_2 \geq \dots \geq \eta_T$, we can derive the following sum as
\begin{align*}
    \sum_{t=1}^{T} \eta_t \sum_{j=1}^{t}  \beta^{t-j} \big\Vert \nabla F(\tilde{w}_{j-1}) \big\Vert ^2
    &= \sum_{t=1}^{T} \eta_t \left[\beta^{t-1} \big\Vert \nabla F(\tilde{w}_{0}) \big\Vert ^2 + \beta^{t-2} \big\Vert \nabla F(\tilde{w}_{1}) \big\Vert ^2 + \dots  +\beta^{0} \big\Vert \nabla F(\tilde{w}_{t-1}) \big\Vert ^2 \right] \\
    &= \quad\eta_1 \beta^{0} \big\Vert \nabla F(\tilde{w}_{0}) \big\Vert ^2  \\
    & \quad + \eta_2 \beta^{1} \big\Vert \nabla F(\tilde{w}_{0}) \big\Vert ^2 
    \ \ + \eta_2\beta^{0} \big\Vert \nabla F(\tilde{w}_1) \big\Vert ^2\\
    &\quad+ \eta_t \beta^{t-1} \big\Vert \nabla F(\tilde{w}_{0}) \big\Vert ^2 
    \ + \eta_t \beta^{t-2} \big\Vert \nabla F(\tilde{w}_1) \big\Vert ^2 + \dots  
    +\eta_t\beta^{0} \big\Vert \nabla F(\tilde{w}_{t-1}) \big\Vert ^2 \\
    &\quad+ \eta_T \beta^{T-1} \big\Vert \nabla F(\tilde{w}_{0}) \big\Vert ^2 
    + \eta_T \beta^{T-2} \big\Vert \nabla F(\tilde{w}_1) \big\Vert ^2 \  + \  \dots  
    \  + \ \eta_T\beta^{0} \big\Vert \nabla F(\tilde{w}_{T-1}) \big\Vert ^2 \\
    &\leq \eta_1 \beta^{0} \big\Vert \nabla F(\tilde{w}_{0}) \big\Vert ^2  \\
    &\quad+ \eta_1 \beta^{1} \big\Vert \nabla F(\tilde{w}_{0}) \big\Vert ^2 
    + \eta_2\beta^{0} \big\Vert \nabla F(\tilde{w}_1) \big\Vert ^2\\
    &\quad+ \eta_1 \beta^{t-1} \big\Vert \nabla F(\tilde{w}_{0}) \big\Vert ^2 
    \ \ + \eta_2 \beta^{t-2} \big\Vert \nabla F(\tilde{w}_1) \big\Vert ^2 + \dots  
    +\eta_t\beta^{0} \big\Vert \nabla F(\tilde{w}_{t-1}) \big\Vert ^2 \\
    &\quad + \eta_1 \beta^{T-1} \big\Vert \nabla F(\tilde{w}_{0}) \big\Vert ^2 
    \ + \eta_2 \beta^{T-2} \big\Vert \nabla F(\tilde{w}_1) \big\Vert ^2 \ + \ \dots  
     \  + \ \eta_T\beta^{0} \big\Vert \nabla F(\tilde{w}_{T-1}) \big\Vert ^2 \\
    &\leq \eta_1\sum_{i=0}^{T-1}\beta^i \big\Vert \nabla F(\tilde{w}_{0}) \big\Vert ^2 + \eta_2 \sum_{i=0}^{T-2}\beta^i \big\Vert \nabla F(\tilde{w}_1) \big\Vert ^2 \ + \ \dots \ + \ \eta_T \beta^{0} \big\Vert \nabla F(\tilde{w}_{T-1}) \big\Vert ^2 \\
    &\leq \frac{1}{1-\beta}  \sum_{t=1}^{T} \eta_t \big\Vert \nabla F(\tilde{w}_{t-1}) \big\Vert ^2,
\end{align*}
where the last inquality follows since $\sum_{i=0}^{t-1}\beta^i \leq \frac{1}{1-\beta}$ for $t = 1, 2 \dots, T$. 
\Eproof
\subsubsection{The Proof of Lemma~\ref{lem_momentum_variance_06}}

We analyze the special case $t = 1$ as follows.
First, letting $i=n$ in Equation \eqref{update_w_i^{T}}, for $t\geq 1$, we have
\begin{align*}
    w_{n}^{(t)} - w_{0}^{(t)} &= - \frac{\eta_t}{n} \sum_{j=0}^{n-1} m_{j+1}^{(t)} 
    \overset{\eqref{eq:main_update}}{=}  - \frac{\eta_t}{n} \sum_{j=0}^{n-1} \big[ \beta m_{0}^{(t)} + (1-\beta) g_{j}^{(t)}\big]
    =  - \frac{\eta_t}{n} \Big[ n \beta m_{0}^{(t)} +  (1-\beta)  \sum_{j=0}^{n-1}  g_{j}^{(t)}  \Big].
\end{align*}
For $t = 1$, since $m_{0}^{(t)} =\mathbf{0}$, we have 
\begin{align*}
    w_{n}^{(t)} - w_{0}^{(t)} =- \frac{\eta_t}{n}  (1-\beta)  \sum_{j=0}^{n-1}  g_{j}^{(t)}.\tagthis \label{update_epoch_03}
\end{align*}
Using the $L$-smoothness of $F$ in Assumption \ref{as:A1}(b), $w_0^{(t+1)} := w_{n}^{(t)}$, and \eqref{update_epoch_03}, we have
\begin{align*} 
    F( w_0^{(t+1)} )  & \overset{\eqref{eq:Lsmooth}}{\leq} F( w_0^{(t)} ) + \nabla F( w_0^{(t)} )^{\top}(w_0^{(t+1)} - w_0^{(t)}) + \frac{L}{2}\norms{w_0^{(t+1)} - w_0^{(t)}}^2  \\
    &\overset{\eqref{update_epoch_03}}{=} F( w_0^{(t)} ) - \eta_t (1-\beta) \nabla F( w_0^{(t)} )^{\top} \left( \frac{1}{n} \sum_{j=0}^{n-1}    g_{j}^{(t)}  \right)  
    + \frac{L \eta_t^2}{2} (1-\beta)^2 \Big \Vert \frac{1}{n} \sum_{j=0}^{n-1} g_{j}^{(t)}  \Big\Vert ^2  \\
    & \overset{\tiny(a)}{=} F( w_0^{(t)} ) - \frac{\eta_t}{2}(1-\beta) \norms{ \nabla F( w_0^{(t)} )}^2 
    + \frac{\eta_t}{2} (1-\beta) \Big\Vert \nabla F( w_0^{(t)} )  - \frac{1}{n} \sum_{j=0}^{n-1}   g_{j}^{(t)}\Big\Vert^2  \\
    & \quad - \frac{\eta_t}{2} (1-\beta) \left( 1 - L\eta_t(1-\beta)  \right) \Big\Vert \frac{1}{n} \sum_{j=0}^{n-1}  g_{j}^{(t)}  \Big\Vert^2   \\
    & \leq F( w_0^{(t)} ) - \frac{\eta_t}{2}(1-\beta) \norms{ \nabla F( w_0^{(t)} )}^2 
    + \frac{\eta_t}{2} (1-\beta) \Big\Vert \nabla F( w_0^{(t)} )  - \frac{1}{n} \sum_{j=0}^{n-1}   g_{j}^{(t)}\Big\Vert^2, \tagthis \label{bound_F(w)_1}
\end{align*}
where (\textit{a}) follows from the equality $u^{\top}v = \frac{1}{2}\left(\norms{u}^2 + \norms{v}^2 - \norms{u - v}^2\right)$ and the last equality comes from the fact that $\eta_t \leq \frac{\sqrt{2}}{\sqrt{5} L}\leq \frac{1}{L(1-\beta)}$.



Next, we bound the last term of \eqref{bound_F(w)_1} as follows:
\begin{align*} 
     \Big\Vert \nabla F( w_0^{(1)} )  - \frac{1}{n} \sum_{j=0}^{n-1}  g_{j}^{(1)} \Big\Vert^2
    &=  \Big\Vert \frac{1}{n} \sum_{j=0}^{n-1}  \Big( \nabla f ( w_{0}^{(1)} ; \pi^{(1)} ( j + 1 ) )  - g_{j}^{(1)} \Big) \Big\Vert^2  \\
    & \overset{(b)}{\leq} \frac{1}{n} \sum_{j=0}^{n-1} \Big\Vert 
     \nabla f ( w_{0}^{(1)} ; \pi^{(1)} ( j + 1 ) )  - g_{j}^{(1)} \Big\Vert^2\\
    & \overset{\eqref{define_I}}{=} \frac{1}{n} I_1,  
\end{align*}
where (\textit{b}) is from the Cauchy-Schwarz inequality. 

Applying this into \eqref{bound_F(w)_1} we get the following estimate
\begin{align*} 
    F( w_0^{(2)} )    & \leq F( w_0^{(1)} ) - \frac{\eta_1}{2}(1-\beta) \norms{ \nabla F( w_0^{(1)} )}^2 
    + \frac{\eta_1}{2} (1-\beta) \frac{1}{n} I_1.\tagthis \label{estimate_alg_1_t=1}
\end{align*}
\begin{remark}
    Note that the derived estimate \eqref{estimate_alg_1_t=1} is true for all shuffling strategies, including Random Reshuffling, Shuffle Once and Incremental Gradient (under the assumptions of Theorem \ref{thm_momentum_variance} and \ref{thm_momentum_variance_RR}, respectively).
\end{remark}

Applying Lemma \ref{lem_momentum_variance_03}, we further have
\begin{align*} 
    F( w_0^{(2)} )    & \leq F( w_0^{(1)} ) - \frac{\eta_1}{2}(1-\beta) \norms{ \nabla F( w_0^{(1)} )}^2 
    + \frac{\eta_1}{2} (1-\beta) \frac{1}{n} I_1\\
    & \overset{\eqref{bound_I=1}}{\leq} F( w_0^{(1)} ) - \frac{\eta_1}{2}(1-\beta) \norms{ \nabla F( w_0^{(1)} )}^2 
    + \frac{\eta_1}{2} (1-\beta) \frac{3}{2n} L^2\xi_1^2  (1-\beta)^2 \Big(I_1 + K_1 \Big)\\
    &\overset{\eqref{bound_sum_I_K=1}}{=} F( w_0^{(1)} ) - \frac{\eta_1}{2}(1-\beta) \norms{ \nabla F( w_0^{(1)} )}^2 
    + \frac{\eta_1}{2} (1-\beta) \frac{3}{2n} L^2\xi_1^2  (1-\beta)^2 \frac{5 - 3\beta}{2} K_1.
\end{align*}


Noting that $\tilde{w}_{t-1} = w_0^{(t)}, 1-\beta < 1$ and $K_{t} := n (\Theta +1) \big\Vert \nabla F(w_{0}^{(t)}) \big\Vert ^2  + n \sigma^2$, we get
\begin{align*} 
    F( \tilde{w}_{1} ) \leq F( \tilde{w}_{0} ) - \frac{\eta_1}{2}(1-\beta) \norms{ \nabla F(\tilde{w}_{0} )}^2 
    + \frac{3\eta_1L^2\xi_1^2(1-\beta) (5 - 3\beta)}{8} \left[ (\Theta +1) \big\Vert \nabla F(\tilde{w}_{0}) \big\Vert ^2  + \sigma^2 \right].
\end{align*}

Rearranging the terms and dividing both sides by $(1-\beta)$, we have
\begin{align*} 
    \frac{\eta_1}{2} \norms{ \nabla F(\tilde{w}_{0} )}^2  \leq \frac{ F( \tilde{w}_{0} ) - F( \tilde{w}_{1} )}{(1-\beta)}
    + \frac{3\eta_1L^2\xi_1^2(5 - 3\beta)}{8} (\Theta +1) \big\Vert \nabla F(\tilde{w}_{0}) \big\Vert ^2  +  \frac{3\eta_1\sigma^2L^2\xi_1^2(5 - 3\beta)}{8}.
\end{align*}

Since $\xi_t$ satisfies $9 L^2 \xi_t^2 (5 -3\beta) (\Theta +1) \leq 1-\beta$ as above, we can deduce from the last estimate that
\begin{align*} 
    \frac{\eta_1}{2} \norms{ \nabla F(\tilde{w}_{0} )}^2  \leq \frac{ F( \tilde{w}_{0} ) - F( \tilde{w}_{1} )}{(1-\beta)}
    + \frac{(1-\beta)\eta_1 }{4} \big\Vert \nabla F(\tilde{w}_{0}) \big\Vert ^2  +  \frac{9\sigma^2L^2(5 - 3\beta)}{4(1-\beta)}\cdot \xi_1^3,
\end{align*}

which proves \eqref{bound_F(w)_1_tilde}.
\Eproof

\section{Convergence Analysis of Algorithm~\ref{sgd_momentum_shuffling_01} for Randomized Reshuffling Strategy}
\label{apdx:sec:convergence_analysis_RR}
In this section, we present the convergence results of Algorithm~\ref{sgd_momentum_shuffling_01} for Randomized Reshuffling strategy. 
Since this analysis is similar to the previous bounds in Theorem~\ref{thm_momentum_variance}, we will highlight the similarities and discuss the differences between the two analyses. 

\subsection{The Proof Sketch of Theorem~\ref{thm_momentum_variance_RR}}\label{subsec:proof_TH1_RR}
The proof of Theorem~\ref{thm_momentum_variance_RR} is similar to Theorem~\ref{thm_momentum_variance}, except for the fact that $\E[\beta I_t + (1-\beta)J_t]$ is upper bounded by a different term.
The proof of Theorem~\ref{thm_momentum_variance_RR} is divided into the following steps.
\begin{compactitem}
\item From Theorem~\ref{thm_momentum_variance} we have 
\begin{equation*}
F(\tilde{w}_{t}) \leq  F(\tilde{w}_{t-1}) - \frac{\eta_t}{2}\norms{\nabla{F}(\tilde{w}_{t-1})}^2 + \frac{\eta_t}{2n}(\beta I_t + (1-\beta)J_t) , 
\end{equation*}
where $I_t$ and $J_t$ are defined by \eqref{define_I} and \eqref{define_J}.

\item Then, $\E[\beta I_t + (1-\beta)J_t]$ can be upper bounded as
\begin{equation*}
\E [\beta I_t + (1-\beta)J_t ]\leq \Ocal\Big(\xi_t^2\sum_{j=0}^t\beta^{t-j} \E \left[\norms{\nabla{F}(w_0^{(j)})}^2 \right]+ \xi^2_t\sigma^2 \Big),
\end{equation*}
as shown in Lemma~\ref{lem_momentum_variance_04_RR}, where $\xi_t := \max\{\eta_{t-1}, \eta_t\}$ defined in \eqref{define_xi}.
\item We further upper bound the sum of the right-hand side of the last estimate as in Lemma~\ref{lem_momentum_variance_05} in terms $\sum_{t=0}^T\eta_t \E \left[\norms{\nabla{F}(\tilde{w}_{t-1})}^2\right]$.
\end{compactitem}
Combining these steps together, and using some simplification, we obtain \eqref{eq_thm_momentum_variance_RR} in Theorem~\ref{thm_momentum_variance_RR}.

\subsection{Technical Lemmas}

In this subsection, we will introduce some  intermediate results used in Theorem~\ref{thm_momentum_variance_RR}. Lemmas \ref{lem_momentum_variance_02_RR} and \ref{lem_momentum_variance_03_RR} in Theorem~\ref{thm_momentum_variance_RR} play a similar role as previous Lemmas \ref{lem_momentum_variance_02} and \ref{lem_momentum_variance_03} in Theorem~\ref{thm_momentum_variance}.

\begin{lem} \label{lem_momentum_variance_02_RR}
Assume that the Randomized Reshuffling strategy is used in Algorithm \ref{sgd_momentum_shuffling_01}. 
Then, under the same setting as of Lemma~\ref{lem_momentum_variance_01} and Assumption~\ref{as:A1}(c), for $t\geq 1$, it holds that
\begin{align*}
    (\textbf{a}) \qquad 
    & \E \big[A_n^{(t)}\big] \leq 2 n \cdot \E \left[I_t + N_t \right] 
    \quad \text{and} \quad \sum_{i=0}^{n-1}\E \big[A_i^{(t)}\big] 
    \leq n^2 \cdot\E \left[I_t + N_t \right], \tagthis \label{bound_A^{T}_RR}\\
    (\textbf{b}) \qquad 
    & \sum_{i=0}^{n-1} \E \Big[B_i^{(t)}\Big] \leq 2n^2 \cdot \E \left[I_t+ N_t\right] . \tagthis \label{bound_B^{T}_RR}
\end{align*}
\end{lem}


The results of Lemma~\ref{lem_momentum_variance_03_RR} below are direct consequences of Lemma \ref{lem_momentum_variance_01}, Lemma~\ref{lem_momentum_variance_02_RR}, and the fact that $L^2\xi_t^2 \leq \frac{3}{5}$.
\begin{lem} \label{lem_momentum_variance_03_RR}
Under the same setting as of Lemma~\ref{lem_momentum_variance_02_RR} and Assumption~\ref{as:A1}(b), it holds that
\begin{align*}
    (\textbf{a}) \ 
    &\E[I_t] \leq \left\{\begin{array}{ll} 
    L^2\xi_t^2 \Big[\frac{2}{3}\beta  \cdot \E[ I_{t-1} + N_{t-1} ]  + (1-\beta) \cdot \E[I_t + N_t ] \Big], &\text{if}~t > 1,  \vspace{1ex}\\
    L^2\xi_t^2  (1-\beta)^2 \cdot \E\left[I_t + N_t \right], &\text{if}~t = 1,
    \end{array}\right. \\
    (\textbf{b}) \ 
    &\E [J_t] \leq 
    \left\{\begin{array}{ll}
    2L^2\xi_t^2 \Big[ \beta \cdot \E [  I_{t-2} + N_{t-2} ]  + (1-\beta) \cdot \E[I_{t-1} + N_{t-1} ] \Big], &\text{if}~ t > 2,  \vspace{1ex}\\
    2L^2\xi_t^2 \cdot \E\left[I_{t-1} + N_{t-1} \right],  &\text{if}~ t =2,
    \end{array}\right.
    \\
    (\textbf{c}) \ &\text{If } L^2\xi_t^2 \leq \frac{3}{5}, \text{ then }  \E[I_t + N_t] \leq 
    \left\{\begin{array}{ll}
    \beta \cdot \E[I_{t-1} + N_{t-1}] + \frac{5 -3\beta}{2}\cdot \E[N_t], &\text{if}~ t > 1,   \vspace{1ex}\\ 
    \frac{5 -3\beta}{2} \cdot \E[N_1], &\text{if}~t = 1.  
    \end{array}\right.
\end{align*}
\end{lem}

We will present below Lemmas \ref{lem_momentum_variance_04_RR} and \ref{lem_momentum_variance_06_RR}, which play a similar role as previous Lemmas \ref{lem_momentum_variance_04} and \ref{lem_momentum_variance_06} in Theorem~\ref{thm_momentum_variance}.

\begin{lem} \label{lem_momentum_variance_04_RR}
Under the same conditions as in Lemma~\ref{lem_momentum_variance_03_RR}, for any $t\geq 2$, we have
\begin{align*}
    \E [\beta J_t + (1-\beta) I_t] & \leq  3(5 -3\beta) L^2\xi_t^2 \Big( (\Theta +n) \sum_{j=1}^{t}  \beta^{t-j} \E \left[ \big\Vert\nabla  F(w_{0}^{(j)}) \big\Vert ^2\right] +  \frac{\sigma^2}{1-\beta} \Big). 
    \tagthis \label{bound_beta_IJ_RR}
\end{align*}
\end{lem}



\begin{lem} \label{lem_momentum_variance_06_RR}
Under the same setting as in Theorem \ref{thm_momentum_variance_RR}, we have
\begin{align*} 
    \frac{\eta_1}{2}\E\left[ \norms{ \nabla F(  \tilde{w}_{0} )}^2 \right]  &\leq \frac{\E \left[F(  \tilde{w}_{0} )- F(\tilde{w}_{1} ) \right]}{(1-\beta)}
    + \frac{(1-\beta)\eta_1 }{4} \E \left[\big\Vert \nabla  F(\tilde{w}_{0})  \big\Vert ^2\right] + \frac{3 \sigma^2 (5 - 3\beta)L^2}{2n(1-\beta)} \cdot \xi_1^3. \tagthis \label{bound_F(w)_1_tilde_RR}
\end{align*}
\end{lem}

\subsection{The Proof of Theorem \ref{thm_momentum_variance_RR}: Key Estimate for Algorithm~\ref{sgd_momentum_shuffling_01}}

First, from the assumption $0 < \eta_t \leq \frac{1}{L\sqrt{D}}$, $t \geq 1$, we have $0 < \eta_t^2 \leq \frac{1}{D L^2}$. 
Next, from \eqref{define_xi}, we have $\xi_t = \max (\eta_t; \eta_{t-1})$ for $t > 1$ and $\xi_1 = \eta_1$, which lead to $0 < \xi_t^2 \leq \frac{1}{D L^2}$ for $t \geq 1$. 
Moreover, from the definition of $D = \max \left(\frac{5}{3}, \frac{6(5 -3\beta) (\Theta +n)}{n(1-\beta)}\right)$ in Theorem~\ref{thm_momentum_variance_RR}, we have $L^2\xi_t^2 \leq \frac{3}{5}$ and $6 L^2 \xi_t^2 (5 -3\beta) (\Theta+n) \leq n(1-\beta)$ for $t\geq 1.$


Similar to the estimate \eqref{estimate_alg_1} in Theorem \ref{thm_momentum_variance}, we get the following result:
\begin{align*} 
    F( w_0^{(t+1)} )& \leq F( w_0^{(t)} ) - \frac{\eta_t}{2} \norms{ \nabla F( w_0^{(t)} )}^2 + \frac{\eta_t}{2n} \Big[ \beta J_t + (1-\beta) I_t \Big].  
\end{align*}

Taking expectation of both sides and applying Lemma~\ref{lem_momentum_variance_04_RR}, we get:
\begin{align*} 
    \Exp{F( w_0^{(t+1)} )} & \leq \Exp{F( w_0^{(t)} ) } - \frac{\eta_t}{2} \Exp{  \norms{ \nabla F( w_0^{(t)} )}^2} + \frac{\eta_t}{2n} \E \Big[ \beta J_t + (1-\beta) I_t \Big] \\
    & \overset{\eqref{bound_beta_IJ_RR}}{\leq} \Exp{F( w_0^{(t)} ) } - \frac{\eta_t}{2} \Exp{  \norms{ \nabla F( w_0^{(t)} )}^2} \\
    &\quad + \frac{\eta_t}{2n} 3(5 -3\beta) L^2\xi_t^2 \bigg[ (\Theta +n) \sum_{j=1}^{t}  \beta^{t-j} \E \left[ \big\Vert\nabla  F(w_{0}^{(j)}) \big\Vert ^2\right] +  \frac{\sigma^2}{1-\beta} \bigg]\\
    & \leq \Exp{F( w_0^{(t)} ) } - \frac{\eta_t}{2} \Exp{  \norms{ \nabla F( w_0^{(t)} )}^2}
    \\
    &\quad+ \frac{ 3(5 -3\beta)\eta_t L^2\xi_t^2 (\Theta +n)}{2n}  \sum_{j=1}^{t}  \beta^{t-j} \E \left[ \big\Vert\nabla  F(w_{0}^{(j)}) \big\Vert ^2\right]  +  \frac{ 3 \sigma^2(5 -3\beta)\eta_t L^2\xi_t^2}{2n (1-\beta)}.
\end{align*}

Since $\xi_t$ satisfies $6 L^2 \xi_t^2 (5 -3\beta) (\Theta+n) \leq n(1-\beta)$ as proved above, we can deduce from the last estimate that
\begin{align*} 
    \Exp{F( w_0^{(t+1)} )}  & \leq \Exp{F( w_0^{(t)} ) } - \frac{\eta_t}{2} \Exp{  \norms{ \nabla F( w_0^{(t)} )}^2} + \frac{(1-\beta)\eta_t }{4} \sum_{j=1}^{t}  \beta^{t-j} \E \left[ \big\Vert\nabla  F(w_{0}^{(j)}) \big\Vert ^2\right] \\
    &\quad +  \frac{ 3 \sigma^2(5 -3\beta)\eta_t L^2\xi_t^2}{2n (1-\beta)}.
\end{align*}
Rearranging this inequality while noting that $\eta_t \leq \xi_t$ and $\tilde{w}_{t-1} = w_0^{(t)}$ we obtain that for $t \geq 2$:
\begin{equation*} 
\arraycolsep=0.2em
\begin{array}{lcl}
    \frac{\eta_t}{2} \Exp{  \norms{ \nabla F( \tilde{w}_{t-1} )}^2}  & \leq  & \Exp{F( \tilde{w}_{t-1}  ) - F( \tilde{w}_{t} )}  + \frac{(1-\beta)\eta_t }{4} \sum_{j=1}^{t}  \beta^{t-j} \E \left[ \big\Vert\nabla  F(\tilde{w}_{j-1}) \big\Vert ^2\right] +  \frac{ 3 \sigma^2(5 -3\beta)L^2}{2n (1-\beta)} \cdot \xi_t^3  \vspace{1ex}\\
    & \overset{1-\beta <1}{\leq} & \frac{\Exp{F( \tilde{w}_{t-1}  ) - F( \tilde{w}_{t} )} }{(1-\beta)} + \frac{(1-\beta)\eta_t }{4} \sum_{j=1}^{t}  \beta^{t-j} \E \left[ \big\Vert\nabla  F(\tilde{w}_{j-1}) \big\Vert ^2\right] +  \frac{ 3 \sigma^2(5 -3\beta)L^2}{2n (1-\beta)} \cdot \xi_t^3.
\end{array}
\end{equation*}

For $t = 1$, since $\xi_1 = \eta_1$ as previously defined in \eqref{define_xi}, from the result of Lemma \ref{lem_momentum_variance_06_RR} we have
\begin{align*} 
    \frac{\eta_1}{2}\E\left[ \norms{ \nabla F(  \tilde{w}_{0} )}^2 \right]  &\leq \frac{\E \left[F(  \tilde{w}_{0} )- F(\tilde{w}_{1} ) \right]}{(1-\beta)}
    + \frac{(1-\beta)\eta_1 }{4} \sum_{j=1}^{1}  \beta^{1-j} \E \left[\big\Vert \nabla  F(\tilde{w}_{j-1})  \big\Vert ^2\right] + \frac{3 \sigma^2 (5 - 3\beta)L^2}{2n(1-\beta)} \cdot \xi_1^3.
\end{align*}

Summing the previous estimate for $t := 2$ to $t := T$, and then using the last one, we obtain
\begin{align*} 
    \sum_{t=1}^{T} \frac{\eta_t}{2} \Exp{  \norms{ \nabla F( \tilde{w}_{t-1} )}^2}  
    & \leq \frac{F(  \tilde{w}_{0} )- F_*}{(1-\beta)} + \frac{(1-\beta) }{4}\sum_{t=1}^{T}\eta_t \sum_{j=1}^{t}  \beta^{t-j} \E \left[ \big\Vert\nabla  F(\tilde{w}_{j-1}) \big\Vert ^2\right] +  \frac{ 3 \sigma^2(5 -3\beta)L^2}{2n (1-\beta)}\sum_{t=1}^{T} \xi_t^3.
\end{align*}

Applying the result of Lemma \ref{lem_momentum_variance_05} to the last estimate, we get
\begin{align*} 
    \frac{1}{2} \sum_{t=1}^{T} \eta_t \Exp{  \norms{ \nabla F( \tilde{w}_{t-1} )}^2}  
    & \leq \frac{F(  \tilde{w}_{0} )- F_*}{(1-\beta)} + \frac{(1-\beta) }{4}\frac{1}{(1-\beta)} \sum_{t=1}^{T} \eta_t \E \left[\norms{ \nabla F( \tilde{w}_{t-1} )}^2 \right]  +  \frac{ 3 \sigma^2(5 -3\beta)L^2}{2n (1-\beta)}\sum_{t=1}^{T} \xi_t^3,
\end{align*}

which is equivalent to 
\begin{align*} 
     \sum_{t=1}^{T} \eta_t \Exp{  \norms{ \nabla F( \tilde{w}_{t-1} )}^2}  
    & \leq \frac{4\left[F(  \tilde{w}_{0} )- F_*\right]}{(1-\beta)}  +  \frac{6 \sigma^2(5 -3\beta)L^2}{n (1-\beta)}\sum_{t=1}^{T} \xi_t^3.
\end{align*}

Dividing both sides of the resulting estimate by $\sum_{t=1}^{T} \eta_t$, we obtain \eqref{eq_thm_momentum_variance_RR}.
Finally, due to the choice of $\hat{w}_T$ at Step~\ref{alg:A1_step6} in Algorithm~\ref{sgd_momentum_shuffling_01}, we have $\mathbb{E}\big[\| \nabla F( \hat{w}_T )  \|^2 \big]  = \frac{1}{\sum_{t=1}^T \eta_t} \sum_{t=1}^T \eta_t \mathbb{E}\big[ \norms{ \nabla F( \tilde{w}_{t-1} )}^2 \big]$.
\Eproof
\subsection{The Proof of Corollaries~\ref{co:constant_LR_RR} and \ref{co:diminishing_LR_RR}: Constant and Diminishing Learning Rates}\label{apdx:subsec:main_result_RR}

\begin{proof}[The proof of Corollary~\ref{co:constant_LR_RR}]
Since $T \geq 1$ and $\eta_t := \frac{\gamma n^{1/3}}{T^{1/3}}$, we have $\eta_t \geq \eta_{t+1}$. We also have $0 < \eta_t \leq \frac{1}{L\sqrt{D}}$ for all $t \geq 1$, and $\sum_{t=1}^T\eta_t = \gamma n^{1/3} T^{2/3} $ and $\sum_{t=1}^T\eta_{t-1}^3 = n \gamma^3$.
Substituting these expressions into \eqref{eq_thm_momentum_variance_RR}, we obtain
\begin{equation*} 
    \mathbb{E}\big[\| \nabla F( \hat{w}_T )  \|^2 \big]    \leq \frac{4\left[F(  \tilde{w}_{0} )- F_*\right]}{(1-\beta)\gamma n^{1/3} T^{2/3}}  +  \frac{6 \sigma^2(5 -3\beta)L^2}{(1-\beta)}\cdot \frac{\gamma^2 }{n^{1/3}T^{2/3}},
\end{equation*}
which is our desired result. 
\end{proof}

\begin{proof}[The proof of Corollary~\ref{co:diminishing_LR_RR}]
For $t \geq 1$, since $\eta_t = \frac{\gamma n^{1/3}}{(t+\lambda)^{1/3}}$, we have $\eta_t \geq \eta_{t+1}$. We also have $0 < \eta_t \leq \frac{1}{L\sqrt{D}}$ for all $t \geq 1$, and
\begin{equation*}
\begin{array}{lll}
& \sum_{t=1}^T\eta_t  & = \gamma n^{1/3} \sum_{t=1}^T\frac{1}{(t+\lambda)^{1/3}} \geq \gamma n^{1/3} \int_1^T\frac{d\tau}{(\tau + \lambda)^{1/3}} \geq \gamma n^{1/3}\left[ (T + \lambda)^{2/3} - (1 + \lambda)^{2/3} \right],  \vspace{1.5ex}\\
& \sum_{t=1}^T\eta_{t-1}^3 & = 2\eta_1^3 + \gamma^3 \sum_{t=3}^{T} \frac{1}{t - 1 + \lambda} 
\leq \frac{2n\gamma^3}{(1+\lambda)} + n\gamma^3 \int_{t=2}^T\frac{d\tau}{\tau - 1 + \lambda} \leq n\gamma^3 \left[ \frac{2}{(1+\lambda)} + \log (T - 1 +\lambda) \right].
\end{array}
\end{equation*}
Substituting these expressions into \eqref{eq_thm_momentum_variance_RR}, we obtain
\begin{equation*} 
    \mathbb{E}\big[\| \nabla F( \hat{w}_T )  \|^2 \big]    \leq \frac{4\left[F(  \tilde{w}_{0} )- F_*\right]}{(1-\beta)\gamma n^{1/3} \left[ (T + \lambda)^{2/3} - (1 + \lambda)^{2/3} \right]}  +  \frac{6 \sigma^2(5 -3\beta)L^2}{ (1-\beta)} \cdot  \frac{\gamma^2 \left[ \frac{2}{(1+\lambda)} + \log (T - 1 +\lambda) \right]}{n^{1/3} \left[ (T + \lambda)^{2/3} - (1 + \lambda)^{2/3} \right]},
\end{equation*}
Let $C_3$ and $C_4$ be defined respectively as
\begin{align*}
    C_3 :=  \frac{4\left[F(  \tilde{w}_{0} )- F_*\right]}{(1-\beta)\gamma }  +  \frac{12 \sigma^2(5 -3\beta)L^2 \gamma^2 }{(1-\beta) (1+\lambda)}, \quad \text{ and } \quad
    C_4 := \frac{6 \sigma^2(5 -3\beta)L^2 \gamma^2 }{ (1-\beta)}.
\end{align*}
Then, the last estimate leads to
\begin{equation*} 
    \mathbb{E}\big[\| \nabla F( \hat{w}_T )  \|^2 \big]    \leq \frac{C_3 + C_4 \log (T - 1 +\lambda)}{n^{1/3}\left[ (T + \lambda)^{2/3} - (1 + \lambda)^{2/3}\right]},
\end{equation*}
which completes our proof.
\end{proof}
\subsection{The Proof of Technical Lemmas}\label{appdx:subec:Tech_lemmas_proof_RR}

We now provide the proof of additional Lemmas that serve for the proof of Theorem~\ref{thm_momentum_variance_RR} above.
\subsubsection{Proof of Lemma~\ref{lem_momentum_variance_02_RR}: Upper Bounding  The Terms $ \E \big[A_n^{(t)}\big]$, $\sum_{i=0}^{n-1}\E \big[A_i^{(t)}\big]$, and $\sum_{i=0}^{n-1}\E \big[B^{(t)}_i\big]$}
In this proof, we will use  \citep{mishchenko2020random}[Lemma 1]  for sampling without replacement.
For the sake of references, we recall it here.

\begin{lem}[Lemma 1 in \cite{mishchenko2020random}]
Let $X_1, \cdots, X_n \in \R^d$ be fixed vectors, $\bar{X} := \frac{1}{n} \sum_{i=1}^n X_i$ be their average and $\sigma^2 := \frac{1}{n} \sum_{i=1}^n \|X_i -\bar{X}\|^2$ 
be the population variance. Fix any $k \in \{1,\cdots, n\}$, let $X_{\pi_1}, \cdots, X_{\pi_k}$ be sampled uniformly without replacement from $\{X_1, \cdots, X_n\}$ and $\bar{X}_\pi$ be their average. 
Then, the sample average and the variance are given, respectively by
\begin{align*}
    \E [\bar{X}_\pi] = \bar{X} \qquad \text{and} \qquad \E \left[ \|\bar{X}_\pi - \bar{X}\|^2 \right] = \frac{n-k}{k(n-1)} \sigma^2.
\end{align*}
\end{lem}
Using this result, we now prove Lemma~\ref{lem_momentum_variance_02_RR} as follows.

(\textbf{a}) Let first upper bound the term $A_i^{(t)}$ defined by \eqref{define_A}.
Indeed, for $i = 0, \cdots, n-1$ and $t \geq 1$, we have
\begin{align*}
    A_i^{(t)} = \Big\Vert  \sum_{j=0}^{i-1}  g_{j}^{(t)} \Big\Vert ^2 
    &\overset{(a)}{\leq} 2\Big\Vert   \sum_{j=0}^{i-1}  \Big(g_{j}^{(t)} - \nabla f ( w_{0}^{(t)} ; \pi^{(t)} ( j + 1 ) ) \Big) \Big\Vert ^2  
    + 2\Big\Vert  \sum_{j=0}^{i-1}  \nabla f ( w_{0}^{(t)} ; \pi^{(t)} ( j + 1 ) )  \Big\Vert ^2  \\
    &\overset{(a)}{\leq} 
    2 i  \sum_{j=0}^{i-1} \Big\Vert g_{j}^{(t)} - \nabla f ( w_{0}^{(t)} ; \pi^{(t)} ( j + 1 ) )  \Big\Vert ^2  + 2\Big\Vert  \sum_{j=0}^{i-1}  \nabla f ( w_{0}^{(t)} ; \pi^{(t)} ( j + 1 ) )  \Big\Vert ^2\\
    &\leq 
    2 i  \sum_{j=0}^{n-1} \Big\Vert g_{j}^{(t)} - \nabla f ( w_{0}^{(t)} ; \pi^{(t)} ( j + 1 ) )  \Big\Vert ^2  + 2\Big\Vert  \sum_{j=0}^{i-1}  \nabla f ( w_{0}^{(t)} ; \pi^{(t)} ( j + 1 ) )  \Big\Vert ^2 
    \\
    & \overset{\eqref{define_I}}{=} 2 i  I_t  + 2\Big\Vert  \sum_{j=0}^{i-1}  \nabla f ( w_{0}^{(t)} ; \pi^{(t)} ( j + 1 ) )  \Big\Vert ^2,
\end{align*}
where we have used the Cauchy-Schwarz inequality in (\textit{a}).

Now taking expectation conditioned on $\mathcal{F}_t$, we get
\begin{align*}
    \E_t \big[A_i^{(t)}\big] &\leq  2 i \cdot \E_t \left[I_t\right]
    + 2 \E_t \left[ \Big\Vert  \sum_{j=0}^{i-1}  \nabla f ( w_{0}^{(t)} ; \pi^{(t)} ( j + 1 ) )  \Big\Vert ^2 \right].
\end{align*}

Note that $\E_t \left[ \nabla f ( w_{0}^{(t)} ; \pi^{(t)} ( j + 1 ) )  \right] = \nabla F ( w_{0}^{(t)})$ for every index $j = 0, 1, \cdots, i-1.$ Hence, we have
\begin{align*}
    \E_t \big[A_i^{(t)}\big] &\leq  2 i \cdot \E_t \left[I_t\right]
    + 2 \E_t \bigg[ \Big\Vert  \sum_{j=0}^{i-1}  \nabla f ( w_{0}^{(t)} ; \pi^{(t)} ( j + 1 ) )  \Big\Vert ^2 \bigg]\\
    &=  2 i \cdot \E_t \left[I_t\right]
    + 2 \big\Vert i \nabla F ( w_{0}^{(t)})\big\Vert^2
    + 2 \E_t \bigg[ \Big\Vert  \sum_{j=0}^{i-1} \nabla f ( w_{0}^{(t)} ; \pi^{(t)} ( j + 1 ) ) - i \nabla F ( w_{0}^{(t)})  \Big\Vert ^2 \bigg]\\
    &=  2 i \cdot \E_t \left[I_t\right]
    + 2 i^2 \big\Vert\nabla F ( w_{0}^{(t)})\big\Vert^2
    + 2 i^2 \E_t \bigg[ \Big\Vert  \frac{1}{i} \sum_{j=0}^{i-1}  \nabla f ( w_{0}^{(t)} ; \pi^{(t)} ( j + 1 ) ) - \nabla F ( w_{0}^{(t)})  \Big\Vert ^2 \bigg]\\
    &=  2 i \cdot \E_t \left[I_t\right]
    + 2 i^2 \big\Vert\nabla F ( w_{0}^{(t)})\big\Vert^2
    + 2 i^2 \frac{n-i}{i(n-1)} \frac{1}{n} \sum_{j=0}^{n-1} \big\Vert \nabla f ( w_{0}^{(t)} ; \pi^{(t)} ( j + 1 ) ) - \nabla F ( w_{0}^{(t)})  \big\Vert ^2 \\
    & \overset{(b)}{\leq}  2 i \cdot \E_t \left[I_t\right]
    + 2 i^2 \big\Vert\nabla F ( w_{0}^{(t)})\big\Vert^2
    + \frac{2i(n-i)}{(n-1)} \Big[\Theta \big\Vert \nabla  F(w_{0}^{(t)})  \big\Vert ^2 + \sigma^2 \Big],
    %
\end{align*}
where we apply the sample variance Lemma from \citep{mishchenko2020random}
and (\textit{b}) follows from Assumption~\ref{as:A1}(c).

Letting $i = n$ in the above estimate and taking total expectation, we obtain the first estimate of Lemma \ref{lem_momentum_variance_02_RR}, i.e.: 
\begin{align*}
    \E \big[A_n^{(t)}\big] &\leq 2  n\cdot \E \left[I_t\right]
    + 2 n^2 \E \left[\big\Vert\nabla F ( w_{0}^{(t)})\big\Vert^2 \right] \leq 2 n \cdot \E\left[I_t + N_t\right].
\end{align*}
Now we are ready to calculate the sum $\sum_{i=0}^{n-1} \E_t \big[A_i^{(t)}\big]$ as
\begin{align*}
    \sum_{i=0}^{n-1}\E_t \big[A_i^{(t)}\big] 
    &\leq 2 \sum_{i=0}^{n-1}i \cdot \E_t \left[I_t\right]
    + 2  \sum_{i=0}^{n-1} i^2 \cdot \big\Vert\nabla F ( w_{0}^{(t)})\big\Vert^2
    + 2 \sum_{i=0}^{n-1} \frac{i (n-i)}{(n-1)} \Big[\Theta \big\Vert \nabla  F(w_{0}^{(t)})  \big\Vert ^2 + \sigma^2 \Big]\\
    &\leq n^2 \cdot \E_t \left[I_t\right]
    + \frac{2n^3}{3} \big\Vert\nabla F ( w_{0}^{(t)})\big\Vert^2
    + \frac{n(n+1)}{3}\Big[\Theta \big\Vert \nabla  F(w_{0}^{(t)})  \big\Vert ^2 + \sigma^2 \Big]\\
    &\leq n^2 \cdot \E_t \left[I_t\right]
    + n^3 \big\Vert\nabla F ( w_{0}^{(t)})\big\Vert^2
    + n^2\Big[\Theta \big\Vert \nabla  F(w_{0}^{(t)})  \big\Vert ^2 + \sigma^2 \Big] \\
    &\leq n^2 \cdot \E_t \left[I_t\right]
   + n^2\Big[ (n +\Theta) \big\Vert \nabla  F(w_{0}^{(t)})  \big\Vert ^2 + \sigma^2 \Big] 
   \overset{\eqref{define_N}}{\leq} n^2 \left(\E_t \left[I_t\right] + N_t \right) , \  t \geq 1,
\end{align*}

where we use the facts that $\sum_{i=0}^{n-1} i =\frac{n(n-1)}{2}\leq \frac{n^2}{2}$ and $\sum_{i=0}^{n-1} i^2 = \frac{n(n-1)(2n-1)}{6} \leq \frac{n^3}{3}$. 

Taking total expectation, we have the second estimate of Lemma \ref{lem_momentum_variance_02_RR}.

(\textbf{b})
Using similar argument as above, for $i =0, 1, \cdots, n-1$ and $t \geq 1$, we can derive
\begin{align*}
    B_i^{(t)} &= \Big\Vert \sum_{j = i}^{n-1} g_{j}^{(t)}\Big\Vert^2 
    \overset{(a)}{\leq}
    2 \Big\Vert   \sum_{j = i}^{n-1}  \Big(g_{j}^{(t)} - \nabla f ( w_{0}^{(t)} ; \pi^{(t)} ( j + 1 ) ) \Big) \Big\Vert ^2 + 2\Big\Vert  \sum_{j=i}^{n-1}  \nabla f ( w_{0}^{(t)} ; \pi^{(t)} ( j + 1 ) )  \Big\Vert ^2 \\
    &\overset{(a)}{\leq}
    2 (n-i)  \sum_{j=i}^{n-1} \Big\Vert g_{j}^{(t)} - \nabla f ( w_{0}^{(t)} ; \pi^{(t)} ( j + 1 ) )  \Big\Vert ^2 + 2\Big\Vert  \sum_{j=i}^{n-1}  \nabla f ( w_{0}^{(t)} ; \pi^{(t)} ( j + 1 ) )  \Big\Vert ^2 \\
    & \overset{\eqref{define_I}}{\leq} 
    2 (n-i)  I_t + 2\Big\Vert  \sum_{j=i}^{n-1}  \nabla f ( w_{0}^{(t)} ; \pi^{(t)} ( j + 1 ) )  \Big\Vert ^2,
\end{align*}
where we use again the Cauchy-Schwarz inequality in (\textit{a}).

Taking expectation conditioned on $\mathcal{F}_t$, we get 
\begin{align*}
    \E_t \big[B_i^{(t)}\big] &\leq  2 (n-i) \cdot \E_t \left[I_t\right]
    + 2 \E_t \bigg[ \Big\Vert  \sum_{j=i}^{n-1}  \nabla f ( w_{0}^{(t)} ; \pi^{(t)} ( j + 1 ) )  \Big\Vert ^2 \bigg].
\end{align*}

Note that $\E_t \big[ \nabla f ( w_{0}^{(t)} ; \pi^{(t)} ( j + 1 ) )  \big] = \nabla F ( w_{0}^{(t)})$ for every index $j = i, \cdots, n-1.$ Hence
\begin{align*}
    \E_t \big[B_i^{(t)}\big] &\leq  2 (n-i) \cdot \E_t \left[I_t\right]
    + 2 \E_t \bigg[ \Big\Vert  \sum_{j=i}^{n-1}  \nabla f ( w_{0}^{(t)} ; \pi^{(t)} ( j + 1 ) )  \Big\Vert ^2 \bigg]\\
    &=  2 (n-i)  \cdot \E_t \left[I_t\right]
    + 2 \big\Vert (n-i) \nabla F ( w_{0}^{(t)})\big\Vert^2
    + 2 \E_t \bigg[ \Big\Vert  \sum_{j=i}^{n-1} \nabla f ( w_{0}^{(t)} ; \pi^{(t)} ( j + 1 ) ) - (n-i)  \nabla F ( w_{0}^{(t)})  \Big\Vert ^2 \bigg]\\
    &=  2 (n-i)  \cdot \E_t \left[I_t\right]
    + 2 (n-i) ^2 \Big\Vert\nabla F ( w_{0}^{(t)})\Big\Vert^2\\
    & \quad
    + 2 (n-i) ^2 \E_t \bigg[ \Big\Vert  \frac{1}{n-i} \sum_{j=i}^{n-1}  \nabla f ( w_{0}^{(t)} ; \pi^{(t)} ( j + 1 ) ) - \nabla F ( w_{0}^{(t)})  \Big\Vert ^2 \bigg]\\
    & \overset{(b)}{\leq}  2 (n-i)  \cdot \E_t \left[I_t\right]
    + 2 (n-i)^2 \Big\Vert\nabla F ( w_{0}^{(t)})\Big\Vert^2
    + \frac{2i(n-i)}{(n-1)} \Big[\Theta \big\Vert \nabla  F(w_{0}^{(t)})  \big\Vert ^2 + \sigma^2 \Big],
\end{align*}
where we apply again the sample variance Lemma from \citep{mishchenko2020random}
and (\textit{b}) follows from Assumption~\ref{as:A1}(c).
Finally, for all $t\geq 1$, we calculate the sum $\sum_{i=0}^{n-1} \E_t \Big[B_i^{(t)}\Big]$ similarly as follows
\begin{align*}
    \sum_{i=0}^{n-1} \E_t \Big[B_i^{(t)}\Big]  &\leq 2 \sum_{i=0}^{n-1}(n-i)  \cdot \E_t \left[I_t\right]
    + 2 \sum_{i=0}^{n-1} (n-i)^2 \Big\Vert\nabla F ( w_{0}^{(t)})\Big\Vert^2
    + \sum_{i=0}^{n-1} \frac{2i(n-i)}{(n-1)} \Big[\Theta \big\Vert \nabla  F(w_{0}^{(t)})  \big\Vert ^2 + \sigma^2 \Big]\\
    &\leq 2 n^2  \cdot \E_t \left[I_t\right]
    + \frac{n(n+1)(2n+1)}{3} \Big\Vert\nabla F ( w_{0}^{(t)})\Big\Vert^2
    + \frac{n(n+1)}{3} \Big[\Theta \big\Vert \nabla  F(w_{0}^{(t)})  \big\Vert ^2 + \sigma^2 \Big]\\
    &\leq 2n^2 \cdot \E_t \left[I_t\right]
    + 2n^3 \Big\Vert\nabla F ( w_{0}^{(t)})\Big\Vert^2
    + 2n^2\Big[\Theta \big\Vert \nabla  F(w_{0}^{(t)})  \big\Vert ^2 + \sigma^2 \Big]\\
    &\leq 2n^2 \cdot \E_t \left[I_t\right]
   + 2n^2\Big[ (n +\Theta) \big\Vert \nabla  F(w_{0}^{(t)})  \big\Vert ^2 + \sigma^2 \Big] \overset{\eqref{define_N}}{\leq} 2n^2 \left(\E_t \left[I_t\right] + N_t \right).
\end{align*}

Taking total expectation, we get \eqref{bound_B^{T}_RR}.
\Eproof

\subsubsection{Proof of Lemma~\ref{lem_momentum_variance_03_RR}: Upper Bounding The Terms $\E [I_t]$, $\E[J_t]$, and $\E[I_t + N_t]$}
(\textbf{a})
First, for every $t \geq 1$ and $g_{i}^{(t)} := \nabla f ( w_{i}^{(t)} ; \sigma^{(t)} ( i + 1 ) )$, we have 
\begin{align*}
    I_t &= \sum_{i=0}^{n-1} \big \Vert \nabla f ( w_{0}^{(t)} ; \pi^{(t)} ( i + 1 ) )  - g_{i}^{(t)} \big\Vert^2 \\
    &\overset{\tiny \textrm{Step~\ref{alg:A1_step3}}}{=} \sum_{i=0}^{n-1} \big \Vert \nabla f ( w_{0}^{(t)} ; \pi^{(t)} ( i + 1 ) )  - \nabla f ( w_{i}^{(t)} ; \pi^{(t)} ( i + 1 ) )  \big\Vert^2 \\
    &\leq L^2 \sum_{i=0}^{n-1} \big\Vert w_{i}^{(t)} - w_{0}^{(t)} \big\Vert^2, 
    \tagthis \label{bound_I_t_02_RR}
\end{align*}
where the last inequality follows from Assumption~\ref{as:A1}(b).

Applying Lemma \ref{lem_momentum_variance_01} to \eqref{bound_I_t_02_RR}, for $t> 1$, we have
\begin{align*}
    I_t &  \leq \frac{L^2\xi_t^2}{n^2} \sum_{i=0}^{n-1}\bigg[ \beta \frac{i^2}{n^2} A_n^{(t-1)} + (1-\beta) A_i^{(t)} \bigg] \\
    &\leq  \frac{L^2\xi_t^2}{n^2}  \bigg[ \beta \frac{ A_n^{(t-1)} }{n^2}  \sum_{i=0}^{n-1} i^2 +  (1-\beta) \sum_{i=0}^{n-1} A_i^{(t)} \bigg]\\
    &\overset{(a)}{\leq} \frac{L^2\xi_t^2}{n^2}  \bigg[ \beta \cdot \frac{n}{3} \cdot A_n^{(t-1)}   +  (1-\beta) \sum_{i=0}^{n-1} A_i^{(t)} \bigg]
\end{align*}
where (\textit{a}) follows from the fact that $\sum_{i=0}^{n-1} i^2 \leq \frac{n^3}{3}$.
Taking expectation and applying result \eqref{bound_A^{T}_RR} in Lemma \ref{lem_momentum_variance_02_RR} we get 
\begin{align*}
    \E[I_t] &\leq  \frac{L^2\xi_t^2}{n^2}  \bigg[ \beta \cdot \frac{n}{3} \cdot \E\Big[A_n^{(t-1)}\Big]   +  (1-\beta) \sum_{i=0}^{n-1} \E \Big[A_i^{(t)}\Big] \bigg]\\
    &\leq \frac{L^2\xi_t^2}{n^2}  \bigg[ \beta \cdot \frac{n}{3} \cdot 2n\cdot \E \big[  I_{t-1} + N_{t-1}\big]  +  (1-\beta) n^2 \E \left[I_t + N_t  \right] \bigg]\\
    & \leq L^2\xi_t^2 \bigg[ \frac{2}{3}\beta \cdot \E \big[ I_{t-1} + N_{t-1}\big] + (1-\beta) \E \left[I_t + N_t \right] \bigg].
\end{align*}

For $t = 1$, by applying Lemma \ref{lem_momentum_variance_01} and \ref{lem_momentum_variance_02_RR} consecutively to \eqref{bound_I_t_02_RR} we have
\begin{align*}
    I_t & \leq L^2 
    \sum_{i=0}^{n-1} \big\Vert w_{i}^{(t)} - w_{0}^{(t)} \big\Vert^2 \overset{\eqref{update_distance_t}}{\leq} \frac{L^2\xi_t^2}{n^2} \sum_{i=0}^{n-1}\left( (1-\beta)^2 A_i^{(t)} \right) 
    \leq  \frac{L^2\xi_t^2 (1-\beta)^2}{n^2} \sum_{i=0}^{n-1} A_i^{(t)}.
\end{align*}
Similarly we get
\begin{align*}
    \E [I_t] &\overset{\eqref{bound_A^{T}_RR}}{\leq}  \frac{L^2\xi_t^2 (1-\beta)^2}{n^2} n^2 \cdot \E \left[I_t+ N_t \right]  
    = L^2\xi_t^2 (1-\beta)^2 \E \left[I_t+ N_t \right]. \tagthis \label{bound_I=1_RR}
\end{align*}

(\textbf{b})
Using the similar argument as above we can derive that 
\begin{align*}
    J_t &= \sum_{i=0}^{n-1} \Big \Vert \nabla f ( w_{0}^{(t)} ; \pi^{(t-1)} ( i + 1 ) )  - g_{i}^{(t-1)} \Big\Vert^2 \\
    &\overset{\tiny\textrm{Step~\ref{alg:A1_step3}}}{=} \sum_{i=0}^{n-1} \Big \Vert \nabla f ( w_{0}^{(t)} ; \pi^{(t-1)} ( i + 1 ) )  - \nabla f ( w_{i}^{(t-1)} ; \pi^{(t-1)} ( i + 1 ) )  \Big\Vert^2 \\
    & \overset{\eqref{eq:Lsmooth_basic}}{\leq} L^2 
    \sum_{i=0}^{n-1} \big\Vert w_{i}^{(t-1)} - w_{0}^{(t)} \big\Vert^2.
\end{align*}
Applying Lemma \ref{lem_momentum_variance_01} for $t > 2$ to the above estimate, we get
\allowdisplaybreaks
\begin{align*}
    J_t
    &\leq \frac{L^2\xi_t^2}{n^2} \sum_{i=0}^{n-1}  \left[ \beta   \frac{(n-i)^2}{n^2} A_n^{(t-2)} + (1-\beta) B_i^{(t-1)} \right] \\
    &\leq \frac{L^2\xi_t^2}{n^2} \left[ \beta \frac{ A_n^{(t-2)} }{n^2}  \sum_{i=0}^{n-1} (n-i)^2 +  (1-\beta) \sum_{i=0}^{n-1} B_i^{(t-1)} \right] \\
    &\leq \frac{L^2\xi_t^2}{n^2} \left[ \beta n \cdot A_n^{(t-2)}  +  (1-\beta) \sum_{i=0}^{n-1} B_i^{(t-1)} \right],
\end{align*}
where the last line follows from the fact that $\sum_{i=0}^{n-1} (n-i)^2 \leq n^3$.
Similar to the previous part, taking total expectation we get
\begin{align*}
    \E \left[J_t\right]&\leq \frac{L^2\xi_t^2}{n^2} \left[ \beta n \cdot \E \Big[ A_n^{(t-2)} \Big]  +  (1-\beta) \sum_{i=0}^{n-1} \E \Big[B_i^{(t-1)} \Big] \right]\\
    &\overset{\eqref{bound_B^{T}_RR}}{\leq} \frac{L^2\xi_t^2}{n^2} \left[ \beta n \cdot \E \Big[ A_n^{(t-2)} \Big]  +  (1-\beta) 2n^2 \cdot \E \left[I_{t-1}+ N_{t-1} \right]  \right]\\
    &\overset{\eqref{bound_A^{T}_RR}}{\leq} \frac{L^2\xi_t^2}{n^2} \Big[ 2\beta n^2 \cdot \E \left[I_{t-2}+ N_{t-2}\right]    +  2(1-\beta) n^2 \cdot \E \left[I_{t-1}+ N_{t-1}\right] \Big]\\
     &\leq 
    2L^2\xi_t^2\Big[ \beta \cdot\E \left[I_{t-2} + N_{t-2}\right] + (1-\beta) \cdot \E \left[I_{t-1}+ N_{t-1}\right] \Big].
\end{align*}


Now applying Lemma \ref{lem_momentum_variance_01} for the special case $t = 2$, we obtain
\begin{align*}
    J_t &\leq \frac{L^2\xi_t^2}{n^2} \sum_{i=0}^{n-1} (1-\beta)^2 B_i^{(t-1)}.
\end{align*}
Therefore, we have
\begin{align*}
    \E[J_t] &\leq \frac{L^2\xi_t^2}{n^2} (1-\beta)^2 \sum_{i=0}^{n-1} \E \left[ B_i^{(t-1)} \right]
    \leq \frac{L^2\xi_t^2}{n^2}  (1-\beta)^2 \cdot 2n^2 \E\left [I_{t-1} + N_{t-1}\right]
    \leq  2L^2\xi_t^2 \cdot \E \left[I_{t-1} + N_{t-1} \right],
\end{align*}
where we have used the result \eqref{bound_B^{T}_RR} from Lemma \ref{lem_momentum_variance_02_RR} and the fact that $1 - \beta \leq 1$. 

(\textbf{c}) 
First, from the result of Part \textbf{(a)}, for $t > 1$, we have
\begin{align*}
    &\E [I_t] \leq L^2\xi_t^2 \left[  \frac{2}{3}\beta \E [I_{t-1}+ N_{t-1}]  +   (1-\beta) \E [N_{t}] \right] +  (1-\beta) L^2\xi_t^2 \E [I_t].
\end{align*}
Since $L^2\xi_t^2 \leq \frac{3}{5}$ (due to the choice of our learning rate) and $1-\beta \leq 1$, we further get
\begin{align*}
    \E [I_t] &\leq   \frac{2}{5} \left[ \beta  \E [I_{t-1}+ N_{t-1}]  + \frac{3}{2} (1-\beta) \E [N_{t}] \right] + \frac{3}{5} \E [I_t],
\end{align*}
which is equivalent to 
\begin{align*}
    \E [I_t] \leq \beta  \E [I_{t-1} + N_{t-1}]+   \frac{3}{2} (1-\beta) \E [N_{t}].
\end{align*}
Adding $\E [N_{t}]$ to both sides of this inequality, for $t > 1$, we get
\begin{align*}
    \E [I_t + N_{t}] \leq \beta \E [I_{t-1}+ N_{t-1}]  +   \frac{5 -3\beta}{2} \E [N_{t}]. \tagthis \label{bound_sum_I_N>1_RR}
\end{align*}
Finally, for $t = 1$, since $L^2\xi_t^2 \leq \frac{3}{5}$ and $1-\beta \leq 1$, we have
\begin{align*}
    \E [I_t] & \leq  L^2\xi_t^2  (1-\beta)^2 \Big(\E [I_t] + \E [N_{t}] \Big) 
    \leq \frac{3}{5} (1-\beta)^2 \E [I_t] + \frac{3}{5} (1-\beta)^2 \E [N_{t}] \leq \frac{3}{5} \E [I_t] + \frac{3}{5} (1-\beta) \E [N_{t}],
\end{align*}
which is equivalent to 
\begin{align*}
    \frac{2}{5} \E [I_t] \leq \frac{3}{5} (1-\beta) \E [N_{t}].
    \tagthis \label{bound_I_1_RR}
\end{align*}
This leads to our desired result for $t=1$ as
\begin{align*}
    \E[I_1 +N_1] \leq \frac{3}{2} (1-\beta) \E[N_1] + \E[N_1] = \frac{5 - 3\beta}{2} \E[N_1]. \tagthis \label{bound_sum_I_N=1_RR}
\end{align*}
Combining two cases, we obtain the desired result in Part (\textbf{c}).
\Eproof

\subsubsection{Proof of Lemma~\ref{lem_momentum_variance_04_RR}: Upper Bounding The Key Quantity $\E [\beta J_t + (1-\beta) I_t]$}
First, we analyze the case $t>2$ using the results of Lemma \ref{lem_momentum_variance_03_RR} as follows:
\begin{align*}
    \beta \E[J_t] + (1-\beta) \E[I_t] & \leq 
    2\beta L^2\xi_t^2 \left[ \beta  \E[I_{t-2} + N_{t-2}] + (1-\beta)  \E [I_{t-1} + N_{t-1}]\right] \\ 
    & \quad + (1-\beta) L^2\xi_t^2 \left[ \frac{2}{3}\beta \E[I_{t-1} + N_{t-1}] +    (1-\beta) \E[I_{t} + N_{t}] \right] \\
    & \leq \frac{2}{3} L^2\xi_t^2 \left[  2 \E[I_{t} + N_{t}]  + 4\beta \E [I_{t-1} + N_{t-1}]
    + 3\beta^2  \E[I_{t-2} + N_{t-2}] \right] \\
    & =: \frac{2}{3} L^2\xi_t^2 S_t, \tagthis \label{define_S_RR}
\end{align*}
where the last line follows since $1 - \beta \leq 1$ and $S_t := 2\E[I_t + N_t] + 4\beta \E[I_{t-1} + N_{t-1}] + 3\beta^2 \E[I_{t-2} + N_{t-2}]$ for $t>2$.

Next, we bound the term $S_t$ in \eqref{define_S_RR} as
\begin{align*}
    S_t &= 2 \E[I_{t} + N_{t}] 
    + 4\beta \E[I_{t-1} + N_{t-1}]
    + 3\beta^2 \E [I_{t-2} + N_{t-2}]\\
    &\leq 2  \left[ \beta \E[I_{t-1} + N_{t-1}] + \frac{5 -3\beta}{2} \E[N_t] \right]
    + 4\beta \E[I_{t-1} + N_{t-1}]
    + 3\beta^2  \E[I_{t-2} + N_{t-2}] &\text{apply \eqref{bound_sum_I_N>1_RR} for} \ t\\
    &= (5 -3\beta) \E[N_t]  
    + 6\beta \E[I_{t-1} + N_{t-1}]
    + 3\beta^2  \E[I_{t-2} + N_{t-2}] \\
    &\leq (5 -3\beta) \E[N_t] + 6\beta \Big[ \beta \E[I_{t-2} + N_{t-2}] + \frac{5 -3\beta}{2} \E[N_{t-1}] \Big]
    + 3\beta^2 \E [I_{t-2} + N_{t-2}] &\text{apply \eqref{bound_sum_I_N>1_RR} for} \ t-1\\
    &= (5 -3\beta) \E[N_t] + 3 \beta (5 -3\beta)\E[N_{t-1}]
    + 9\beta^2  \E[I_{t-2} + N_{t-2}].
\end{align*}
We consider the last term, which can be bounded as
\begin{align*}
    \beta^2  \E[I_{t-2} + N_{t-2}]
    &= \beta^3 \E[I_{t-3} + N_{t-3}] + \frac{5 -3\beta}{2} \beta^2 \E[N_{t-2}] &\text{apply \eqref{bound_sum_I_N>1_RR} recursively for} \ t-2 > 1\\
    &= \beta^{t-1} \E[I_{1} + N_{1}] +\frac{5 -3\beta}{2}   \sum_{j=2}^{t-2}  \beta^{t-j} \E[N_j]  \\
    &\leq \beta^{t-1} \frac{5 -3\beta}{2} \E[N_1] +\frac{5 -3\beta}{2}   \sum_{j=2}^{t-2}  \beta^{t-j}\E [N_j] 
    &\text{apply \eqref{bound_sum_I_N=1_RR}} \\
    &\leq \frac{5 -3\beta}{2}   \sum_{j=1}^{t-2}  \beta^{t-j}\E[ N_j].
\end{align*}

Note that this bound is also true for the case $t-2 = 1$:
\begin{align*}
    \beta^2  \E[I_{t-2} + N_{t-2}] &= \beta^2  \E[I_{1} + N_{1}] \leq \beta^2  \frac{5 -3\beta}{2} \E[ N_j]  &\text{apply \eqref{bound_sum_I_N=1_RR}} 
\end{align*}

Substituting this inequality into $S_t$, for $t > 2$, we get
\begin{align*}
    S_t &\leq (5 -3\beta) \E[N_t] + 3 \beta (5 -3\beta)\E[N_{t-1}]
    + 9\beta^2  \E[I_{t-2} + N_{t-2}] \\
    &\leq \frac{9}{2}(5 -3\beta) \E[ N_t] + \frac{9}{2} (5 -3\beta)  \beta \E[N_{t-1}]
    +  \frac{9}{2} (5 -3\beta)  \sum_{j=1}^{t-2}  \beta^{t-j} \E[N_j] \\
    & \leq \frac{9}{2} (5 -3\beta)   \sum_{j=1}^{t}  \beta^{t-j}\E[ N_j].
\end{align*}
Now we analyze similarly for the case $t=2$ as follows:
\begin{align*}
    \beta \E[J_2] + (1-\beta) \E[I_2] & \leq 
    2\beta L^2\xi_2^2  \E [I_{1} + N_{1}]
    + (1-\beta) L^2\xi_2^2 \left[ \frac{2}{3}\beta \E [I_{1} + N_{1}] +    (1-\beta)\E [I_{2} + N_{2}] \right]\\
    & \leq \frac{2}{3} L^2\xi_2^2 \big[  2 \E[I_{2} + N_{2}]  + 4\beta \E[I_{1} + N_{1}] \big] \\
    & =: \frac{2}{3}L^2\xi_2^2 S_2, 
\end{align*}
where the last line follows since $1 - \beta \leq 1$ and $S_2 := 2 \E[I_{2} + N_{2}] + 4\beta \E [I_{1} + N_{1}]$. 

Next, we bound the term $S_2$ as follows:
\begin{align*}
    S_2 &= 2 \E[I_{2} + N_{2}]  + 4\beta \E[I_{1} + N_{1}] \\
    & \leq 2  \left(\beta \E[I_{1} + N_{1}] + \frac{5 -3\beta}{2} \E[N_2] \right)  + 4\beta \E[I_{1} + N_{1}] &\text{apply \eqref{bound_sum_I_N>1_RR} for} \ t=2\\
    &= (5 -3\beta) \E[N_2]   + 6\beta \E[I_{1} + N_{1}] \\
    & \leq (5 -3\beta) \E[N_2] + 6\beta  \frac{5 -3\beta}{2} \E[N_1]  &\text{apply \eqref{bound_sum_I_N=1_RR}}\\
    & = (5 -3\beta)\E[ N_2] + 3 \beta (5 -3\beta)\E[N_1] \\
    & \leq \frac{9}{2} (5 -3\beta)   \sum_{j=1}^{2}  \beta^{2-j} \E[N_j].
\end{align*}
Hence, the statements $\E[\beta J_t + (1-\beta)I_t] \leq \frac{2}{3}L^2\xi_t^2 S_t$ and $S_t \leq \frac{9}{2} (5 -3\beta)   \sum_{j=1}^{t}  \beta^{t-j} \E[N_j] $ are true for every $t \geq 2$. 

Combining these two cases, we have the following estimate
\begin{align*}
    \E[\beta J_t + (1-\beta)I_t] \leq \frac{2}{3}L^2\xi_t^2 S_t
    &\leq L^2\xi_t^2 \cdot \frac{2}{3} \cdot \frac{9}{2} (5 -3\beta)   \sum_{j=1}^{t}  \beta^{t-j} \E[N_j]\\
    &\overset{\eqref{define_N}}{\leq} 3 L^2\xi_t^2 (5 -3\beta)   \sum_{j=1}^{t}  \beta^{t-j} \E \Big[(\Theta+n) \big\Vert \nabla  F(w_{0}^{(j)}) \big\Vert ^2  + \sigma^2 \Big] \\
    &\leq 3 L^2\xi_t^2  (5 -3\beta)\bigg[ (\Theta+n) \sum_{j=1}^{t}  \beta^{t-j} \E \left[ \big\Vert \nabla F(w_{0}^{(j)}) \big\Vert ^2 \right] +  \sum_{j=1}^{t}  \beta^{t-j} \sigma^2 \bigg] \\
    &\leq 3 L^2\xi_t^2 (5 -3\beta)  \bigg[ (\Theta+n) \sum_{j=1}^{t}  \beta^{t-j} \E \left[\big\Vert \nabla F(w_{0}^{(j)}) \big\Vert ^2 \right] +  \frac{ \sigma^2}{1-\beta} \bigg],  
\end{align*}
where the last line follows since $\sum_{j=1}^{t}  \beta^{t-j} \leq \frac{1}{1-\beta}$ for every $t\geq2$.
Hence, we have proved \eqref{bound_beta_IJ_RR}.
\Eproof
\subsubsection{The Proof of Lemma~\ref{lem_momentum_variance_06_RR}}
First, from the assumption $0 < \eta_t \leq \frac{1}{L\sqrt{D}}$, $t \geq 1$, we have $0 < \eta_t^2 \leq \frac{1}{D L^2}$. 
Next, from \eqref{define_xi}, we have $\xi_t = \max (\eta_t; \eta_{t-1})$ for $t > 1$ and $\xi_1 = \eta_1$, which lead to $0 < \xi_t^2 \leq \frac{1}{D L^2}$ for $t \geq 1$. 
Moreover, from the definition of $D = \max \left(\frac{5}{3}, \frac{6(5 -3\beta) (\Theta +n)}{n(1-\beta)}\right)$ in Theorem~\ref{thm_momentum_variance_RR}, we have $L^2\xi_t^2 \leq \frac{3}{5}$ and $6 L^2 \xi_t^2 (5 -3\beta) (\Theta+n) \leq n(1-\beta)$ for $t\geq 1.$

Similar to the estimate \eqref{estimate_alg_1_t=1} in Theorem \ref{thm_momentum_variance_RR}, we get the following result:
\begin{align*} 
    F( w_0^{(2)} )    & \leq F( w_0^{(1)} ) - \frac{\eta_1}{2}(1-\beta) \norms{ \nabla F( w_0^{(1)} )}^2 
    + \frac{\eta_1}{2} (1-\beta) \frac{1}{n} I_1.
\end{align*}
Taking total expectation and applying Lemma~\ref{lem_momentum_variance_03_RR} we further have:
\begin{align*} 
    \E \left[F( w_0^{(2)} ) \right]   & \leq \E \left[F( w_0^{(1)} )\right] - \frac{\eta_1}{2}(1-\beta)\E\left[ \norms{ \nabla F( w_0^{(1)} )}^2 \right]
    + \frac{\eta_1}{2n} (1-\beta)  \E[I_1]\\
    &\overset{\eqref{bound_I=1_RR}}{\leq} \E \left[F( w_0^{(1)} )\right] - \frac{\eta_1}{2}(1-\beta)\E\left[ \norms{ \nabla F( w_0^{(1)} )}^2 \right]
    + \frac{\eta_1}{2n} (1-\beta)  L^2\xi_1^2  (1-\beta)^2 \cdot \E\left[I_1 + N_1 \right]\\
    &\overset{\eqref{bound_sum_I_N=1_RR}}{\leq} \E \left[F( w_0^{(1)} )\right] - \frac{\eta_1}{2}(1-\beta)\E\left[ \norms{ \nabla F( w_0^{(1)} )}^2 \right]
    + \frac{\eta_1}{2n} (1-\beta)^3  L^2\xi_1^2   \cdot \frac{5 - 3\beta}{2} \E[N_1].
\end{align*}

Noting that $\tilde{w}_{t-1} = w_0^{(t)}, 1-\beta < 1$ and $N_t =(\Theta+n) \big\Vert \nabla  F(w_{0}^{(t)})  \big\Vert ^2 + \sigma^2$, we get
\begin{align*} 
    \E \left[ F(\tilde{w}_{1} ) \right] &\leq \E \left[F(  \tilde{w}_{0} )\right] - \frac{\eta_1}{2}(1-\beta)\E\left[ \norms{ \nabla F(  \tilde{w}_{0} )}^2 \right]
    +\frac{\eta_1  L^2\xi_1^2   (1-\beta) (5 - 3\beta)}{4n}\Big[(\Theta+n) \E \left[\big\Vert \nabla  F(\tilde{w}_{0})  \big\Vert ^2\right] + \sigma^2\Big].
\end{align*}
Rearranging the terms and dividing both sides by $(1-\beta)$, we have
\begin{align*} 
    \frac{\eta_1}{2}\E\left[ \norms{ \nabla F(  \tilde{w}_{0} )}^2 \right]  &\leq \frac{\E \left[F(  \tilde{w}_{0} )- F(\tilde{w}_{1} ) \right]}{(1-\beta)}
    + \frac{\eta_1  L^2\xi_1^2   \cdot (5 - 3\beta)}{4n}(\Theta+n) \E \left[\big\Vert \nabla  F(\tilde{w}_{0})  \big\Vert ^2\right] + \frac{\eta_1 \sigma^2 L^2\xi_1^2 (5 - 3\beta)}{4n}.
\end{align*}
Since $\xi_t$ satisfies $6 L^2 \xi_t^2 (5 -3\beta) (\Theta+n) \leq n(1-\beta)$ as above, we can deduce from the last estimate that
\begin{align*} 
    \frac{\eta_1}{2}\E\left[ \norms{ \nabla F(  \tilde{w}_{0} )}^2 \right]  &\leq \frac{\E \left[F(  \tilde{w}_{0} )- F(\tilde{w}_{1} ) \right]}{(1-\beta)}
    + \frac{(1-\beta)\eta_1 }{4} \E \left[\big\Vert \nabla  F(\tilde{w}_{0})  \big\Vert ^2\right] + \frac{3 \sigma^2 (5 - 3\beta)L^2}{2n(1-\beta)} \cdot \xi_1^3,
\end{align*}

which proves \eqref{bound_F(w)_1_tilde_RR}.
\Eproof

\section{The Proof of Technical Results in Section~\ref{sec:single_shuffling}}\label{apdx:sec:single_shuffling}
This section provides the full proofs of the results in Section~\ref{sec:single_shuffling}.
Before proving Theorem~\ref{thm_momentum}, we introduce some common quantities and provide four technical lemmas.

\subsection{Technical Lemmas}
For the sake of our analysis, we introduce the following quantities:
\begin{align*}
    & h_{i}^{(t)}    := g_{i}^{(t)} - \nabla f ( w_{0}^{(t)} ; \pi ( i + 1 ) ), \hspace{1.42cm} t \geq 1, \ i = 0, \dots, n-1; \tagthis \label{define_h_t} \\
    & k_{i}^{(t-1)} := g_{i}^{(t-1)} - \nabla f ( w_{0}^{(t)} ; \pi ( i + 1 ) ), \qquad t \geq 2, \ i = 0, \dots, n-1; \tagthis \label{define_k_t}\\
    & p_t := \frac{1}{n} \sum_{i=0}^{n-1} m_{i+1}^{(t-1)}, \qquad t \geq 2; \tagthis \label{define_p_t} \\ 
    & q_t := \frac{1-\beta}{n} \sum_{i=0}^{n-1} \left[  \left( \sum_{j=n-i}^{n-1} \beta^{j} \right) g_{i}^{(t-1)} +  \left( \sum_{j=0}^{n-i-1} \beta^{j} \right) g_{i}^{(t)}  \right],  \qquad t \geq 2 . \tagthis \label{define_q_t}
\end{align*}
The proof of Theorem~\ref{thm_momentum} requires the following lemmas as intermediate steps of its proof.

\begin{lem}\label{lem_bounded_momentum}
Let $\{w_i^{(t)}\}$ be generated by Algorithm \ref{sgd_momentum_shuffling2} with $0\leq \beta < 1$ and $\eta_i^{(t)} := \frac{\eta_t}{n}$ for every $t \geq 1$. Suppose that Assumption~\ref{as:A2} holds. Then we have
\begin{align}
    (\textbf{a}) \quad &\norms{ m_{i+1}^{(t)} }^2 \leq G^2,  \quad \text{ for } t \geq 1 \text{ and } i=0,\dots,n-1;
    \label{eq_bounded_momentum}\\
    (\textbf{b}) \quad &\big\Vert p_t \big\Vert^2 \leq G^2, \text{ for } t \geq 2 \quad  \text{and} \quad \big\Vert\nabla F( w_0^{(t)} ) \big\Vert^2 \leq G^2, \text{ for } t \geq 1 . \label{eq_bounded_p_t_and_F}
\end{align}

\end{lem}

\begin{lem} \label{lem_bounded_02}
Under the same setting as of Lemma~\ref{lem_bounded_momentum} and Assumption~\ref{as:A1}(b), if $\xi_t$ is defined by \eqref{define_xi}, then 
for $i = 0, 1, \dots ,n-1$ and $t\geq 2$, it holds that
\begin{align*}
    \max \Big( \big \Vert h_{i}^{(t)} \big \Vert^2, \big \Vert k_{i}^{(t-1)} \big \Vert^2 \Big) \leq L^2 G^2 \xi_t^2. \tagthis \label{eq_bounded_02}
\end{align*}
\end{lem}

\begin{lem} \label{lem_momentum_formula}
Under the same setting as of Lemma~\ref{lem_bounded_momentum}, for $i = 0, 1, \dots, n-1$ and $t\geq 2$, we have 
\begin{align*}
    m_{i+1}^{(t)} 
    &= \beta^{n} m_{i+1}^{(t-1)} + (1-\beta) \Big(\beta^{n-1} g_{i+1}^{(t-1)} + \dots + \beta^{i+1} g_{n-1}^{(t-1)} + \beta^i g_{0}^{(t)} + \dots +  \beta g_{i-1}^{(t)} + g_{i}^{(t)} \Big).
\end{align*}
From this expression we obtain the following sum: 
\begin{align*} 
    \frac{1}{n} \sum_{i=0}^{n-1} m_{i+1}^{(t)} = \beta^n p_t +  q_t, \text{ for } t \geq 2; \tagthis \label{eq_formula_sum}
\end{align*}
where $p_t$ and $q_t$ are defined in \eqref{define_p_t} and \eqref{define_q_t}, respectively. 
\end{lem}

\begin{lem} \label{lem_bounded_first_epoch}
Under the same setting as in Theorem \ref{thm_momentum}, the initial objective value $F(\tilde{w}_{1})$ is upper bounded by
\begin{align*} 
    F( \tilde{w}_{1} )  & \leq F( \tilde{w}_{0} ) + \frac{1}{2 L} \norms{ \nabla F( \tilde{w}_{0} ) }^2 + L \eta_1^2 G^2. \tagthis \label{bound_F(w)_1_tilde_2}
\end{align*}
\end{lem}

\subsection{The Proof of Theorem \ref{thm_momentum}: Key Bound for Algorithm~\ref{sgd_momentum_shuffling2}}

From the update $w_{i+1}^{(t)} := w_{i}^{(t)} - \eta_i^{(t)} m_{i+1}^{(t)}$ at Step \ref{alg:A2_step4} and Step \ref{alg:A2_step5} of Algorithm~\ref{sgd_momentum_shuffling2} with $\eta_i^{(t)} := \frac{\eta_t}{n}$, for $i =  1, ... ,n-1$, we have
\begin{align*}
    w_{i}^{(t)} &= w_{i-1}^{(t)} - \frac{\eta_t}{n} m_{i}^{(t)} = w_{0}^{(t)} - \frac{\eta_t}{n} \sum_{j=0}^{i-1} m_{j+1}^{(t)}.
    \tagthis \label{update_momentum_epoch}
\end{align*}
Now, letting $i=n$ in the estimate and noting that $w_{0}^{(t+1)} = w_{n}^{(t)}$ for all $t\geq 1$, we obtain
\begin{align*}
    w_{0}^{(t+1)} -  w_{0}^{(t)} = w_{n}^{(t)} - w_{0}^{(t)} = - \frac{\eta_t}{n} \sum_{j=0}^{n-1} m_{j+1}^{(t)} 
    \tagthis \label{update_momentum_epoch_alg2}
\end{align*}

From this update and the $L$-smoothness of $F$ from Assumption \ref{as:A1}(c), for $t\geq 2$, we can derive 
\begin{align*}
    F( w_0^{(t+1)} )  & \overset{\eqref{eq:Lsmooth}}{\leq} F( w_0^{(t)} ) + \nabla F( w_0^{(t)} )^{\top}(w_0^{(t+1)} - w_0^{(t)}) + \frac{L}{2}\norms{w_0^{(t+1)} - w_0^{(t)}}^2  \\
    &\overset{\tiny\eqref{update_momentum_epoch_alg2}}{=} F( w_0^{(t)} ) - \eta_t \nabla F( w_0^{(t)} )^\top \bigg( \frac{1}{n} \sum_{j=0}^{n-1} m_{j+1}^{(t)} \bigg) + \frac{L \eta_t^{2}}{2} \Big\Vert \frac{1}{n}\sum_{j=0}^{n-1} m_{j+1}^{(t)} \Big\Vert^2  \\
    &\overset{\tiny\eqref{eq_formula_sum}}{=} F( w_0^{(t)} ) - \eta_t \nabla F( w_0^{(t)} )^\top \left( \beta^n p_{t} +  q_{t} \right) + \frac{L \eta_t^{2}}{2} \big\Vert \beta^n p_{t} +  q_{t} \big\Vert^2  \\
    &\overset{(a)}{\leq} F( w_0^{(t)} ) - \eta_t  \beta^n \nabla F( w_0^{(t)} )^\top  p_{t}  - \eta_t \nabla F( w_0^{(t)} )^\top   q_{t}  + \frac{L \eta_t^{2}}{2} \beta^n \norms{p_{t}}^2 +  \frac{L \eta_t^{2}}{2(1-\beta^n)} \norms{q_{t}}^2  \\
    & \overset{(b)}{\leq} F( w_0^{(t)} ) + \eta_t  \beta^n G^2  - \eta_t \nabla F( w_0^{(t)} )^\top   q_{t}  +  \frac{L \eta_t^{2}}{2} \beta^n G^2 +  \frac{L \eta_t^{2}}{2(1-\beta^n)} \norms{q_{t}}^2  \\
    & \overset{(c)}{\leq} F( w_0^{(t)} )  + 2 \eta_t  \beta^n G^2  - \eta_t \nabla F( w_0^{(t)} )^\top   q_{t}   +  \frac{L \eta_t^{2}}{2(1-\beta^n)}  \norms{q_{t}}^2    
\end{align*}
where (\textit{a}) comes from the convexity of squared norm: $\norms{\beta^n p_{t} +  q_{t}}^2 \leq \beta^n \norms{u_{t}}^2 +  (1-\beta^n)  \big\Vert \frac{q_{t}}{1-\beta^n} \big\Vert ^2 $ for $0 \leq \beta^n < 1$, (\textit{b}) is from the inequalities $\| F( w_0^{(t)} ) \| \leq G$ and $\norms{p_{t}} \leq G$  noted in Lemma~\ref{lem_bounded_momentum}, and (\textit{c}) follows from the fact that $L \eta_t \leq 1$.

Next, we focus on processing the term $\nabla F( w_0^{(t)} )^\top   q_{t}$.
We can further upper bound the above expression as
\begin{align*}
    F( w_0^{(t+1)} ) &\leq F( w_0^{(t)} )  + 2 \eta_t  \beta^n G^2  - \eta_t \nabla F( w_0^{(t)} )^\top   q_{t}  + \frac{L \eta_t^{2}}{2(1-\beta^n)}  \norms{q_{t}}^2 \\
    &= F( w_0^{(t)} )  + 2 \eta_t  \beta^n G^2 - \frac{\eta_t}{1 - \beta^n} (1-\beta^n)\nabla F( w_0^{(t)} )^\top   q_{t}   + \frac{L \eta_t^{2}}{2(1-\beta^n)}  \norms{q_{t}}^2 \\ 
    &\overset{(d)}{=} F( w_0^{(t)} )  + 2 \eta_t  \beta^n G^2 + \frac{\eta_t}{2(1 - \beta^n)} \big \Vert (1-\beta^n)\nabla F( w_0^{(t)} ) -   q_{t} \big \Vert^2 \\
    & \quad - \frac{\eta_t}{2(1 - \beta^n)} \big \Vert (1-\beta^n)\nabla F( w_0^{(t)} )  \big \Vert^2 
    - \frac{\eta_t}{2(1 - \beta^n)} \norms{q_{t}}^2
    + \frac{L \eta_t^{2}}{2(1-\beta^n)}  \norms{q_{t}}^2 \\
    &\overset{(e)}{\leq} F( w_0^{(t)} )  + 2 \eta_t  \beta^n G^2 +   \frac{\eta_t}{2(1 - \beta^n)} M_t
    - \frac{\eta_t(1 - \beta^n)}{2} \big \Vert \nabla F( w_0^{(t)} )  \big \Vert^2, \tagthis \label{eq_key_thm_01}
\end{align*}
where $M_t := \Vert (1-\beta^n)\nabla F( w_0^{(t)} ) -   q_{t}  \Vert^2$.
Here, (\textit{d}) follows from the equality $u^{\top}v = \frac{1}{2}\left(\norms{u}^2 + \norms{v}^2 - \norms{u - v}^2\right)$, and  (\textit{e}) comes from the fact that $\eta_t \leq \frac{1}{L}$. 

Note that $1-\beta^n = (1- \beta) \sum_{j=0}^{n-1}\beta^{j}$ and $\nabla F( w_0^{(t)} ) = \frac{1}{n} \sum_{i=0}^{n-1} \nabla f ( w_{0}^{(t)} ; \pi ( i + 1 ) )$, we can rewrite
\begin{align*}
    (1-\beta^n)\nabla F( w_0^{(t)} ) &= \frac{1-\beta}{n} \sum_{i=0}^{n-1} \Bigg[ \Big( \sum_{j=0}^{n-1}\beta^{j} \Big) \nabla f ( w_{0}^{(t)} ; \pi ( i + 1 )) \Bigg]. \tagthis \label{eq_rewrite_F}
\end{align*} 
Recall the definition of $q_{t}$, $h_{i}^{(t)}$, and $k_{i}^{(t-1)}$ from \eqref{define_q_t}, \eqref{define_h_t}, and \eqref{define_k_t}, respectively, we can bound $M_t$ as
\begin{align*}
    M_t &:=  \big\Vert q_t - (1-\beta^n)\nabla F( w_0^{(t)} )  \big\Vert^2 \\
    &\overset{\eqref{define_q_t}, \eqref{eq_rewrite_F}}{=}  \left\Vert \frac{1-\beta}{n} \sum_{i=0}^{n-1} \left[  \bigg( \sum_{j=n-i}^{n-1} \beta^{j} \bigg) g_{i}^{(t-1)} +  \bigg( \sum_{j=0}^{n-i-1} \beta^{j} \bigg) g_{i}^{(t)} - \Big( \sum_{j=0}^{n-1}\beta^{j} \Big) \nabla f ( w_{0}^{(t)} ; \pi ( i + 1 ))  \right]   \right\Vert^2 \\
    &\overset{\eqref{define_h_t}, \eqref{define_h_t}}{=} \left\Vert \frac{1-\beta}{n} \sum_{i=0}^{n-1} \left[  \bigg( \sum_{j=n-i}^{n-1} \beta^{j} \bigg) k_{i}^{(t-1)} +  \bigg( \sum_{j=0}^{n-i-1} \beta^{j} \bigg) h_{i}^{(t)}  \right] \right\Vert ^2\\
    &\leq \frac{(1-\beta)^2}{n} \sum_{i=0}^{n-1} \left\Vert  \bigg( \sum_{j=n-i}^{n-1} \beta^{j} \bigg) k_{i}^{(t-1)} +  \bigg( \sum_{j=0}^{n-i-1} \beta^{j} \bigg) h_{i}^{(t)}  \right\Vert ^2,
\end{align*}
where the last line comes from the Cauchy-Schwarz inequality. We normalize the last squared norms as follows:
\begin{align*}
    M_t &\leq \frac{(1-\beta)^2}{n} \bigg(\sum_{j=0}^{n-1}\beta^{j} \bigg)^2  \sum_{i=0}^{n-1} \left\Vert \frac{\sum_{j=n-i}^{n-1} \beta^{j}}{\sum_{j=0}^{n-1}\beta^{j}}  k_{i}^{(t-1)} + \frac{\sum_{j=0}^{n-i-1} \beta^{j}}{\sum_{j=0}^{n-1}\beta^{j}}   h_{i}^{(t)}  \right\Vert ^2\\
    & \leq \frac{(1-\beta)^2}{n} \bigg(\sum_{j=0}^{n-1}\beta^{j} \bigg)^2 \sum_{i=0}^{n-1}   \left [ \frac{\sum_{j=n-i}^{n-1} \beta^{j}}{\sum_{j=0}^{n-1}\beta^{j}}   \norms{k_{i}^{(t-1)}} ^2
    + \frac{\sum_{j=0}^{n-i-1} \beta^{j}}{\sum_{j=0}^{n-1}\beta^{j}}  \norms{h_{i}^{(t)}}^2 \right] \qquad \text{ by convexity of } \norms{\cdot}^2\\
    &\leq \frac{(1-\beta)^2}{n} \bigg(\sum_{j=0}^{n-1}\beta^{j} \bigg)^2 \sum_{i=0}^{n-1} \max \Big( \norms{k_{i}^{(t-1)}} ^2; \norms{h_{i}^{(t)}}^2 \Big) \\
    &= \frac{(1-\beta^n)^2}{n}  \sum_{i=0}^{n-1} \max \Big( \norms{k_{i}^{(t-1)}} ^2; \norms{h_{i}^{(t)}}^2 \Big) \hspace{4cm} \text{ since } (1- \beta) \sum_{j=0}^{n-1}\beta^{j} = 1-\beta^n \\
    &\leq  (1-\beta^n)^2 L^2 G^2 \xi_t^2,
\end{align*}
where the last estimate follows from Lemma~\ref{lem_bounded_02} with $\xi_t^2 = \max (\eta_t^2; \eta_{t-1}^2)$ for $t\geq 2$.
Substituting the last estimate into \eqref{eq_key_thm_01} for $t\geq 2$, we obtain
\begin{align*}
    F( w_0^{(t+1)} ) 
    \leq F( w_0^{(t)} )  + 2 \eta_t  \beta^n G^2 +   \frac{\eta_t (1 - \beta^n)}{2}  L^2 G^2 \xi_t^2  - \frac{\eta_t(1 - \beta^n)}{2} \big \Vert \nabla F( w_0^{(t)} )  \big \Vert^2.  
\end{align*}
Since $\tilde{w}_{t-1} = w_0^{(t)}$ and $\tilde{w}_{t} = w_0^{(t+1)}$, this inequality leads to  
\begin{align*}
     \eta_t \big \Vert \nabla F( \tilde{w}_{t-1} )  \big \Vert^2 \leq  \frac{2 [F( \tilde{w}_{t-1} ) - F( \tilde{w}_{t})]} {(1 - \beta^n)}   + \frac{4\beta^n G^2}{1 - \beta^n} \eta_t  +   L^2 G^2 \xi_t^3, \ \  t \geq 2.
\end{align*}
Now, using the fact that $\xi_1 = \eta_1$, we can derive the following statement for $t=1$: 
\begin{align*}
     \eta_1 \big \Vert \nabla F( \tilde{w}_{0} )  \big \Vert^2 \leq  \frac{\eta_1 \norms{ \nabla F( \tilde{w}_{0} ) }^2 }{1 - \beta^n}   + \frac{4\beta^n G^2}{1 - \beta^n} \eta_1  +   L^2 G^2 \xi_1^3.
\end{align*} 
Summing the previous statement for $t = 2, 3, \dots, T$, and using the last one, we can deduce that 
\begin{align*}
     \sum_{t=1}^T \left(\eta_t \big \Vert \nabla F( \tilde{w}_{t-1} )  \big \Vert^2  \right) & \leq   \frac{2 [F( \tilde{w}_{1} ) - F_* ] + \eta_1 \norms{ \nabla F( \tilde{w}_{0} ) }^2 } {(1 - \beta^n)}   + \frac{4\beta^n G^2}{1 - \beta^n}  \sum_{t=1}^T \eta_t  +   L^2 G^2  \sum_{t=1}^T \xi_t^3 \\
     & \overset{\eqref{bound_F(w)_1_tilde_2}}{\leq}  \frac{2 [F( \tilde{w}_{0} ) - F_* ] + \left( \frac{1}{L} + \eta_1 \right) \norms{ \nabla F( \tilde{w}_{0} ) }^2 + 2 L \eta_1^2 G^2 } {(1 - \beta^n)}   + \frac{4\beta^n G^2}{1 - \beta^n}  \sum_{t=1}^T \eta_t  +   L^2 G^2  \sum_{t=1}^T \xi_t^3. 
\end{align*}
Finally, dividing both sides of this estimate by $\sum_{t=1}^T \eta_t$, we obtain \eqref{eq_thm_momentum}.
\Eproof

\subsection{The Proof of Corollaries~\ref{co:SSM_constant_LR} and \ref{co:SSM_diminishing_LR}: Constant and Diminishing Learning Rates}\label{apdx:subsec:main_result_SSM}

\begin{proof}[The proof of Corollary~\ref{co:SSM_constant_LR}]
First, from $\beta = \left(\frac{\nu}{T^{2/3}}\right)^{1/n}$ we have $\beta^n = \frac{\nu}{T^{2/3}}$. Since $\beta \leq \left(\frac{R-1}{R}\right)^{1/n}$ for some $R \geq 1$, we have $\frac{1}{1 - \beta^n} \leq R$. We also have $\eta_t := \frac{\gamma}{T^{1/3}}$ and $\eta_t \leq \frac{1}{L}$ for all $t \geq 1$.
Moreover, $\sum_{t=1}^T\eta_t = \gamma T^{2/3}$ and $\sum_{t=1}^T\xi_t^3 = \gamma^3$.
Substituting these expressions into \eqref{eq_thm_momentum}, we obtain
\begin{align*} 
\mathbb{E}\big[\| \nabla F( \hat{w}_T )  \|^2 \big] & \leq  \displaystyle\frac{\Delta_1}{\gamma (1 - \beta^n)T^{2/3}} + \frac{L^2 G^2 \gamma^2}{T^{2/3}} + \frac{4\beta^n G^2}{1 - \beta^n} \\
&  \leq \frac{\Delta_1 R}{\gamma T^{2/3}} + \frac{L^2 G^2 \gamma^2}{T^{2/3}} + 4\beta^n G^2 R \\ 
&=  \displaystyle\frac{\Delta_1 R}{\gamma T^{2/3}} + \frac{L^2 G^2 \gamma^2}{T^{2/3}} + \frac{4 \nu G^2 R}{T^{2/3}} \\
& \leq \frac{\frac{R}{\gamma} \Delta_1 + L^2 G^2 \gamma^2 + 4\nu  G^2 R }{T^{2/3}}. 
\end{align*}
Here, we have $\Delta_1 = 2 [F( \tilde{w}_{0} ) - F_* ] + \left( \frac{1}{L} + \frac{\gamma}{T^{1/3}} \right) \norms{ \nabla F( \tilde{w}_{0} ) }^2 + \frac{2 L G^2 \gamma^2}{T^{2/3}}$ as defined in Theorem~\ref{thm_momentum}. 
Substituting this expression into the last inequality, we obtain the main inequality in Corollary \ref{co:SSM_constant_LR}.
\end{proof}

\begin{proof}[The proof of Corollary~\ref{co:SSM_diminishing_LR}]
We have $\eta_t = \frac{\gamma}{(t + \lambda)^{1/3}}$ and $\eta_t \leq \frac{1}{L}$ for all $t \geq 1$.
Moreover, we have

\begin{equation*}
\begin{array}{lll}
& \sum_{t=1}^T\eta_t  & = \gamma \sum_{t=1}^T\frac{1}{(t+\lambda)^{1/3}} \geq \gamma \int_1^T\frac{d\tau}{(\tau + \lambda)^{1/3}} \geq \gamma \left[ (T + \lambda)^{2/3} - (1 + \lambda)^{2/3} \right],  \\
& \sum_{t=1}^T\xi_{t}^3 & = 2\eta_1^3 + \gamma^3 \sum_{t=3}^{T} \frac{1}{t - 1 + \lambda} 
\leq \frac{2\gamma^3}{1+\lambda} + \gamma^3 \int_{t=2}^T\frac{d\tau}{\tau - 1 + \lambda} \leq \gamma^3 \left[ \frac{2}{1+\lambda} + \log (T - 1 +\lambda) \right].
\end{array}
\end{equation*}
Substituting these expressions with $\beta = \left(\frac{\nu}{T^{2/3}}\right)^{1/n}$ into \eqref{eq_thm_momentum}, we obtain
\begin{align*} 
\mathbb{E}\big[\| \nabla F( \hat{w}_T )  \|^2 \big] & \leq  \displaystyle\frac{\Delta_1}{\gamma \left[ (T+\lambda)^{2/3} - (1+\lambda)^{2/3}\right](1 - \beta^n)} + \frac{L^2 G^2 \gamma^2 \left( \frac{2}{1+\lambda} + \log (T - 1 +\lambda)  \right)}{(T+\lambda)^{2/3} - (1+\lambda)^{2/3}} + \frac{4\beta^n G^2}{1 - \beta^n} \\ 
& \leq \frac{\Delta_1 R}{\gamma \left[ (T+\lambda)^{2/3} - (1+\lambda)^{2/3}\right]} + \frac{L^2 G^2 \gamma^2 \left( \frac{2}{1+\lambda} + \log (T - 1 +\lambda)  \right)}{(T+\lambda)^{2/3} - (1+\lambda)^{2/3}} + 4\beta^n G^2 R \\
& = \frac{\frac{R}{\gamma} \Delta_1}{ (T+\lambda)^{2/3} - (1+\lambda)^{2/3}} + \frac{L^2 G^2 \gamma^2 \left( \frac{2}{1+\lambda} + \log (T - 1 +\lambda)  \right)}{(T+\lambda)^{2/3} - (1+\lambda)^{2/3}} + \frac{4\nu G^2 R}{T^{2/3}} \\
%
&= \frac{\frac{R}{\gamma} \Delta_1 + \frac{2}{1+\lambda} L^2 G^2 \gamma^2 }{(T+\lambda)^{2/3} - (1+\lambda)^{2/3}} + \frac{L^2 G^2 \gamma^2 \log (T-1+\lambda)}{(T+\lambda)^{2/3} - (1+\lambda)^{2/3}} + \frac{4\nu G^2 R}{T^{2/3}}.
\end{align*}
Substituting $\Delta_1$ from Theorem~\ref{thm_momentum} and $\eta_1 = \frac{\gamma}{(1+\lambda)^{1/3}}$ into the last estimate, we obtain the inequality in Corollary~\ref{co:SSM_diminishing_LR}.
\end{proof}

\subsection{The Proof of Technical Lemmas}\label{appdx:subec:Tech_lemmas_proof_2}
We now provide the full proof of lemmas that serve for the proof of Theorem~\ref{thm_momentum} above.

\subsubsection{The Proof of Lemma~\ref{lem_bounded_momentum}}
We prove the first estimate of Lemma~\ref{lem_bounded_momentum} by induction. 
\begin{compactitem}
\item For $t = 1$, it is obvious that $\norms{ m_{0}^{(1)} }^2 = \norms{ \tilde{m}_{0}}^2 = 0 \leq G^2$. 
\item For $t = 1$, assume that \eqref{eq_bounded_momentum} holds for $i-1$, that is, $\norms{ m_{i}^{(1)} }^2 \leq G^2$.   
From the update $m_{i+1}^{(1)} = \beta m_{i}^{(1)} + (1-\beta) g_{k}^{(1)}$, by our induction hypothesis, convexity of $\norms{\cdot}^2$, $0 \leq \beta < 1$, and Assumption~\ref{as:A2}, we have
\begin{equation*}
    \big\Vert m_{i+1}^{(1)} \big\Vert^2  \leq  \beta \big\Vert m_{i}^{(1)} \big\Vert^2 + (1-\beta) \big\Vert g_{i}^{(1)} \big\Vert^2  \leq  \beta G^2 + (1-\beta) G^2 = G^2. 
\end{equation*}
Consequently, for $t=1$, we have $\norms{ m_{i}^{(1)} }^2 \leq G^2$ for all $i = 0, \dots, n - 1$.

\item Now, we prove \eqref{eq_bounded_momentum} for $t > 1$.
Assume that \eqref{eq_bounded_momentum} holds for epoch $t - 1$, that is $\norms{ m_{i+1}^{(t-1)} }^2 \leq G^2$ for every $i = 0, \dots, n - 1$. 
Since $m_0^{(t)} := m^{(t-1}_n$, we have $ \norms{ m_{0}^{(t)} }^2 = \norms{ \tilde{m}_{t - 1}}^2 = \norms{ m_{n}^{(t - 1)} }^2  \leq G^2$. 
Using the previous argument from the case $t=1$, we also get $\norms{ m_{i+1}^{(t)} }^2 \leq G^2$ for all $i = 0, \dots, n - 1$.
\end{compactitem}
Combining all the cases above, we have shown that \eqref{eq_bounded_momentum} holds for all $i = 0, \dots, n - 1$ and $t = 1, \dots, T$.

Now using the Cauchy-Schwarz inequality, we get
\begin{align*}
    \big\Vert p_t \big\Vert^2 = \left\Vert\frac{1}{n} \sum_{i=0}^{n-1} m_{i+1}^{(t-1)} \right\Vert^2 \leq \frac{1}{n} \sum_{i=0}^{n-1} \big\Vert m_{i+1}^{(t-1)} \big\Vert^2 \leq G^2, \text{ for } t \geq 2.
\end{align*}
Similarly, for $t \geq 1$, we also have
\begin{align*}
    \big\Vert\nabla F( w_0^{(t)} ) \big\Vert^2 = \left\Vert \frac{1}{n} \sum_{i=0}^{n-1} \nabla f ( w_{0}^{(t)} ; \pi ( i + 1 ) ) \right\Vert^2 
    \leq \frac{1}{n} \sum_{i=0}^{n-1} \big\Vert  \nabla f ( w_{0}^{(t)} ; \pi ( i + 1 ) ) \big\Vert^2 \leq G^2,
\end{align*}
which proves the second estimate \eqref{eq_bounded_p_t_and_F} of Lemma~\ref{lem_bounded_momentum}.
\Eproof

\subsubsection{The Proof of Lemma~\ref{lem_bounded_02}}
Using the $L$-smoothness assumption of $f(\cdot;i)$, we derive the following for $t \geq 2$:
\begin{align*}
   & \big\Vert h_{i}^{(t)} \big\Vert^2  \overset{\eqref{define_h_t}}{=} \big\Vert\nabla f ( w_{i}^{(t)} ; \pi ( i + 1 ) ) - \nabla f ( w_{0}^{(t)} ; \pi ( i + 1 ) ) \big\Vert^2 \leq L^2 \big \Vert w_{i}^{(t)} - w_{0}^{(t)} \big \Vert^2, \\
    & \big \Vert k_{i}^{(t-1)} \big \Vert^2  \overset{\eqref{define_k_t}}{=} \big \Vert \nabla f ( w_{i}^{(t-1)} ; \pi ( i + 1 ) ) - \nabla f ( w_{0}^{(t)} ; \pi ( i + 1 ) ) \big \Vert^2 \leq L^2 \big \Vert w_{i}^{(t-1)} - w_{0}^{(t)} \big \Vert^2.
\end{align*}
Now from the update $w_{i+1}^{(t)} := w_{i}^{(t)} - \eta_i^{(t)} m_{i+1}^{(t)}$ in Algorithm~\ref{sgd_momentum_shuffling2} with $\eta_i^{(t)} := \frac{\eta_t}{n}$, for $i =  1, \dots, n-1$ and $t \geq 2$, we have
\begin{align*}
    &w_{i}^{(t)} = w_{i-1}^{(t)} - \frac{\eta_t}{n} m_{i}^{(t)} = w_{0}^{(t)} - \frac{\eta_t}{n} \sum_{j=0}^{i-1} m_{j+1}^{(t)} , \\ 
    &w_{0}^{(t)} = w_{n}^{(t-1)} = w_{n-1}^{(t-1)} - \frac{\eta_{t-1}}{n} m_{n}^{(t-1)} = w_{i}^{(t-1)} - \frac{\eta_{t-1}}{n} \sum_{j = i}^{n-1} m_{j+1}^{(t-1)} .  
\end{align*}
Therefore, for $i =  0, \dots ,n-1$ and $t \geq 2$, using these expressions and \eqref{eq_bounded_momentum}, we can bound
\begin{align*}
    & \norms{ w_{i}^{(t)} - w_{0}^{(t)} }^2 = \eta_t^2 \Big\Vert \frac{1}{n} \sum_{j=0}^{i-1} m_{j+1}^{(t)}  \Big\Vert^2  
     \overset{\tiny\eqref{eq_bounded_momentum}}{\leq}  \eta_t^2 G^2,   \\
    & \norms{ w_{i}^{(t-1)} - w_{0}^{(t)} }^2 = \eta_{t-1}^2 \Big\Vert \frac{1}{n} \sum_{j = i}^{n-1} m_{j+1}^{(t-1)} \Big\Vert^2  
     \overset{\tiny\eqref{eq_bounded_momentum}}{\leq}  \eta_{t-1}^2 G^2.   
\end{align*}
Substituting these inequalities into the first two expressions, for $i =  0, \dots ,n-1$ and $t \geq 2$, we get
\begin{align*}
    \max \Big( \big \Vert h_{i}^{(t)} \big \Vert^2, \big \Vert k_{i}^{(t-1)} \big \Vert^2 \Big) &\leq L^2 \max \Big( \eta_t^2 G^2; \eta_{t-1}^2 G^2\Big)  \leq  L^2 G^2 \max \big(\eta_t^2; \eta_{t-1}^2\big)  = L^2 G^2 \max (\eta_t; \eta_{t-1})^2 = L^2 G^2 \xi_t^2,
\end{align*}
which proves \eqref{eq_bounded_02}, where the last line follows from the facts that $\eta_t > 0$ and $\eta_{t-1}>0$.
\Eproof

\subsubsection{The Proof of Lemma~\ref{lem_momentum_formula}}
The first statement follows directly from the momentum update rule, i.e.:
\begin{align*}
    m_{i+1}^{(t)} &= \beta m_{i}^{(t)} + (1-\beta) g_{i}^{(t)} 
    = \beta^2 m_{i-1}^{(t)} + (1-\beta) \beta g_{i-1}^{(t)} + (1-\beta) g_{i}^{(t)} \\ 
    &= \beta^{i+1} m_{0}^{(t)} + (1-\beta) \Big( \beta^i g_{0}^{(t)} + \dots +  \beta g_{i-1}^{(t)} + g_{i}^{(t)} \Big) \\ 
    &\overset{(a)}{=} \beta^{i+1} m_{n}^{(t-1)} + (1-\beta) \Big( \beta^i g_{0}^{(t)} + \dots +  \beta g_{i-1}^{(t)} + g_{i}^{(t)} \Big) \\
    &= \beta^{n} m_{i+1}^{(t-1)} + (1-\beta) \Big(\beta^{n-1} g_{i+1}^{(t-1)} + \dots + \beta^{i+1} g_{n-1}^{(t-1)} + \beta^i g_{0}^{(t)} + \dots +  \beta g_{i-1}^{(t)} + g_{i}^{(t)} \Big), 
\end{align*}
where (\textit{a}) follows from the rule $m_{0}^{(t)} = \tilde{m}_{t-1} =  m_{n}^{(t-1)} $. 

Next, summing up this expression from $i := 0$ to $i := n-1$ and averaging, we get
\begin{align}\label{eq:m_sum}
    \frac{1}{n} \sum_{i=0}^{n-1} m_{i+1}^{(t)} =  \beta^n p_t +  q_t,
\end{align}
where
\begin{align*}
    p_t := \frac{1}{n} \sum_{i=0}^{n-1} m_{i+1}^{(t-1)} \quad \text{and}\quad
    q_t := \frac{1-\beta}{n} \sum_{i=0}^{n-1} \Big(\beta^{n-1} g_{i+1}^{(t-1)} + \dots + \beta^{i+1} g_{n-1}^{(t-1)} + \beta^i g_{0}^{(t)} + \dots  + g_{i}^{(t)}  \Big).
\end{align*}
Now we consider the term $q_t$ in the last expression:
\begin{align*}
    q_t &=   \frac{1 - \beta}{n} \sum_{i=0}^{n-1} \Big(\beta^{n-1} g_{i+1}^{(t-1)} + \dots + \beta^{i+1} g_{n-1}^{(t-1)} + \beta^i g_{0}^{(t)} + \dots +  \beta g_{i-1}^{(t)} + g_{i}^{(t)} \Big)\\
    &=  \ \frac{1-\beta}{n}  \Big(\beta^{n-1} g_{1}^{(t-1)} + \beta^{n-2} g_{2}^{(t-1)}
    \dots + \beta^{1} g_{n-1}^{(t-1)} + g_{0}^{(t)} \Big) + \dots\\
    & \quad + \frac{1-\beta}{n}  \Big(\beta^{n-1} g_{i+1}^{(t-1)} 
    + \beta^{n-2} g_{i+2}^{(t-1)}
    + \dots + \beta^{i+1} g_{n-1}^{(t-1)} + \beta^i g_{0}^{(t)}  \dots +  \beta g_{i-1}^{(t)} + g_{i}^{(t)} \Big) + \dots \\
    & \quad +\frac{1-\beta}{n}  \Big( \beta^{n-1} g_{0}^{(t)} + 
    \beta^{n-2} g_{1}^{(t)}
    \dots +  \beta g_{n-2}^{(t)} + g_{n-1}^{(t)} \Big).
\end{align*}

Reordering the terms $g_{i}^{(t)}$ and $g_{i}^{(t-1)}$ and noting that $\sum_{j=n}^{n-1} \beta^{j} = 0$ by convention, we get
\begin{align*}
    q_t &= \frac{1-\beta}{n} \sum_{i=1}^{n-1} \left[ \bigg( \sum_{j=n-i}^{n-1} \beta^{j} \bigg) g_{i}^{(t-1)} \right]
    + \frac{1-\beta}{n} \sum_{i=0}^{n-1} \left[ \bigg( \sum_{j=0}^{n-i-1} \beta^{j} \bigg) g_{i}^{(t)} \right]\\
    & = \frac{1-\beta}{n} \sum_{i=0}^{n-1} \left[  \bigg( \sum_{j=n-i}^{n-1} \beta^{j} \bigg) g_{i}^{(t-1)} +  \Bigg( \sum_{j=0}^{n-i-1} \beta^{j} \Bigg) g_{i}^{(t)}  \right]. 
\end{align*}
Substituting this expression of $q_t$ into \eqref{eq:m_sum}, we obtain Equation \eqref{eq_formula_sum} of Lemma~\ref{lem_momentum_formula}. 
\Eproof

\subsubsection{The Proof of Lemma~\ref{lem_bounded_first_epoch}: Upper Bounding The Initial Objective Value}
Using \eqref{eq_bounded_momentum}, we have
\begin{align*}
    \norms{w_0^{(2)} - w_0^{(1)}}^2 = \norms{w_n^{(1)} - w_0^{(1)}}^2 = \eta_1^2 \Big \Vert \frac{1}{n} \sum_{i=0}^{n-1} m_{i+1}^{(1)} \Big \Vert^2 \overset{\eqref{eq_bounded_momentum}}{\leq} \eta_1^2 G^2. \tagthis \label{eq_lem_bound_first_01}
\end{align*}
Since $F$ is $L$-smooth, we can derive 
\begin{align*}
    F( w_0^{(2)} )  & \overset{\eqref{eq:Lsmooth}}{\leq} F( w_0^{(1)} ) + \nabla F( w_0^{(1)} )^{\top}(w_0^{(2)} - w_0^{(1)}) + \frac{L}{2}\norms{w_0^{(2)} - w_0^{(1)}}^2 \\
    & \overset{(a)}{\leq} F( w_0^{(1)} ) + \frac{1}{2L} \norms{ \nabla F( w_0^{(1)} ) }^2 + L \norms{w_0^{(2)} - w_0^{(1)}}^2, \\
    & \overset{\eqref{eq_lem_bound_first_01}}{\leq} F( w_0^{(1)} ) + \frac{1}{2L} \norms{ \nabla F( w_0^{(1)} ) }^2 + L \eta_1^2 G^2,
\end{align*}
where $(a)$ follows since $u^{\top}v \leq \frac{\norms{u}^2}{2} + \frac{\norms{v}^2}{2}$. 

Substituting  $w_0^{(2)} := \tilde{w}_{1}$ and  $w_0^{(1)} := \tilde{w}_{0}$ into this estimate, we obtain \eqref{bound_F(w)_1_tilde_2}. 
\Eproof

\section{Detailed Implementation and Additional Experiments}\label{sec:experiments_details}
In this Supp. Doc., we explain the detailed hyper-parameter tuning strategy in Section~\ref{sec:experiments} and provide additional experiments for our proposed methods.

\subsection{Comparing SMG with Other Methods}\label{sec:compare1}
We compare our SMG algorithm with Stochastic Gradient Descent (SGD) and two other methods: SGD with Momentum (SGD-M) \citep{polyak1964some} and Adam \citep{Kingma2014}. 
To have a fair comparison, a random reshuffling strategy is applied to all methods. 
We tune each algorithm with a constant learning rate using grid search and select the hyper-parameters that perform best according to their results. 
We report the additional results on the squared norm of gradients of this experiment in Figure~\ref{fig_appendix_01}. 
The hyper-parameters tuning strategy for each method is given below:
\begin{compactitem}
    \item For SGD, the first (coarse) searching grid is $\{0.1, 0.01, 0.001\}$ and the fine grid is like $\{0.5, 0.4, 0.2, 0.1, 0.08, 0.06, 0.05\}$ if $0.1$ is the best value we get after the first stage.
    \item For our new SMG algorithm, we fixed the parameter $\beta := 0.5$ since it usually performs the best in our experiments. 
    We let the coarse searching grid be $\{1, 0.1, 0.01\}$ and the fine grid be $\{0.5, 0.4, 0.2, 0.1, 0.08, 0.06, 0.05\}$ if $0.1$ is the best value of the first stage. 
    Note that our SMG algorithm may work with a bigger learning rate than the traditional SGD algorithm.
    \item For SGD-M, we update the weights using the following rule: 
    \begin{align*}
        &m_{i+1}^{(t)} := \beta m_{i}^{(t)} + g_i^{(t)}\\
        &w_{i+1}^{(t)} := w_{i}^{(t)} - \eta_{i}^{(t)} m_{i+1}^{(t)},
    \end{align*}
    where $g_i^{(t)}$ is the $(i+1)$-th gradient at epoch $t$. 
    Note that this momentum update is implemented in PyTorch with the default value $\beta = 0.9$. 
    Hence, we choose this setting for SGD-M, and we tune the learning rate using the two searching grids $\{0.1, 0.01, 0.001\}$ and $\{0.5, 0.4, 0.2, 0.1, 0.08, 0.06, 0.05\}$ as in the SGD algorithm. 
    \item For Adam, we fixed two hyper-parameters $\beta_1 := 0.9$, $\beta_2 := 0.999$ as in the original paper. 
    Since the default learning rate for Adam is $0.001$, we let our coarse searching grid be $\{0.01, 0.001, 0.0001\}$, and the fine grid be $\{0.002, 0.001, 0.0005\}$ if $0.001$ performs the best in the first stage. 
    We note that since the best learning rate for Adam is usually $0.001$, its hyper-parameter tuning process requires little effort than other algorithms in our experiments. 
\end{compactitem}

\begin{figure}[ht!] 
\begin{center}
\includegraphics[width=0.245\textwidth]{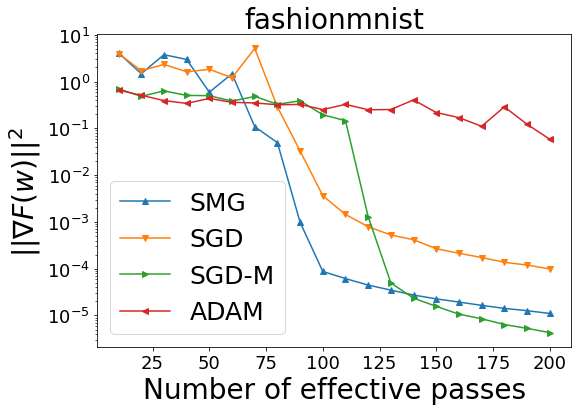}
\includegraphics[width=0.245\textwidth]{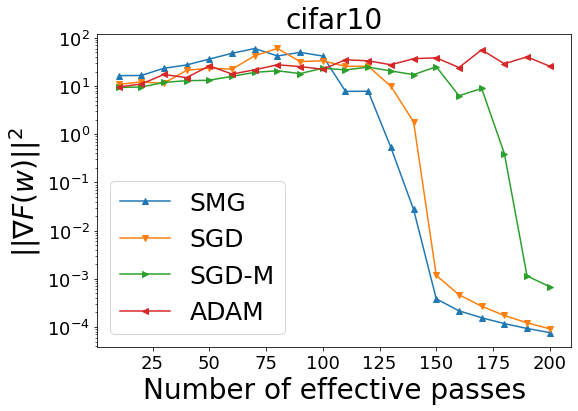}
\includegraphics[width=0.245\textwidth]{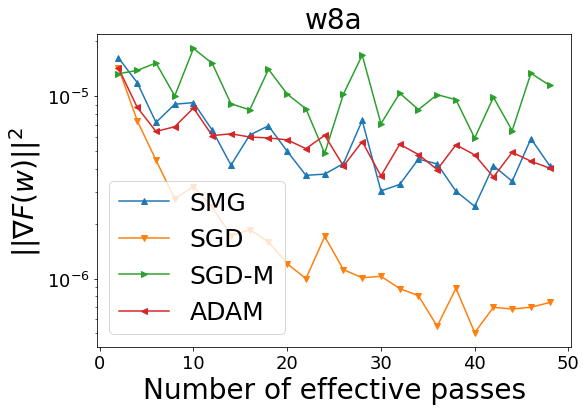}
\includegraphics[width=0.245\textwidth]{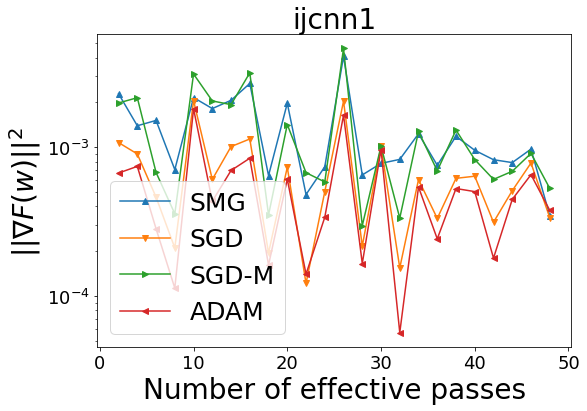}
\caption{The squared norm of gradient produced by SMG, SGD, SGD-M, and Adam for different datasets.}
\label{fig_appendix_01}
\end{center}
\end{figure}

\subsection{The Choice of Hyper-parameter $\beta$}
\begin{figure}[ht!] 
\begin{center}
\includegraphics[width=0.245\textwidth]{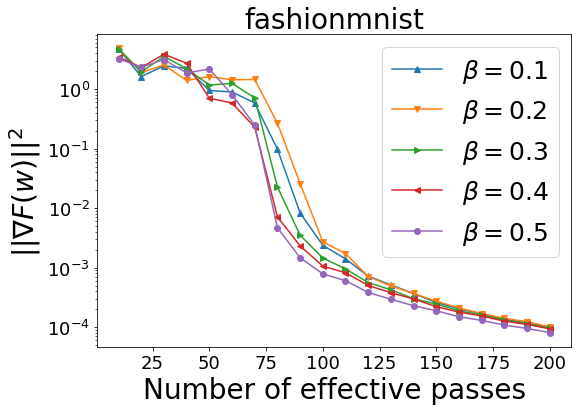}
\includegraphics[width=0.245\textwidth]{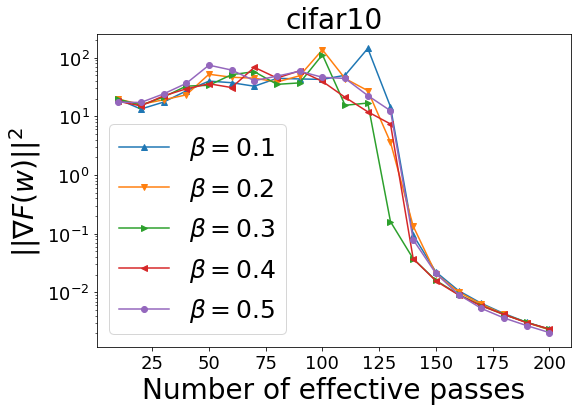}
\includegraphics[width=0.245\textwidth]{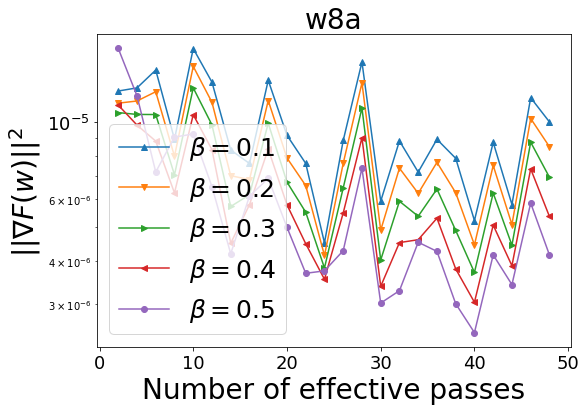}
\includegraphics[width=0.245\textwidth]{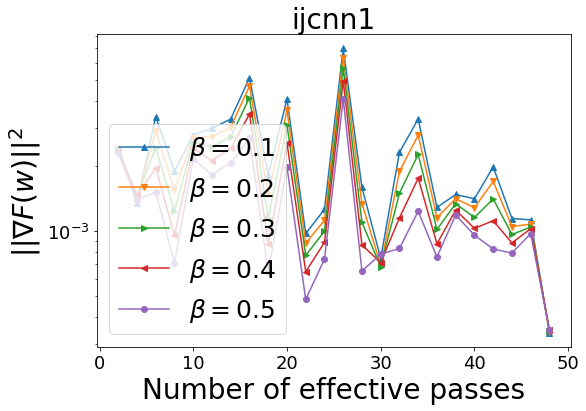}
\caption{The squared norm of gradient reported by SMG with different $\beta$ on \texttt{FashionMNIST}, \texttt{CIFAR-10}, \texttt{w8a}, and \texttt{ijcnn1} datasets}
\label{fig_appendix_02}
\end{center}
\end{figure}

As presented in Section \ref{sec:experiments}, in this experiment, we investigate the sensitivity of our proposed SMG algorithm to the hyper-parameter $\beta$. 
We choose a constant learning rate for each dataset and run the algorithm for different values of $\beta$ in the linear grid $\{0.1, 0.2, 0.3, 0.4, 0.5, 0.6, 0.7, 0.8, 0.9\}$. 
However, the choice of $\beta \geq 0.6 $ does not lead to a good performance, and, therefore, we omit to report them in our results. 
Figure~\ref{fig_appendix_02} shows the comparison of $\Vert \nabla F (w) \Vert^2$ on different values of $\beta$ for \texttt{FashionMNIST}, \texttt{CIFAR-10}, \texttt{w8a}, and \texttt{ijcnn1} datasets. 
Our empirical study here shows that the choice of $\beta$ in $[0.1, 0.5]$ works reasonably well, and $\beta := 0.5$ seems to be the best in our test. 

\subsection{Different Learning Rate Schemes} 



\begin{figure}[ht!] 
\begin{center}
\includegraphics[width=0.245\textwidth]{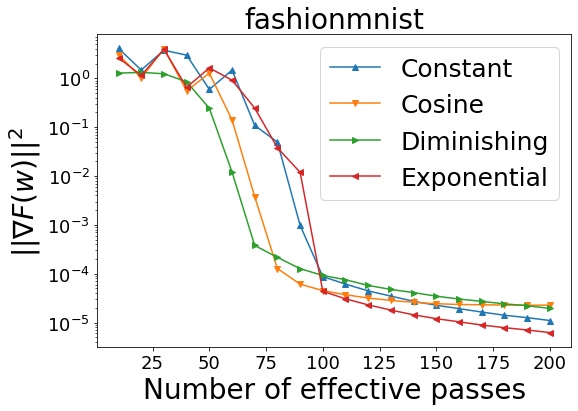}
\includegraphics[width=0.245\textwidth]{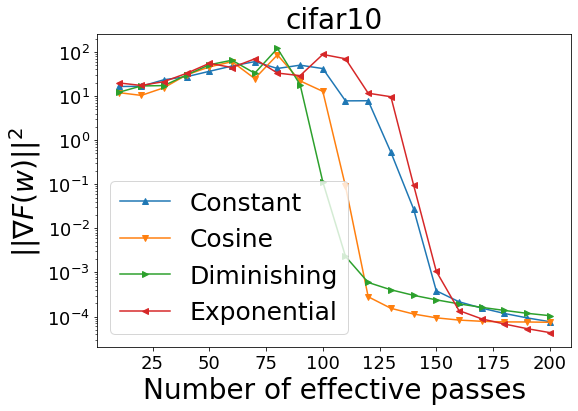}
\includegraphics[width=0.245\textwidth]{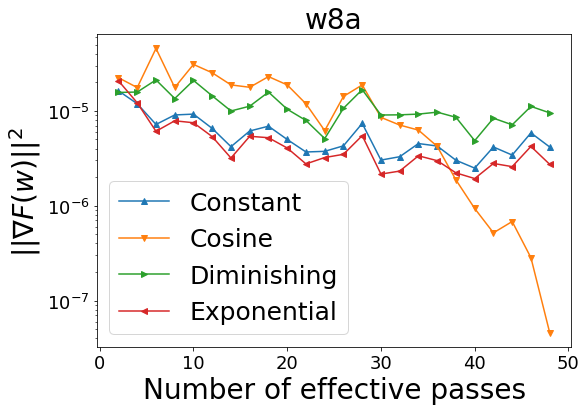}
\includegraphics[width=0.245\textwidth]{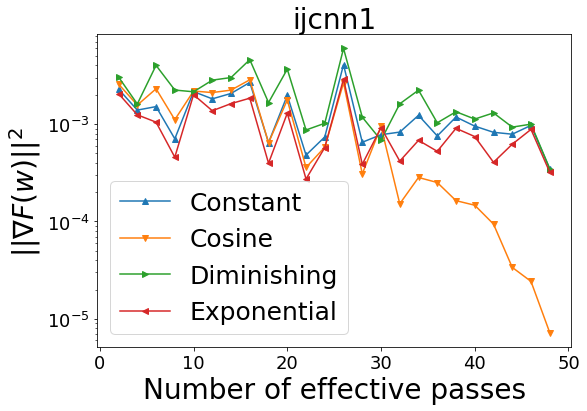}
\caption{The squared norm of gradient produced by SMG under four different learning rate schemes on \texttt{FashionMNIST}, \texttt{CIFAR-10}, \texttt{w8a}, and \texttt{ijcnn1} datasets.}
\label{fig_appendix_04}
\end{center}
\end{figure}

As presented in Section \ref{sec:experiments}, in the last experiment, we examine the performance of our SMG method with the effect of four different learning rate schemes. 
We choose the hyper-parameter $\beta = 0.5$ since it tends to give the best results in our test. 
The additional result on the squared norm gradient is reported in Figure~\ref{fig_appendix_04}.
\begin{compactitem}
    \item \textit{Constant learning rate}: Similar to the first section, we set the coarse searching grid as $\{1, 0.1, 0.01\}$ and let the fine grid be $\{0.5, 0.4, 0.2, 0.1, 0.08, 0.06, 0.05\}$ if $0.1$ is the best value of the first stage.
    
    \item \textit{Cosine annealing learning rate}: For a fixed number of epoch $T$, we need one hyper-parameter $\eta$ for the cosine learning rate scheme: $\eta_t = \eta (1 + \cos (t \pi / T))$. 
    We choose a coarse grid $\{1, 0.1, 0.01, 0.001\}$ and a fine grid $\{0.5, 0.4, 0.2, 0.1, 0.08, 0.06, 0.05\}$ if $0.1$ is the best value for the parameter $\eta$ in the first stage.
    
    \item \textit{Diminishing learning rate}: In order to apply the diminishing scheme $\eta_t = \frac{\gamma}{(t +\lambda)^{1/3}}$, we need two hyper-parameters $\gamma$ and $\lambda$. 
    At first, we let the searching grid for $\lambda$ be $\{1, 2, 4, 8\}$. 
    For the $\gamma$ value, we set its searching grid so that the initial value $\eta_1$ lies in the first grid of $\{1, 0.1, 0.01\}$, and then lies in a fine grid centered at the best one of the coarse grid. 
    
    \item \textit{Exponential learning rate}: We choose two hyper-parameters for the exponential scheme $\eta_t = \eta \alpha^t$. 
    First, let the searching grid for decay rate $\alpha$ be $\{0.99, 0.995, 0.999 \}$. 
    Then we set the searching grid for $\eta$ similarly to the diminishing case (such that initial learning rate $\eta_1$ is in the coarse grid $\{1, 0.1, 0.01\}$, and then a second grid centered at the best value of the first).
\end{compactitem}




\subsection{Additional Experiment with Single Shuffling Scheme}
Algorithm 2 (Single Shuffling Momentum Gradient - SSMG) uses the single shuffling strategy which covers incremental gradient as a special case. Therefore in this experiment we compare our proposed methods (SMG and SSMG) with SGD, SGD with Momentum (SGD-M) and Adam using the single shuffling scheme. 

\begin{figure}[hpt!] 
\includegraphics[width=0.245\textwidth]{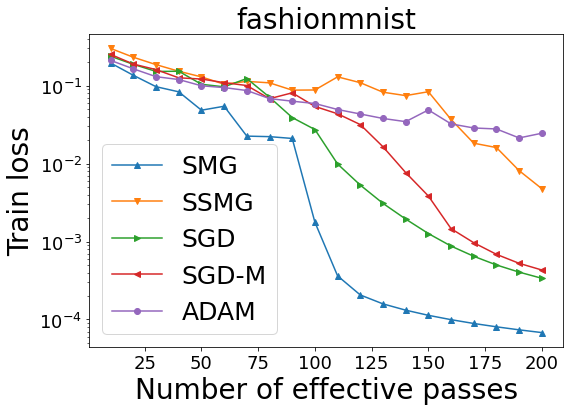}
\includegraphics[width=0.245\textwidth]{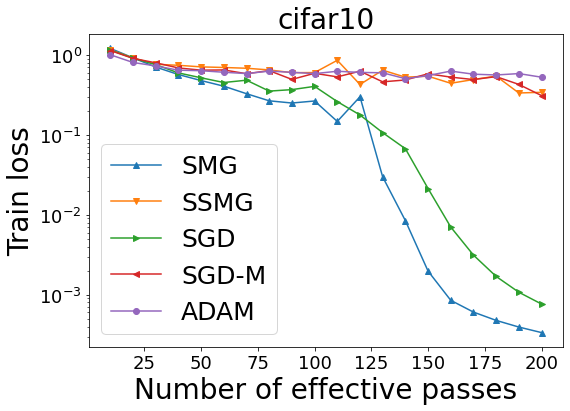}
\includegraphics[width=0.245\textwidth]{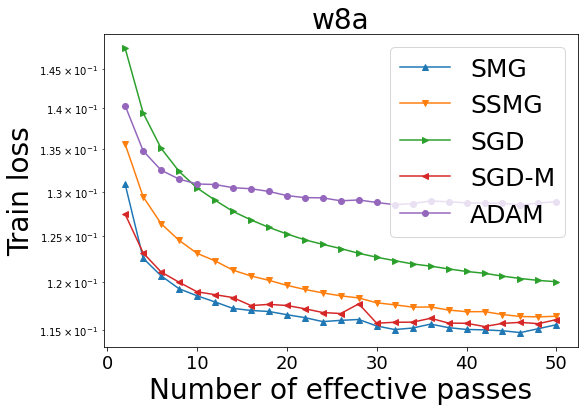}
\includegraphics[width=0.245\textwidth]{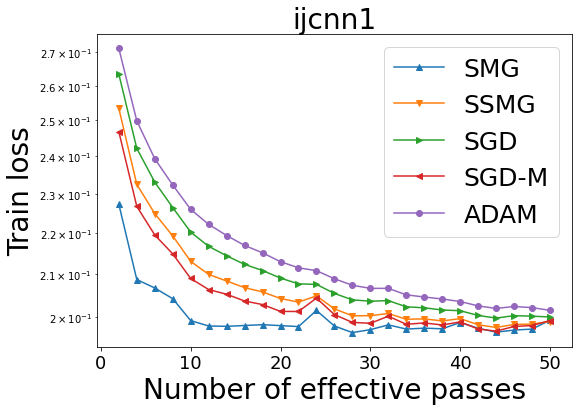}
\caption{The train loss produced by SMG, SSMG, SGD, SGD-M, and Adam on the four different datasets: \texttt{Fashion-MNIST}, \texttt{CIFAR-10}, \texttt{w8a}, and \texttt{ijcnn1}.}
\label{fig_add}
\end{figure}

\begin{figure}[ht!] 
\includegraphics[width=0.245\textwidth]{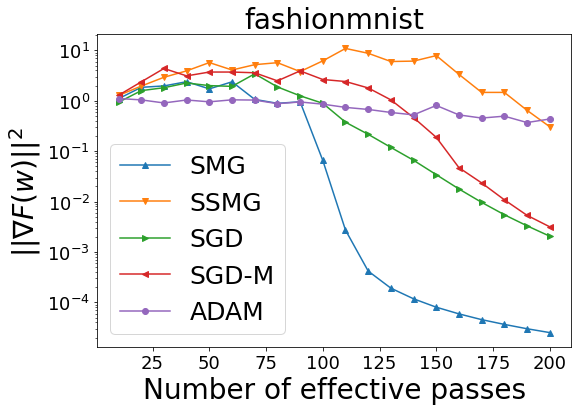}
\includegraphics[width=0.245\textwidth]{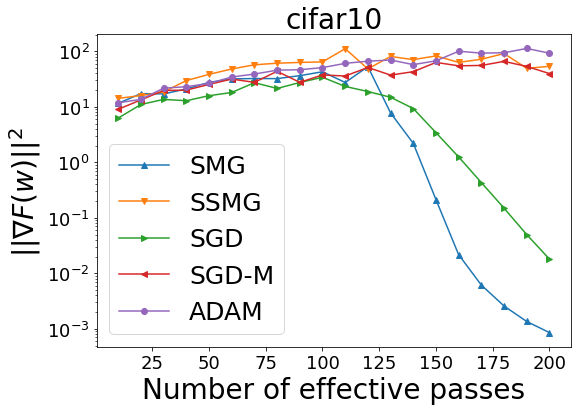}
\includegraphics[width=0.245\textwidth]{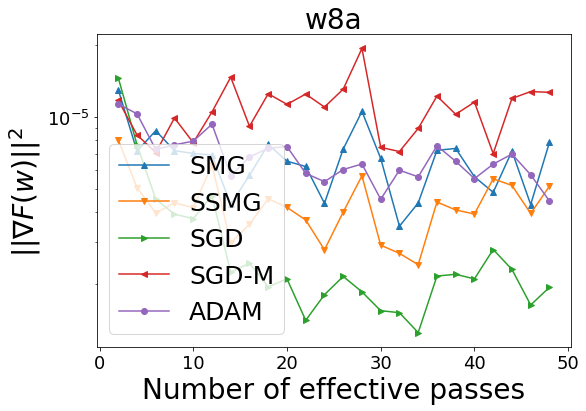}
\includegraphics[width=0.245\textwidth]{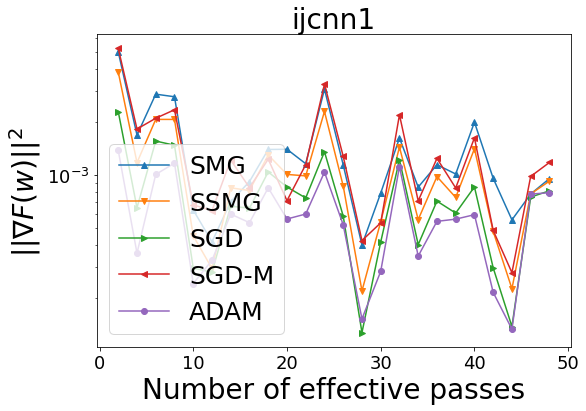}
\caption{The squared norm of gradient produced by SMG, SSMG, SGD, SGD-M, and Adam on the four different datasets: \texttt{Fashion-MNIST}, \texttt{CIFAR-10}, \texttt{w8a}, and \texttt{ijcnn1}.}
\label{fig_add2}
\end{figure}    

For the SSMG algorithm, the hyper-parameter $\beta$ is chosen in the grid $\{0.1, 0.5, 0.9\}$. 
For the learning rate, we let the coarse searching grid be $\{0.1, 0.01, 0.001\}$ and the fine grid be $\{0.5, 0.4, 0.2, 0.1, 0.08, 0.06, 0.05\}$ if $0.1$ is the best value of the first stage. 
For other methods, the hyper-parameters tuning strategy is similar to the settings in subsection \ref{sec:compare1}. 

The train loss $F(w)$ and the squared norm of gradient $\Vert \nabla F(w)\Vert^2$ of each methods  are reported in Figures~\ref{fig_add} and~\ref{fig_add2}. 
We observe that SSMG does not work well in the first two datasets compared to SMG, SGD, and SGD-M, but it performs reasonably well in the last two datasets on the binary classification problem.

\end{document}